\documentclass{article}
\usepackage{amsmath,amssymb,amsthm,cite,needspace,xcolor}
\PassOptionsToPackage{hyphens}{url}
\usepackage{hyperref}
\theoremstyle{plain}
\newtheorem{thm}{Theorem}[section]
\newtheorem{prop}[thm]{Proposition}
\newtheorem{lem}[thm]{Lemma}
\newtheorem{cor}[thm]{Corollary}
\theoremstyle{definition}
\newtheorem{dfn}[thm]{Definition}
\newtheorem{rem}[thm]{Remark}
\newtheorem{ex}[thm]{Example}
\numberwithin{equation}{section}
\newcommand{\C}{\mathbb C}
\newcommand{\RR}{\mathbb R}
\newcommand{\ZZ}{\mathbb Z}
\newcommand{\bX}{\boldsymbol X}
\newcommand{\bY}{\boldsymbol Y}
\newcommand{\bL}{\mathbb L}
\newcommand{\bT}{\mathbb T}
\newcommand{\cO}{\mathcal O}
\newcommand{\cH}{\mathcal H}
\newcommand{\R}{\mathbf R}
\newcommand{\an}{\mathrm{an}}
\newcommand{\hol}{\mathrm{hol}}
\newcommand{\Lotimes}{\mathbin{\otimes^{\mathbf L}}}
\newcommand{\cofib}{\operatorname{cofib}}
\newcommand{\Specan}{\operatorname{Spec}^{\an}}

\newcommand{\im}{\operatorname{im}}
\newcommand{\rk}{\operatorname{rank}}
\newcommand{\vdim}{\operatorname{vdim}}
\newcommand{\DR}{\operatorname{DR}}

\newcommand{\id}{\mathrm{id}}
\newcommand{\Dhat}{\widehat{\DR}^{\an}}
\title{GW/PT-correspondence for CY $4$-fold-II: Virtual Classes}
\author{Ivan Karpov\thanks{Department of Mathematics, MIT}}
\date{2 November 2026}

\makeatletter
\renewcommand{\@maketitle}{%
  \newpage\null
  \vskip 2em%
  \begin{center}%
    {\LARGE\@title\par}%
    \vskip 1.5em%
    {\large\@author\par}%
    \vskip 1em%
    {\large\@date\par}%
    \vskip 1.5em%
    \noindent\hfill
    \begin{minipage}{0.66\textwidth}
      \small\raggedright\itshape
      We believe that there is no essential difficulty in translating everything
      to the complex analytic setting by substituting the theory of Douady spaces
      in place of the theory of Hilbert schemes. Unfortunately, we do not know
      a reference for this.
      \par\smallskip
      {\upshape\raggedleft J.~Pardon,~\cite{PardonCurves}\par}
    \end{minipage}\par
  \end{center}%
  \par\vskip 1.5em%
}

\begin{document}
\maketitle

\begin{abstract}
Let $B$ be a complex manifold and $f:\bX\to B$ a derived \textit{complex-analytic} (not necessarily \textit{algebraic}) space with an oriented relative $(-2)$-shifted symplectic
structure. 

We suppose that the closed symplectic form is locally exact
in the completed relative de Rham complex, and construct a derived
smooth manifold over $B$, independent of additional choices. When $f$ is proper and oriented, this gives
integral virtual fundamental classes on its fibres. 

 Our argument is heavily based on the strategy of Borisov-Joyce who already solved the \textit{algebraic} version of this problem, and uses \textit{derived differential geometry}. The proof also contains a relative
analytic Darboux theorem a-la Brav-Bussi-Joyce, and has a natural analogue in Borel-Moore homology.

We provide applications to complex-analytic enumerative geometry: in particular, to counting Pandharipande-Thomas pairs on compact and non-compact complex-analytic Calabi-Yau $4$-folds.

This is part of a joint project with Felix Thimm dedicated to the proof of the GW/PT-correspondence for Calabi-Yau $4$-folds by implementing Pardon's programme in dimension $4$.
\end{abstract}

\section{Introduction}\label{sec:1}

\subsection{Complex-analytic enumerative geometry}

\subsubsection{}
Enumerative geometry studies the counting of geometric objects, such as, say, curves which lie on Calabi-Yau manifolds.

Surprising as it may sound to a complex geometer, until very recently this part of mathematics was mostly focused on the \textit{algebraic} (or even quasi-projective), rather than the \textit{complex-analytic}, context: very few people really cared about counting curves on non-algebraic complex manifolds. \footnote{See, nevertheless, the important work of Junho Lee~\cite{LeeFamilyGW}
and Bryan--Leung~\cite{BryanLeungK3}, which also treats the non-algebraic
case. See also ~\cite[\S\S~3,7]{PTKKVNotes}.}

One of the reasons was, perhaps, that it was not entirely evident what fundamentally new understanding the study of complex-analytic enumerative geometry might bring, and of what use this understanding might be.

We find that the seminal work \cite{PardonCurves} has provided some answers to these questions: John Pardon essentially uses there the \textit{existence} of well-behaved enumerative theories on non-algebraic complex manifolds to prove a deep enumerative conjecture about the geometry of Calabi-Yau $3$-folds: \textit{the MNOP conjecture} of \cite{MNOPI,MNOPII}, which relates two different ways of curve-counting.

We would like to emphasise that the MNOP conjecture can be formulated entirely in the algebraic realm, so Pardon's use of transcendental methods in \cite{PardonCurves} made his breakthrough even more interesting.

\subsubsection{}
As the reader may learn from the epigraph to the present work, even in Pardon's \textit{$3$-dimensional context}, the foundations needed for complex-analytic enumerative geometry were, until very recently, somewhat shaky.

The situation in \textit{$4$-dimensional} geometry is even worse.
Indeed, curve counting in $4$ dimensions is intricate, and requires methods much more sophisticated than those needed in one dimension lower: as written in the literature, many of the non-trivial foundational arguments do not apply in the \textit{complex-analytic setting} at all (see below).

The present work provides the foundations which make sense of such counting in the non-algebraic situation. It serves as a technical prerequisite to the upcoming joint work of the author with Felix Thimm, which will provide the analogue of Pardon's argument in $4$ complex dimensions, thus proving the conjecture of Maulik, Cao, and Toda \cite{CMT}.

\subsubsection{}
There are currently two ways of providing a set-up for counting suitable objects on a given Calabi-Yau $4$-fold $X$: the original foundations of Borisov-Joyce \cite{BJ}, and the later development due to Oh and Thomas \cite{OT} (which simplified the matter significantly). They are essentially equivalent by another result of Oh and Thomas \cite[Theorem~4]{OTII}, and both of them are written down \textit{only for an algebraic $X$}: moreover, Oh and Thomas assume in their paper that $X$ is quasi-projective.

Since the framework of Oh-Thomas is somewhat more comprehensible,
one would naturally like to choose it for generalisation to our complex-analytic setting. Unfortunately, however, their proofs make significant use of \textit{the resolution property of $X$}: the fact that any perfect complex on a quasi-projective manifold has a resolution by vector bundles. And this fails for complex-analytic manifolds.

We were not able to rewrite the argument of Oh and Thomas without using the resolution property, and so had to resort to a complex-analytic generalisation of the argument in \cite{BJ}, which can be found below.

We have, nevertheless, little doubt that in future a way to get rid of the resolution property in \cite{OT} will be found: maybe by using derived algebraic geometry even more sophsiticatedly than it has been used for these purposes till now?

\subsubsection{}
The main idea of Borisov and Joyce is to note that moduli spaces $M$ of curves which are of interest to us have a natural \textit{derived scheme} structure and, moreover, a \textit{$-2$-symplectic} structure, so that one may use a theorem of \cite[Theorem~5.18]{BBJ} which produces a local description of them in terms of certain explicit commutative differential graded algebras.

One may then invoke the theory of the \textit{category of derived manifolds $\mathbf{dMan}$} due to Joyce \cite{JoyceSummary}: gluing these local descriptions turns $M$ into an object of $\mathbf{dMan}$ and, thus, endows it with a virtual fundamental class, once an orientation is provided. This makes counting curves on $X$ possible.

We essentially follow here their strategy, fighting on the way with inevitable technical difficulties. For example, instead of derived schemes, we are naturally forced into the category of \textit{derived complex-analytic spaces} (due to Porta and Lurie, see \cite[\S12]{LurieIX} and \cite[\S2]{PortaGAGA}), and the objects of this category do \textit{not}, in general, have a natural description in terms of \textit{finitely generated} cdgas.

Indeed, even the algebra of analytic functions on $\mathbb C$ is not finitely generated! One has, thus, to work with \textit{derived analytic rings} which provide an effective replacement for cdgas in the complex-analytic context, and this somewhat muddies the subject.

There are all sorts of other technical difficulties as well: for example, in order to perform Pardon's programme in dimension $4$, the author and Thimm have to use the existence of curve-counting theories \textit{in families}, and the part of the work \cite[Theorem~3.22]{BJ} dedicated to this contained an incorrect argument. In order to overcome this, we have to use methods invented by Park \cite[Remark~5.2.6, Lemma~6.1.2 and Proposition~6.2.1]{Park}.

Moreover, our main curve-counting theory below is
\textit{Pandharipande--Thomas theory}. And even such basic
assertions as the representability of the moduli functor of
compactly supported PT-stable pairs by second countable
Hausdorff complex-analytic spaces have not, to our knowledge,
been written down in the generality required here
(nevertheless, see the important related work in
\cite{BingenerRepresentability, FlennerLubke, Fujiki}).
We therefore give the details below.

Finally, we also provide a \textit{Borel-Moore} analogue of our theory: it will be of much use in the joint project with Thimm when applied to families of \textit{non-compact} Calabi-Yau $4$-folds.

\subsection{Main statements}

\subsubsection{The analytic Darboux theorem}

A derived scheme is locally described by a commutative
differential graded algebra $A$, with $\operatorname{Spec}H^0(A)$
as its underlying classical scheme. A $(-2)$-shifted symplectic
structure is a closed two-form of degree $-2$ whose underlying
form identifies the tangent complex with the cotangent complex
shifted by $-2$: see \cite{PTVV}, and contrast with the usual definition of a classical (not derived) symplectic form.

The classical Darboux theorem states that a symplectic form
can be written in standard form in suitable local coordinates.
Brav, Bussi and Joyce~\cite[Theorem~5.18]{BBJ} proved an
analogue for negatively shifted symplectic derived schemes,
using differential graded algebra. We below prove a relative \textit{complex-analytic
version} of their theorem in degree $-2$.

For this purpose, we use \emph{cell presentations} over a
holomorphic submersion $U\to B$. These are objects that behave like cdgas, and that we use to overcome the difficulty that derived complex-analytic spaces can not, a-priori, be locally described as spectra of \textit{finitely generated} cdgas.
Section~\ref{ss:2.3} associates to each such presentation $A$
a derived complex analytic space $\Specan A$.

\begin{thm}[Theorem~\ref{thm:darboux}]
Let $B$ be a complex manifold and $\bX\to B$ a derived
complex analytic space over $B$.
Suppose that $\bX/B$ has a relative $(-2)$-shifted symplectic
form $\omega$ which is locally exact in the completed relative
de Rham complex, in the sense of Definition~\ref{def:exactness}.

Then every point of $\bX$ has an analytic neighbourhood
equivalent over $B$ to $\Specan A$, where $A$ is a cell
presentation on a holomorphic submersion $U\to B$ with relative
coordinates $x_1,\ldots,x_m$. Its underlying graded algebra
and differential have the form
\[
\begin{gathered}
 A^\sharp=\cO_U[z_1,\ldots,z_r,y_1,\ldots,y_m],
 \qquad |z_a|=-1,\quad |y_i|=-2,\\
 \delta z_a=\tfrac12s_a,\qquad
 \delta y_i=\sum_{a=1}^r
       \frac{\partial s_a}{\partial x_i}z_a,
 \qquad \sum_{a=1}^r s_a^2=0,
\end{gathered}
\]
for holomorphic functions $s_1,\ldots,s_r$ on $U$.
Under this equivalence, $\omega$ is homotopic, as a closed
relative form, to $(\omega_{\mathrm{st}},0,\ldots)$, where
\[
 \omega_{\mathrm{st}}
 =\sum_{i=1}^m d_Bx_i\,d_By_i
  +\sum_{a=1}^r d_Bz_a\,d_Bz_a,
\]
and $d_B$ denotes differentiation relative to $B$.
\end{thm}
\subsubsection{Virtual fundamental classes}\label{derived}

Using the analytic Darboux theorem, we extend the construction
of Borisov and Joyce~\cite{BJ} to derived complex analytic
spaces. In particular, we obtain the following result.

\begin{thm}
Let $\bX$ be a derived complex analytic space, with compact Hausdorff classical truncation $X$.
Suppose that $\bX$ has an oriented $(-2)$-shifted symplectic
structure $(\omega,o)$ and constant complex virtual dimension $d$.
Then there is a canonical virtual fundamental class
\[
 [\bX,\omega,o]^{\mathrm{vir}}\in H_d(X;\ZZ).
\]
\end{thm}

For families, we require the relative symplectic form to be
\emph{locally exact}: near each point it must have a primitive
in the completed relative de Rham complex
(Definition~\ref{def:exactness}).
For a single space this condition is automatic, by
Lemma~\ref{pt:de-rham-detection}(iii).

\begin{thm}\label{thm:main}\label{r:1.1}
Let $B$ be a complex manifold of dimension $b$, and let
$f:\bX\to B$ be a derived complex-analytic space over $B$, with Hausdorff second countable
classical truncation $X$. Suppose that $\bX/B$ has an oriented
locally exact relative $(-2)$-shifted symplectic structure
$(\omega,o)$ and constant relative complex virtual dimension $d$.
Then:
\begin{enumerate}
\renewcommand{\labelenumi}{(\roman{enumi})}

\item Our construction defines an oriented derived smooth
manifold $M\to B_{\RR}$ of virtual dimension $d+2b$, whose
underlying space and map to $B$ are $X\to B$.

\item If $M_0,M_1$ arise from two choices in the construction,
there is an oriented bordism between them over $B_{\RR}$,
with underlying space $X\times[0,1]$.

\item For every $u\in B$ with $X_u$ compact, the derived fibre
$\bX_u=\bX\times_B^{\R}\{u\}$, with the pulled back form
$\omega_u$ and orientation $o_u$, has a canonical class
\[
 [\bX_u,\omega_u,o_u]^{\mathrm{vir}}\in H_d(X_u;\ZZ).
\]

\item If $f$ is proper and $u_0,u_1$ lie in the same connected
component of $B$, then
\[
 (i_{u_0})_*[\bX_{u_0},\omega_{u_0},o_{u_0}]^{\mathrm{vir}}
 =
 (i_{u_1})_*[\bX_{u_1},\omega_{u_1},o_{u_1}]^{\mathrm{vir}}
 \quad\text{in }H_d(X;\ZZ),
\]
where $i_u:X_u\hookrightarrow X$ is the inclusion.
\end{enumerate}
\end{thm}

\subsubsection{Stable pairs}

Our main application concerns Pandharipande--Thomas stable
pairs~\cite{PT,CMT,CaoToda}.
Let $\pi:Y\to B$ be a holomorphic submersion of relative
complex dimension four between
complex manifolds, with a nowhere vanishing relative
holomorphic volume form.
A stable pair on $Y_u$ is a morphism $s:\cO_{Y_u}\to F$,
where $F$ is a nonzero coherent pure one-dimensional sheaf
with compact support and $\operatorname{coker}s$ has support
of dimension at most zero.

In Appendix~\ref{app:pt-moduli}, we prove that the functor
of flat families of stable pairs with proper support is
represented by a Hausdorff second countable complex analytic
space $P(Y/B)$, separated over $B$, with a universal stable
pair. Applying the preceding theorem and using canonical orientations from \cite{KarpovOrientations} and \cite{JU}, we obtain:

\begin{thm}[Theorem~\ref{thm:pt-classes}]
Let $P\subseteq P(Y/B)$ be an open subspace on which
$\chi(F)=n$ is constant. Then $P$ has a derived analytic
enhancement $\mathbf P$ with a canonical orientation and
a locally exact relative $(-2)$-shifted symplectic form.
Every compact fibre $P_u$ has a canonical virtual
fundamental class
\[
 [P_u]^{\mathrm{vir}}\in H_{2n}(P_u;\ZZ).
\]
These classes are compatible with holomorphic base change
between complex manifold bases.
If $P\to B$ is proper, their pushforwards under
$P_u\hookrightarrow P$ agree for points $u$ in the same
connected component of $B$.
\end{thm}

In Section~\ref{sec:noncompact} we extend the construction to
noncompact moduli spaces. We use Pardon's relative homology
$H_*^{\mathrm{rel}\infty}(P/B;\ZZ)$
\cite[Definition~2.2]{PardonCurves}.

\begin{cor}[Corollary~\ref{cor:noncompact-pt}]
With $Y/B$, $P$ and $n$ as above, there is a canonical
relative virtual fundamental class
\[
 [\mathbf P/B]^{\mathrm{vir}}
 \in H_{2n}^{\mathrm{rel}\infty}(P/B;\ZZ)
 \cong H_{2n+2\dim_{\C}B}^{\mathrm{BM}}(P;\ZZ).
\]

This class is compatible with holomorphic base change between
complex manifold bases. Its pullback to $u\in B$ is a class
\[
 [P_u]^{\mathrm{vir}}_{\mathrm{BM}}
 \in H_{2n}^{\mathrm{BM}}(P_u;\ZZ),
\]
which agrees with $[P_u]^{\mathrm{vir}}$ when $P_u$ is compact.
\end{cor}

\subsection{Acknowledgements}

I would like to thank Davesh Maulik for useful encouragement, Dominic Joyce for even more useful criticism, and Felix Thimm for joys of our collaboration.

\subsubsection{AI Disclosure}
The theorems and proofs in this paper are entirely human-generated, and I bear the full responsibility for any mistakes.

However, I used ChatGPT 6 and 6.1 to find several references used in Section \ref{sec:pt-exactness}, Section \ref{sec:noncompact}, and in Appendix \ref{app:pt-moduli}.

Once I finished the draft of this paper, I also used ChatGPT 6.1 for editing it and for proofreading several arguments.

\section{Background material}\label{sec:background}\label{sec:2}
We explain here the parts of derived complex analytic geometry used in this paper. 

We first briefly recall cdgas, their K\"ahler differentials, derived schemes,
and cotangent complexes, following \cite[\S\S2.1--2.2]{BJ}. We then discuss
sheaves of cdgas on complex manifolds, holomorphic presentations, and their
holomorphic differentials. Next we explain the basics of derived analytic
geometry in the sense of Lurie \cite[\S12]{LurieIX}, Porta \cite[\S2]{PortaGAGA},
and Porta-Yu \cite[\S5]{PortaYu}. Finally we discuss closed forms, symplectic
structures, and orientations. 

Throughout this section, $B$ is a complex manifold. 
\subsection{Commutative differential graded algebras and derived schemes}\label{ss:2.1}
We follow \cite[\S\S2.1--2.2]{BJ}, with using notation $\delta$ for the
differential and $A^\sharp$ for the underlying graded algebra of a cdga $A$.

\begin{dfn}\label{r:2.1}
A \emph{commutative differential graded algebra}, or \emph{cdga}, over $\C$ is a
complex of vector spaces
\[
 \cdots\longrightarrow A^{-2}\xrightarrow{\delta}A^{-1}\xrightarrow{\delta}A^0,
\]
with an associative multiplication, bilinear over $\C$, and a unit $1\in A^0$, such that
\[
 ab=(-1)^{ij}ba,\qquad \delta(ab)=(\delta a)b+(-1)^ia\delta b
 \quad\text{for }a\in A^i,\ b\in A^j.
\]
As it was said above, we write $A^\sharp$ for the graded algebra obtained by forgetting $\delta$. \textit{Morphisms}
of cdgas are maps preserving the grading, multiplication, unit, and differential. A morphism is a
\emph{quasiisomorphism} if it induces isomorphisms on cohomology in every degree.
\end{dfn}

Write $\operatorname{cdga}_{\C}$ for this category. A cdga $A$ is of
\emph{standard form} if $A^0$ is a smooth finitely generated $\C$-algebra
of pure dimension and $A^\sharp$ is freely generated over $A^0$ by
finitely many generators in each negative degree.

A differential graded $A$-module is a graded $A^\sharp$-module $M$
with a differential $\delta_M:M^j\to M^{j+1}$ satisfying $\delta_M^2=0$ and
$\delta_M(am)=(\delta a)m+(-1)^{|a|}a\delta_Mm$.
Write $D(\operatorname{mod}A)$ for the derived category of these modules.

For a morphism $R\to A$ of cdgas, the \emph{module of K\"ahler differentials}
$\Omega^1_{A/R}$ is the differential graded $A$-module whose underlying classical module is generated by symbols $\mathrm da$,
for homogeneous $a\in A$, of degree $|a|$, subject to $\C$-linearity,
vanishing on the image of $R$, and the Leibniz relation
\[
 \mathrm d(ab)=a\,\mathrm db+(-1)^{|a||b|}b\,\mathrm da.
\]
Its differential is determined by
$\delta_\Omega(\mathrm da)=\mathrm d(\delta a)$.
Equivalently, $\Omega^1_{A/R}=I/I^2$, where
$I=\ker(A\otimes_R A\to A)$ is the kernel of multiplication and
$\mathrm da=a\otimes1-1\otimes a$ modulo $I^2$.
The tensor product here is the ordinary tensor product of cdgas.
Every derivation $\partial:A\to M$ over $R$ of degree zero which commutes
with differentials factors uniquely through an $A$-module map
$\Omega^1_{A/R}\to M$ commuting with differentials. We write
$\Omega^1_A=\Omega^1_{A/\C}$.

The module $\Omega^1_{A/R}$ need not be invariant under
quasiisomorphisms. The \emph{relative cotangent complex}
corrects this by using a resolution. Choose a quasiisomorphism
$Q\to A$ over $R$, where $Q$ is obtained from $R$ by successively
adjoining free generators in nonpositive degrees, each with differential
in the algebra already constructed. Such a resolution always exists
(we allow infinitely many generators and degrees here). Put
\[
 \bL_{A/R}\simeq\Omega^1_{Q/R}\otimes_Q A.
\]
This ordinary tensor product computes the derived tensor product, and the resulting object of $D(\operatorname{mod}A)$ is independent of the
resolution and invariant under quasiisomorphisms. Write
$\bL_A=\bL_{A/\C}$. 

Localising $\operatorname{cdga}_{\C}$ at quasiisomorphisms results in
the $\infty$-category of derived algebras over $\C$. Its opposite is
the $\infty$-category of \emph{affine derived $\C$-schemes}, with
the object corresponding to $A$ denoted by $\operatorname{Spec}A$.

A \emph{derived $\C$-scheme} $\bY$ is a topological space $|\bY|$
with a sheaf $\cO_{\bY}$ of derived algebras, locally equivalent to
these affine models. Its \emph{classical truncation} $t_0(\bY)$ is the
ordinary scheme with the same underlying space and structure sheaf
$H^0(\cO_{\bY})$, and, in particular,
\[
 t_0(\operatorname{Spec}A)=\operatorname{Spec}H^0(A).
\]
An ordinary scheme is regarded as a derived scheme whose structure
sheaf is concentrated in degree zero. A derived scheme is equivalent
to its classical truncation if and only if its structure sheaf has
zero cohomology in all negative degrees.

The cotangent complex $\bL_{\bY}$ is a quasicoherent complex on $\bY$
whose restriction to each affine open is the cotangent complex of its cdga.
For $\bY=\operatorname{Spec}A$, the
equivalence between quasicoherent complexes and $D(\operatorname{mod}A)$
identifies $\bL_{\bY}$ with $\bL_A$. 

The dual of the cotangent complex dual is the \emph{tangent
complex} $\bT_{\bY}$. For a morphism $f:\bY\to\boldsymbol Z$, the
relative cotangent complex fits into the triangle
\[
 f^*\bL_{\boldsymbol Z}\longrightarrow\bL_{\bY}
 \longrightarrow\bL_{\bY/\boldsymbol Z}
 \longrightarrow f^*\bL_{\boldsymbol Z}[1].
\]
If $\bY$ is locally finitely presented, then $\bL_{\bY}$ is perfect,
and every point has a Zariski open neighbourhood of the form
$\operatorname{Spec}A$ with $A$ of standard form
\cite[Theorem~2.5]{BJ}.

\subsection{Sheaves of cdgas and holomorphic differentials}\label{ss:2.2}
For a complex manifold $U$,
a \emph{sheaf of cdgas on $U$} means a complex of sheaves of vector
spaces over $\C$ on the usual topology of $U$, with multiplication and unit
satisfying Definition~\ref{r:2.1} on each open subset. A morphism of these sheaves is a
quasiisomorphism if it induces isomorphisms of cohomology sheaves.

We will use the following terminology.

\begin{dfn}\label{def:cell}\label{r:2.2}
(a) Let $U$ be a complex manifold and $\pi:U\to B$ a holomorphic
submersion to the fixed complex manifold $B$, that is, we assume that 
$\mathrm d\pi_u:T_uU\to T_{\pi(u)}B$ is surjective for every $u\in U$.

Write $\cO_U$ for the sheaf of holomorphic functions on $U$ (we use the analytic topology on $U$ everywhere here).
A \emph{finite cell presentation}, or simply a \emph{cell presentation}, over $U$ is a sheaf of cdgas $A$
on $U$ with
\[
 A^\sharp=\cO_U[e_1,\ldots,e_N],\qquad |e_i|<0,
\]
freely generated as a graded commutative algebra, and with $\delta|_{\cO_U}=0$. 

\noindent (b) Let $A$ be a finite cell presentation over $U$. A \emph{cell
extension} of $A$ is an inclusion $A\hookrightarrow A'$ constructed
by choosing $r\ge0$ and a sequence
\[
 A=A_0\hookrightarrow A_1\hookrightarrow\cdots\hookrightarrow A_r=A',
\]
so that the following holds.

For each $j=1,\ldots,r$, there exists a negative integer $m_j$ and a section
\[
 b_j\in\Gamma(U,A_{j-1}^{m_j+1}),\qquad \delta b_j=0,
\]
and the graded algebra $A_j^{\sharp}=A_{j-1}^{\sharp}[u_j]$ is freely
generated over $A_{j-1}^{\sharp}$ by one generator $u_j$ of degree
$m_j$, so that its differential extends that of $A_{j-1}$ and is determined
by $\delta u_j=b_j$ and the graded Leibniz rule.
\end{dfn}

Notice that if $Z\subset U$ is defined
by the ideal $\delta(A^{-1})$, then $H^0(A)=\cO_Z$, regarded as a sheaf on $U$, and each
$H^j(A)$ is a coherent $\cO_Z$-module. Indeed, multiplication by $\delta a$, for $a\in A^{-1}$, is
homotopic to zero by multiplication by $a$.

A morphism from a presentation $A$ over $U$ to a presentation $A'$ over $U'$ consists
of a holomorphic map $\varphi:U'\to U$ over $B$ and a morphism $\varphi^{-1}A\to A'$
whose degree zero part is pullback of holomorphic functions. 

\begin{ex}\label{r:2.3}
Given $f_1,\ldots,f_r\in\Gamma(U,\cO_U)$, the Koszul cdga is
\[
 A=\cO_U[z_1,\ldots,z_r],\qquad |z_a|=-1,\qquad\delta z_a=f_a.
\]
Its zeroth cohomology is $\cO_U/(f_1,\ldots,f_r)$. If the functions form a regular sequence,
its negative cohomology vanishes. If $r=1$ and $f_1=0$, then $H^{-1}(A)=\cO_Uz_1$.
\end{ex}

Let $A$ be a finite cell presentation over a holomorphic submersion
$\pi:U\to B$ of complex manifolds. 

For a differential graded
$A$-module $M$, a \emph{holomorphic derivation} $\partial:A\to M$
over $B$ is a derivation of degree zero which commutes with
differentials, vanishes on functions pulled back from $B$, and
satisfies the holomorphic chain rule:  in local relative
coordinates $x_1,\ldots,x_m$,
\[
 \partial h=\sum_{a=1}^m
 \frac{\partial h}{\partial x_a}\,\partial x_a
 \qquad\text{for }h\in\mathcal O_U.
\]
The \emph{module of holomorphic differentials}
$\Omega^1_{\mathrm{hol}}(A/B)$ is characterised by a universal
holomorphic derivation $\mathrm d_{\mathrm{hol}}$:
every holomorphic derivation $A\to M$ over $B$ factors uniquely
through a morphism of differential graded $A$ modules
$\Omega^1_{\mathrm{hol}}(A/B)\to M$.
Writing $e_i$ for the generators of $A$, its underlying
graded module is
\[
 \Omega^1_{\mathrm{hol}}(A/B)^\sharp
 =\bigl(A^\sharp\otimes_{\mathcal O_U}\Omega^1_{U/B}\bigr)
 \oplus\bigoplus_i A^\sharp\mathrm d_{\mathrm{hol}}e_i,
 \qquad |\mathrm d_{\mathrm{hol}}e_i|=|e_i|.
\]
Here $\Omega^1_{U/B}$ is the usual sheaf of relative holomorphic
one-forms. The differential is determined by
$\delta_\Omega(\mathrm d_{\mathrm{hol}}a)
=\mathrm d_{\mathrm{hol}}(\delta a)$.

Note now that the algebraic construction of K\"ahler differentials in
\S\ref{ss:2.1} applies to sheaves by sheafification.
The module $\Omega^1_{\hol}(A/B)$ is then the quotient of
$\Omega^1_{A/\pi^{-1}\cO_B}$ by the differential graded submodule
generated by the holomorphic chain rule relations displayed above.
Here $\pi^{-1}\cO_B$ is the inverse image on $U$ of the sheaf
$\cO_B$. In \S\ref{ss:2.3} we explain how
$\Omega^1_{\hol}(A/B)$ computes the analytic cotangent complex.

\subsection{Derived complex-analytic spaces}\label{ss:2.3}
Having recalled derived schemes in \S\ref{ss:2.1}, we now describe
derived \textit{complex-analytic} spaces. 

The main idea is as follows: for an ordinary
complex analytic space $X$, holomorphic maps to $\mathbb C$ can be substituted
into a holomorphic function whenever their values lie in its domain.

Let $\mathcal T_{\an}$ be the category of complex manifolds and holomorphic maps. We
call local biholomorphisms \emph{admissible morphisms}, and take  surjective
families of them as \textit{coverings}. This category, with its admissible morphisms
and coverings, is Lurie's analytic \emph{pregeometry} \cite[Definition 11.1]{LurieIX}. 

For a topological space $X$, write $\operatorname{Sh}(X)$ for the $\infty$-category of sheaves of
spaces on $X$. 
\begin{dfn}\label{r:2.4}
An \emph{analytic ring} on $X$ is a functor $\mathcal A:\mathcal T_{\an}\to\operatorname{Sh}(X)$ satisfying:
\begin{enumerate}
\renewcommand{\labelenumi}{(\roman{enumi})}
\item $\mathcal A$ preserves finite products.
\item $\mathcal A$ takes a pullback square along an admissible morphism to a homotopy
pullback square.
\item For a covering $\{U_a\to U\}$, the map $\coprod_a\mathcal A(U_a)\to\mathcal A(U)$ is an effective
epimorphism, that is, it is surjective on sheaves of connected components.
\end{enumerate}
A morphism $\mathcal A\to\mathcal A'$ is \emph{local} if the square
\[
 \begin{array}{ccc}
 \mathcal A(V)&\longrightarrow&\mathcal A'(V)\\
 \big\downarrow&&\big\downarrow\\
 \mathcal A(U)&\longrightarrow&\mathcal A'(U)
 \end{array}
\]
is a homotopy pullback for every admissible $V\to U$. Write $\operatorname{AnRing}_{\C}(X)$ for
the $\infty$-category of analytic rings and local morphisms.
\end{dfn}

\begin{dfn}\label{r:2.5}
A \emph{derived complex-analytic space} in this paper is a pair $\bX=(X,\cO_{\bX})$, where $X$ is a topological space and $\cO_{\bX}$ is an analytic ring on $X$, such
that
\[
 t_0(\bX)=\bigl(X,\pi_0(\cO_{\bX}^{\mathrm{alg}})\bigr)
\]
is a classical complex analytic space and every $\pi_j(\cO_{\bX}^{\mathrm{alg}})$ is coherent over its
structure sheaf. A morphism $f:\bX\to\bY$ consists of a continuous map $f:X\to Y$
and a local morphism $f^\#:f^{-1}\cO_{\bY}\to\cO_{\bX}$ of analytic rings. 

Write $\operatorname{dAn}_{\C}$ for
their $\infty$-category.
\end{dfn}

Classical analytic spaces form a full subcategory. 

A derived space is Stein if its truncation is Stein,
and a morphism is proper if its truncation is proper. Derived fibre products
exist, and truncation preserves them \cite[Corollary 12.16]{LurieIX}. We write $\times^{\R}$ or $\times^h$
for a derived fibre product.

The main point for us is that finite cell presentation $A$ over $U$ determines a derived analytic space
$\Specan A$ whose truncation is the analytic subspace $j:Z\hookrightarrow U$ defined by
$\delta(A^{-1})$, with
\[
 \cO^{\mathrm{alg}}_{\Specan A}\simeq j^{-1}A.
\]
For the Koszul cdga of Example~\ref{r:2.3}, this is $U\times^{\R}_{\C^r}\{0\}$, using the map $(f_1,\ldots,f_r)$.
In particular, for the sheaf of holomorphic functions on $U$, $\Specan\cO_U=U$.

For a general cell presentation, the construction uses derived fibre
products. Put $K_0=\C$, pointed at $0$, and $K_{r+1}=\{0\}\times^h_{K_r}\{0\}$. By \cite[Example
12.17 and Lemma 12.18]{LurieIX},
\[
 \operatorname{Map}(\bY,K_r)\simeq
 \operatorname{Map}_{D(\C)}\bigl(\C[r],\R\Gamma(Y,\cO^{\mathrm{alg}}_{\bY})\bigr)
\]
for any space $\bY$.
Here $D(\C)$ denotes the derived $\infty$-category of complexes of
vector spaces over $\C$. 

Starting with $\bY_0=U$, the cocycle $\delta e_i$, for $|e_i|=-n_i$, determines $\bY_{i-1}\to K_{n_i-1}$.
Set $\bY_i=\bY_{i-1}\times^h_{K_{n_i-1}}\{0\}$. The chain of spaces $\bY_i$ gives the desired cell presentation, that is, $\bY_N=:\operatorname{Spec}^{an} A$.

(The same construction can be used to describe morphisms into $\Specan A$. To lift $\boldsymbol W\to\Specan P$
through $\Specan P[e]$, one must choose a homotopy from the pullback
of $\delta e$ to zero. And in a cdga representing derived sections this has a concrete description as an element $h$
with $\delta h=\delta e|_{\boldsymbol W}$.)

We now recall the cotangent complex for $f:\bX\to B$. By \cite[Theorem
4.5]{PortaYu}, modules over $\bX$ may be identified with modules over $\cO^{\mathrm{alg}}_{\bX}$. A module is
\emph{connective} if its cohomology vanishes in positive degrees, and \emph{$n$-connective} if it
vanishes in degrees greater than $-n$. 

A
module is \emph{perfect} if locally it is a retract of a module obtained by finitely many
extensions of finite sums of shifts of $\cO^{\mathrm{alg}}_{\bX}$.

Let $M$ be a connective $\cO^{\mathrm{alg}}_{\bX}$-module.
The analytic \textit{square-zero extension }$\cO_{\bX}\oplus M$ has underlying
algebra $\cO^{\mathrm{alg}}_{\bX}\oplus M$. Its multiplication extends
that of $\cO^{\mathrm{alg}}_{\bX}$ and its action on $M$, with
$M\cdot M=0$ \cite[Corollary~5.16]{PortaYu}.

An \emph{analytic derivation} over $B$ is a section of the projection
$\cO_{\bX}\oplus M\to\cO_{\bX}$ in the $\infty$-category of analytic
rings under $f^{-1}\cO_B$.

We denote the space of these sections by
$\operatorname{Der}^{\an}_B(\cO_{\bX},M)$.
The \emph{analytic cotangent complex} $\bL^{\an}_{\bX/B}$ is the
connective $\cO^{\mathrm{alg}}_{\bX}$-module characterised by a
natural equivalence
\[
 \operatorname{Map}_{\cO^{\mathrm{alg}}_{\bX}}
   (\bL^{\an}_{\bX/B},M)
 \simeq \operatorname{Der}^{\an}_B(\cO_{\bX},M)
\]
for every connective module $M$
\cite[Definition~5.17 and Proposition~5.18]{PortaYu}. 

For $\bX\xrightarrow{f}\bY\to B$ there is a transitivity triangle, and for $\bX'=\bX\times^{\R}_{\bY}\bY'$
with projection $g$ there is a base change equivalence:
\[
 \begin{gathered}
 f^*\bL^{\an}_{\bY/B}\longrightarrow\bL^{\an}_{\bX/B}
 \longrightarrow\bL^{\an}_{\bX/\bY}\longrightarrow f^*\bL^{\an}_{\bY/B}[1],\\
 \bL^{\an}_{\bX'/\bY'}\simeq g^*\bL^{\an}_{\bX/\bY}.
 \end{gathered}
\]
See \cite[Propositions 5.10 and 5.27]{PortaYu}.  

For a holomorphic submersion $U\to B$ one has $\bL^{\an}_{U/B}=\Omega^1_{U/B}$. For a cell
presentation, the holomorphic differentials of \S\ref{ss:2.2} give the model
\begin{equation}\label{e:2.1}
 \bL^{\an}_{\Specan A/B}\simeq j^{-1}\Omega^1_{\hol}(A/B).
\end{equation}

For a perfect module put $M|_x=M\Lotimes\C_x$. Its \emph{Tor-amplitude} lies in $[a,b]$ if
tensoring with any module in degree zero has cohomology only in $[a,b]$, and this condition
can be checked on the fibres $M|_x$. Its restriction to the truncation is locally
a bounded complex of holomorphic vector bundles. In particular, if $\bL^{\an}_{\bX/B}$ is
perfect, its dual is the tangent complex $\bT_{\bX/B}$, and the \emph{complex virtual dimension}
is the locally constant function
\[
 \vdim_{\C}(\bX/B)(x)=\sum_j(-1)^j\dim_{\C}H^j(\bL^{\an}_{\bX/B}|_x).
\]

We also use the following consequence of \cite[Proposition 5.34 and Corollary
5.35]{PortaYu}. For a morphism of analytic rings $\mathcal A\to\mathcal A'$ on the same space, if
$\cofib(\mathcal A^{\mathrm{alg}}\to\mathcal A'^{\mathrm{alg}})$ is $n$-connective for $n\ge1$, then $\bL^{\an}_{\mathcal A'/\mathcal A}$ is $n$-connective. 

The
converse holds when the map on $\pi_0$ is an isomorphism. In particular, a morphism
of derived analytic spaces inducing an isomorphism on truncations and a morphism
with zero relative cotangent complex is an equivalence \cite[Corollary 5.37]{PortaYu}.

\subsection{Closed forms and shifted symplectic structures}\label{ss:2.4}
We describe the derived version of \textit{differential forms} using the local models above.

Let $\pi:U\to B$ be a holomorphic submersion and $A$ a cell
presentation on $U$, as in Definition~\ref{def:cell}. Pullback of
holomorphic functions gives a morphism
$\pi^{-1}\cO_B\to\cO_U=A^0$, so $A$ is a sheaf of
$\pi^{-1}\cO_B$-algebras on $U$. As in \S\ref{ss:2.2},
denote by $\Omega^1_{\hol}(A/B)$  the sheaf of holomorphic differentials
relative to this morphism, and note that its universal derivation vanishes on
$\pi^{-1}\cO_B$. Put also
\[
 \Omega^p_{\hol}(A/B)=\bigwedge^p_A\Omega^1_{\hol}(A/B).
\]

 Write $\mathrm d_B$ for the holomorphic exterior derivative, and note 
it has internal degree zero and commutes with $\delta_\Omega$. We use $\delta=(-1)^p\delta_\Omega$ on $p$-forms, so that $\delta\mathrm d_B+\mathrm d_B\delta=0$.

\Needspace{8\baselineskip}
\begin{dfn}\label{def:forms}\label{r:2.6}
In the preceding notation, suppose that $U$ is Stein and put
$\bY=\Specan A$, with its induced morphism $\bY\to B$.
Let $p\ge0$ and $k\in\ZZ$. A \emph{relative $p$-form of degree
$k$} is a cocycle $\omega_0\in\Gamma(U,\Omega^p_{\hol}(A/B))^k$. A \emph{closed relative $p$-form of degree $k$} is a
sequence $\omega=(\omega_0,\omega_1,\ldots)$ such that
\[
 \omega_i\in\Gamma(U,\Omega^{p+i}_{\hol}(A/B))^{k-i},\qquad
 \delta\omega_0=0,\qquad \delta\omega_{i+1}+\mathrm d_B\omega_i=0.
\]
The form $\omega_0$ is called the \emph{underlying form}.
\end{dfn}

To describe homotopies
between closed forms, assign total degree $k+p$ to a $p$ form of internal degree $k$
and use the complex
\begin{equation}\label{e:2.2}
 C_B(A;U)^j=\prod_{p\ge0}\Gamma(U,\Omega^p_{\hol}(A/B))^{j-p},
 \qquad D=\delta+\mathrm d_B.
\end{equation}

On the complex $C_B(A;U)$ in \eqref{e:2.2}, the multiplication is
\[
 \alpha\cdot\beta=(-1)^{kq}\alpha\wedge\beta
\]
for $\alpha\in\Gamma(U,\Omega^p_{\hol}(A/B))^k$ and
$\beta\in\Gamma(U,\Omega^q_{\hol}(A/B))^l$.`

Let $F^p C_B(A;U)$ consist of components of degree at least $p$. Two closed $p$-forms $\omega,\omega'$ of degree $k$ are homotopic
if and only if $\omega-\omega'=D\eta$ for some
$\eta\in F^p C_B(A;U)^{k+p-1}$.

To define the analytic
derived de Rham complex for a general $\bX\to B$, one may use analytic resolutions and descent. The filtration by form degree on the resulting complex has graded pieces
$\bigwedge^p\bL^{\an}_{\bX/B}[-p]$. Write $\DR^{<N}_{\bX/B}$ for the quotient in form degrees $0,\ldots,N-1$, and
define its completion by
\[
 \Dhat_{\bX/B}=\mathop{\varprojlim}\limits_N^{h}\DR^{<N}_{\bX/B},
 \qquad F^p\Dhat_{\bX/B}=\operatorname{fib}(\Dhat_{\bX/B}\to\DR^{<p}_{\bX/B}).
\]

\begin{rem} The construction roughly goes as follows: on a sufficiently small Stein presentation, derived sections of
this completion are computed by \eqref{e:2.2}.
Forms and closed forms on $\bX$ are then defined by analytic descent from these
local complexes. 
The sets of homotopy classes of relative $p$-forms and closed
relative $p$-forms of degree $k$ are, respectively,
\[
 H^k\R\Gamma(X,\bigwedge^p\bL^{\an}_{\bX/B}),\qquad
 H^{k+p}\R\Gamma(X,F^p\Dhat_{\bX/B}).
\]

\end{rem}

A closed relative $2$-form of degree $k$ is \emph{symplectic} if the underlying cotangent complex
is perfect and contraction with its underlying form defines an equivalence $\bT_{\bX/B}\simeq\bL^{\an}_{\bX/B}[k]$.
For $k=-2$ one obtains
\begin{equation}\label{e:2.3}
 \omega^\sharp:\bT_{\bX/B}\xrightarrow{\sim}\bL^{\an}_{\bX/B}[-2].
\end{equation}
Since the cotangent complex is connective, this implies Tor amplitude in $[-2,0]$.
At $x\in X$ put
\begin{equation}\label{e:2.4}
 T_x=H^0(\bT_{\bX/B}|_x),\qquad O_x=H^1(\bT_{\bX/B}|_x).
\end{equation}
The symplectic structure determines a nondegenerate symmetric form $Q_x$ on $O_x$ and
identifies $H^2(\bT_{\bX/B}|_x)$ with $T_x^\vee$. Therefore
\begin{equation}\label{e:2.5}
 \vdim_{\C}(\bX/B)=2\dim_{\C}T_x-\dim_{\C}O_x.
\end{equation}

\begin{dfn}\label{def:exactness}\label{r:2.7}
A closed relative two form $\omega$ of degree $-2$ is \emph{locally exact} if
its image in $H^0\R\Gamma(X,\Dhat_{\bX/B})$ vanishes on an analytic neighbourhood of each
point. Equivalently, locally $\omega=D\beta$ for $\beta\in C_B(A;U)^{-1}$.
\end{dfn}

\subsection{Orientations}\label{subsec:orientations}\label{ss:2.5}
For a perfect complex $E$ on $X$, represented locally by vector bundles $E^j$, its determinant
line is $\det E=\bigotimes_j(\det E^j)^{(-1)^j}$, where exponent $-1$ means dual. Put
$K_{\bX/B}=\det(\bL^{\an}_{\bX/B}|_X)$. A symplectic structure of degree $-2$ produces a
trivialisation $q_K:K_{\bX/B}^{\otimes2}\xrightarrow{\sim}\cO_X$. We use the convention of \cite[\S2.4]{BJ}: under the
identification $K_{\bX/B}|_x\cong\det O_x$ obtained by cancelling tangent factors, $q_K|_x$ is
$\det Q_x$.

\begin{dfn}\label{def:orientation}\label{r:2.8}
An \emph{orientation} is an isomorphism $o:\cO_X\xrightarrow{\sim}K_{\bX/B}$ satisfying
$q_K\circ o^{\otimes2}=\id$. Its inverse is the orientation isomorphism used in \cite[Definition
2.12]{BJ}.
\end{dfn}
Orientations exist locally, since a nowhere vanishing holomorphic function
locally has a square root. Two orientations differ by a locally constant sign.

\section{Kuranishi spaces and derived manifolds}\label{sec:kuranishi}\label{sec:3}

\subsection{Kuranishi neighbourhoods and coordinate changes}\label{ss:3.1}
\begin{dfn}\label{r:3.1}
Let $X$ be a topological space. A \emph{Kuranishi neighbourhood} on $X$ is a
quadruple $\mathcal U=(V,E,s,\psi)$, where $V$ is a manifold, $E\to V$ is a
vector bundle, $s\in C^\infty(V,E)$, and $\psi:s^{-1}(0)\to X$ is a
homeomorphism onto an open subset $S=\im\psi$ of $X$.
Its \emph{virtual dimension} is $\dim V-\operatorname{rank}E$,
which we require to be constant. For an open $W\subseteq V$, its restriction
to $W$ is $(W,E|_W,s|_W,\psi|_{s^{-1}(0)\cap W})$.
\end{dfn}

Let $\mathcal U=(V,E,s,\psi)$ be a Kuranishi neighbourhood, and let
$v\in s^{-1}(0)$. Choose a smooth local frame $e_1,\ldots,e_r$ of $E$
near $v$, and write $s=\sum_{a=1}^r s^a e_a$, where the $s^a$ are smooth
real functions. We define
\[
 \mathrm ds|_v:T_vV\longrightarrow E_v,\qquad
 \mathrm ds|_v(w)=\sum_{a=1}^r\mathrm d s^a|_v(w)e_a(v).
\]
This map is independent of the frame.

We can also compute this map using any smooth linear connection $\nabla$
on $E$. Such a connection is an $\mathbb R$-linear map
\[
 \nabla:C^\infty(V,E)\longrightarrow C^\infty(V,T^*V\otimes E)
\]
satisfying $\nabla(ha)=\mathrm dh\otimes a+h\nabla a$ for every smooth
function $h$ on $V$ and every smooth section $a$ of $E$. In the frame above,
\[
 \nabla s=\sum_{a=1}^r\mathrm d s^a\otimes e_a
             +\sum_{a=1}^r s^a\nabla e_a.
\]
Thus $(\nabla s)|_v=\mathrm ds|_v$ whenever $s(v)=0$. 

For $x=\psi(v)$, define the \emph{tangent space} and \emph{obstruction
space} of $\mathcal U$ at $x$ by
\begin{equation}\label{e:3.1}
 T_x\mathcal U=\ker(\mathrm ds|_v),\qquad
 O_x\mathcal U=\operatorname{coker}(\mathrm ds|_v).
\end{equation}
These are the cohomology groups in degrees $0,1$ of the complex
$[T_vV\xrightarrow{\mathrm ds|_v}E_v]$. In particular,
$\dim T_x\mathcal U-\dim O_x\mathcal U=\operatorname{vdim}\mathcal U$.

We next define the notation $O(s)$ and $O(s^2)$ used in the definitions
of morphisms, following Joyce \cite[Definition~2.1]{JoyceKur}.

\begin{enumerate}
\renewcommand{\labelenumi}{(\roman{enumi})}
\item Let $F\to V$ be a vector bundle, and let $a,b\in C^\infty(V,F)$.
We write $a=b+O(s)$ if there is a bundle map $A:E\to F$ over $V$ such that
$a-b=A(s)$. We write $a=b+O(s^2)$ if there is a bundle map
$B:E\otimes E\to F$ over $V$ such that $a-b=B(s\otimes s)$.
In local frames, these conditions respectively mean
\[
 a^\mu-b^\mu=\sum_\alpha A^\mu_\alpha s^\alpha,
 \qquad
 a^\mu-b^\mu=\sum_{\alpha,\beta}
                 B^\mu_{\alpha\beta}s^\alpha s^\beta,
\]
with smooth coefficient functions $A^\mu_\alpha$ and
$B^\mu_{\alpha\beta}$. We write $a=O(s)$ or $a=O(s^2)$ when $b=0$.
For scalar functions, we take $F=V\times\mathbb R$.

The condition $a=b+O(s)$ implies $a=b$ on $s^{-1}(0)$, but the converse
need not hold. For example, if $s(t)=t^2$ on $\mathbb R$, the function
$a(t)=t$ vanishes on $s^{-1}(0)$ but cannot be written as $A(t)t^2$ with
$A$ smooth near $0$.

\item Let $f,g:V\to W$ be smooth maps of manifolds. We write
$g=f+O(s)$ if
\[
 h\circ g-h\circ f=O(s)
\]
for every smooth function $h:W\to\mathbb R$. The notation
$g=f+O(s^2)$ is defined in the same way, with $O(s^2)$ in place of
$O(s)$. In particular, either condition implies $g=f$ on $s^{-1}(0)$.

If $\Lambda:E\to f^*TW$ is a bundle map over $V$, we write
$g=f+\Lambda s+O(s^2)$ if, for every smooth function $h:W\to\mathbb R$,
\[
 h\circ g-h\circ f-(f^*\mathrm dh)(\Lambda s)=O(s^2).
\]
Here $(f^*\mathrm dh)(\Lambda s)$ is the smooth function whose value at
$u\in V$ is $\mathrm dh|_{f(u)}\bigl(\Lambda_u(s(u))\bigr)$.
 In local coordinates $y^1,\ldots,y^n$ on $W$, the condition reads
\[
 y^j\circ g-y^j\circ f=(f^*\mathrm dy^j)(\Lambda s)+O(s^2)
 \qquad\text{for }j=1,\ldots,n,
\]
where both maps take values in the coordinate domain.

\item Suppose $g=f+O(s)$, and let $F\to W$ be a vector bundle. We also
need to compare sections $a\in C^\infty(V,f^*F)$ and
$b\in C^\infty(V,g^*F)$. Their values lie in different fibres of $F$
when $f(u)\ne g(u)$, so their difference is not defined as a section.
We define $b=a+O(s)$ as follows.

For each $v\in s^{-1}(0)$, choose an open neighbourhood $W_0$ of
$f(v)=g(v)$ on which $F$ has a smooth frame
$\epsilon_1,\ldots,\epsilon_m$. Choose an open neighbourhood $D$ of $v$
such that $f(D),g(D)\subseteq W_0$, and write
\[
 a|_D=\sum_{\mu=1}^m a^\mu f^*\epsilon_\mu,
 \qquad
 b|_D=\sum_{\mu=1}^m b^\mu g^*\epsilon_\mu.
\]
We require $b^\mu-a^\mu=O(s|_D)$ for every $\mu$, in the sense of
\textup{(i)}. It is easy to see that this condition is independent of the frame.

The same convention applies to bundle maps $H\to f^*F$ and
$H\to g^*F$ with a common source bundle $H\to V$: after choosing a
local frame of $H$, we compare their matrix entries modulo $O(s)$.
For each $a\in C^\infty(V,f^*F)$, there exists
$b\in C^\infty(V,g^*F)$ with $b=a+O(s)$. Any two such sections $b$
differ by $O(s)$ in the sense of \textup{(i)}
\cite[Definition~2.1(vi)]{JoyceKur}.
\end{enumerate}

All these conditions are local on $V$.

\begin{dfn}\label{r:3.2}
Let $\mathcal U_i=(V_i,E_i,s_i,\psi_i)$ and
$\mathcal U_j=(V_j,E_j,s_j,\psi_j)$ be Kuranishi neighbourhoods on $X$, and
let $S\subseteq\im\psi_i\cap\im\psi_j$ be open. A \emph{1-morphism}
$\Phi_{ij}:\mathcal U_i\to\mathcal U_j$ over $S$ consists of an open
neighbourhood $V_{ij}\subseteq V_i$ of $\psi_i^{-1}(S)$, a smooth map
$\phi_{ij}:V_{ij}\to V_j$, and a bundle map
$\chi_{ij}:E_i|_{V_{ij}}\to\phi_{ij}^*E_j$, such that
\begin{equation}\label{e:3.2}
 \chi_{ij}s_i=\phi_{ij}^*s_j+O(s_i^2),\qquad
 \psi_j\circ\phi_{ij}=\psi_i\quad\text{on }s_i^{-1}(0)\cap V_{ij}.
\end{equation}
We write $\Phi_{ij}=(V_{ij},\phi_{ij},\chi_{ij})$. If $\Phi_{ij}$ is defined
over $S_{ij}$ and $\Phi_{jk}$ over $S_{jk}$, their composite over
$S_{ij}\cap S_{jk}$ is represented by
\begin{equation}\label{e:3.3}
 \bigl(V_{ij}\cap\phi_{ij}^{-1}(V_{jk}),\quad
 \phi_{jk}\circ\phi_{ij},\quad(\phi_{ij}^*\chi_{jk})\circ\chi_{ij}\bigr)
\end{equation}
The identity is $(V_i,\id_{V_i},\id_{E_i})$.
\end{dfn}

Differentiating \eqref{e:3.2} at a zero, one get
$\chi_{ij}\mathrm ds_i=\mathrm ds_j\mathrm d\phi_{ij}$. Thus a 1-morphism
induces a map of the complexes defining \eqref{e:3.1}. A
\emph{coordinate change} will be a 1-morphism which is an equivalence, in
the sense explained below (see the next subsection). The following useful criterion is
\cite[Theorem 4.16]{JoyceKur}.

\begin{thm}\label{r:3.3}
A 1-morphism $\Phi_{ij}$ over $S$ is a coordinate change if and only if, for
every $x\in S$, with $v=\psi_i^{-1}(x)$ and $w=\psi_j^{-1}(x)$, the sequence
\begin{equation}\label{e:3.4}
 0\longrightarrow T_vV_i
 \xrightarrow{(\mathrm ds_i,\mathrm d\phi_{ij})}E_{i,v}\oplus T_wV_j
 \xrightarrow{-\chi_{ij}+\mathrm ds_j}E_{j,w}\longrightarrow0
\end{equation}
is exact. Equivalently, its maps on tangent and obstruction spaces are
isomorphisms. In particular, coordinate changes preserve virtual dimension.
\end{thm}

\begin{ex}\label{r:3.4}
For $r\geqslant0$, replace $(V,E,s,\psi)$ by
\[
 \bigl(V\times\RR^r,\operatorname{pr}_V^*E\oplus\RR^r,
 (v,u)\longmapsto(s(v),u),\psi\circ\operatorname{pr}_V\bigr).
\]
Projection to $V$ and to $E$ defines a coordinate change. On the complex at
a zero, its kernel is $[\RR^r\xrightarrow{\id}\RR^r]$, which has zero
cohomology. This is called \emph{stabilisation}.
\end{ex}

\subsection{Morphisms and cocycle conditions}\label{ss:3.2}

\begin{dfn}\label{r:3.5}
Let $\Phi_{ij}$ and $\Phi'_{ij}$ be 1-morphisms over the same open set $S$.
A \emph{2-morphism} $\Lambda:\Phi_{ij}\Rightarrow\Phi'_{ij}$ is represented
by an open neighbourhood $\dot V_{ij}\subseteq V_{ij}\cap V'_{ij}$ of
$\psi_i^{-1}(S)$ and a bundle map
$\lambda:E_i|_{\dot V_{ij}}\to\phi_{ij}^*TV_j$ satisfying
\begin{equation}\label{e:3.5}
 \begin{aligned}
 \phi'_{ij}&=\phi_{ij}+\lambda s_i+O(s_i^2),\\
 \chi'_{ij}&=\chi_{ij}+(\phi_{ij}^*\nabla s_j)\circ\lambda+O(s_i).
 \end{aligned}
\end{equation}
Here $\nabla$ is any connection on $E_j$ and the second equation uses the
preceding convention for comparing pullbacks. Two representatives define
the same 2-morphism if, on a smaller neighbourhood of $\psi_i^{-1}(S)$,
their bundle maps differ by $O(s_i)$. We write $\Lambda=[\dot V_{ij},\lambda]$.
\end{dfn}

This definition is independent of the connection.
 At a zero, \eqref{e:3.5} implies
\[
 \mathrm d\phi'_{ij}-\mathrm d\phi_{ij}=\lambda\mathrm ds_i,
 \qquad\chi'_{ij}-\chi_{ij}=\mathrm ds_j\lambda.
\]
These are the equations for a homotopy of maps of complexes. Consequently,
2-isomorphic 1-morphisms induce the same maps on tangent and obstruction
spaces.

Vertical composition, written $\Lambda'\odot\Lambda$, is represented by
$\lambda+\lambda'$ on a common smaller neighbourhood, after identifying
the pullback bundles modulo $O(s_i)$. The identity is represented by zero,
and the inverse by $-\lambda$. All 2-morphisms are therefore invertible.
Restricting a 1-morphism to a smaller neighbourhood of the same zero locus
gives a canonical 2-isomorphism, also represented by zero.

We write $*$ for horizontal composition. If
$\Lambda:\Phi_{ij}\Rightarrow\Phi'_{ij}$ is represented by $\lambda$, then
its composition on either side with a fixed 1-morphism is represented by
\begin{equation}\label{e:3.6}
 \begin{aligned}
 \Lambda*\id_{\Phi_{hi}}&:\quad(\phi_{hi}^*\lambda)\circ\chi_{hi},\\
 \id_{\Phi_{jk}}*\Lambda&:\quad(\phi_{ij}^*\mathrm d\phi_{jk})\circ\lambda.
 \end{aligned}
\end{equation}
(We restrict here domains wherever necessary.)

One may obtain general horizontal composition from these operations and vertical composition: if
$\Lambda_{ij}:\Phi_{ij}\Rightarrow\Phi'_{ij}$ and
$\Lambda_{jk}:\Phi_{jk}\Rightarrow\Phi'_{jk}$, then
\[
 \Lambda_{jk}*\Lambda_{ij}
 =(\id_{\Phi'_{jk}}*\Lambda_{ij})\odot(\Lambda_{jk}*\id_{\Phi_{ij}}).
\]
These operations give a \textit{strict}  2-category of Kuranishi neighbourhoods over a fixed open
subset of $X$. 

A \textit{coordinate change} is defined as an equivalence in this category: it has
an inverse up to invertible 2-morphisms on both sides. 

\begin{dfn}\label{r:3.6}
Let $X$ be Hausdorff and second countable. A \emph{Kuranishi structure with
trivial isotropy} of virtual dimension $n$ on $X$ consists of:
\begin{enumerate}
\item[(i)] Kuranishi neighbourhoods $\mathcal U_i$ of virtual dimension $n$,
indexed by a set $I$, whose images
 $S_i$ cover $X$;
\item[(ii)] coordinate changes $\Phi_{ij}:\mathcal U_i\to\mathcal U_j$ over
$S_i\cap S_j$, for all $i,j\in I$, with $\Phi_{ii}=\id$;
\item[(iii)] 2-morphisms $\Lambda_{ijk}:\Phi_{jk}\circ\Phi_{ij}\Rightarrow
\Phi_{ik}$ over $S_i\cap S_j\cap S_k$, for all $i,j,k\in I$.
\end{enumerate}
We require $\Lambda_{iij}$ and $\Lambda_{ijj}$ to be identities and, on every
intersection $S_i\cap S_j\cap S_k\cap S_l$, require
\begin{equation}\label{e:3.7}
 \Lambda_{ikl}\odot(\id_{\Phi_{kl}}*\Lambda_{ijk})
 =\Lambda_{ijl}\odot(\Lambda_{jkl}*\id_{\Phi_{ij}}).
\end{equation}
(Both sides are 2-morphisms from $\Phi_{kl}\circ\Phi_{jk}\circ\Phi_{ij}$ to
$\Phi_{il}$.) The resulting space is an object of Joyce's 2-category
$\mathbf{mKur}$ \cite[Theorem 4.59]{JoyceKur}.
\end{dfn}

Local 1-morphisms and 2-morphisms satisfy the following gluing theorem
of Joyce \cite[Theorem~4.14(a)]{JoyceKur}. Fix Kuranishi neighbourhoods
$\mathcal U=(V,E,s,\psi)$ and
$\mathcal U'=(V',E',s',\psi')$ on $X$, and an open subset
\[
 S\subseteq\operatorname{im}\psi\cap\operatorname{im}\psi'.
\]
A 1-morphism $\Phi:\mathcal U\to\mathcal U'$ over $S$ is also a
1-morphism over any open subset $T\subseteq S$, since its domain
contains $\psi^{-1}(T)$. We denote this restriction by $\Phi|_T$.

Let $\{S_\alpha:\alpha\in A\}$ be an open cover of $S$, and put
$S_{\alpha\beta}=S_\alpha\cap S_\beta$ and
$S_{\alpha\beta\gamma}=S_\alpha\cap S_\beta\cap S_\gamma$.
Then the following statements hold.
\begin{enumerate}
\renewcommand{\labelenumi}{(\roman{enumi})}
\item Let $\Phi,\Psi:\mathcal U\to\mathcal U'$ be 1-morphisms over $S$.
Suppose that for each $\alpha\in A$ we are given a 2-morphism
\[
 \eta_\alpha:\Phi|_{S_\alpha}\Longrightarrow\Psi|_{S_\alpha},
\]
and that
\[
 \eta_\alpha|_{S_{\alpha\beta}}
 =\eta_\beta|_{S_{\alpha\beta}}
 \qquad\text{for all }\alpha,\beta\in A.
\]
There is a unique 2-morphism $\eta:\Phi\Rightarrow\Psi$ over $S$
such that $\eta|_{S_\alpha}=\eta_\alpha$ for every $\alpha\in A$.

\item Suppose that for each $\alpha\in A$ we are given a 1-morphism
$\Phi_\alpha:\mathcal U\to\mathcal U'$ over $S_\alpha$, and for each
$\alpha,\beta\in A$ a 2-morphism
\[
 \theta_{\alpha\beta}:
 \Phi_\alpha|_{S_{\alpha\beta}}
 \Longrightarrow\Phi_\beta|_{S_{\alpha\beta}}.
\]
Assume that $\theta_{\alpha\alpha}=\operatorname{id}_{\Phi_\alpha}$
and that the cocycle condition
\[
 \theta_{\beta\gamma}|_{S_{\alpha\beta\gamma}}
 \odot\theta_{\alpha\beta}|_{S_{\alpha\beta\gamma}}
 =\theta_{\alpha\gamma}|_{S_{\alpha\beta\gamma}}
\]
holds for all $\alpha,\beta,\gamma\in A$. Here $\odot$ denotes
vertical composition of 2-morphisms.

There exist a 1-morphism $\Phi:\mathcal U\to\mathcal U'$ over $S$
and 2-morphisms
\[
 \eta_\alpha:\Phi_\alpha\Longrightarrow\Phi|_{S_\alpha}
 \qquad(\alpha\in A)
\]
such that
\[
 \eta_\beta|_{S_{\alpha\beta}}\odot\theta_{\alpha\beta}
 =\eta_\alpha|_{S_{\alpha\beta}}
 \qquad\text{for all }\alpha,\beta\in A.
\]
If $\widetilde\Phi$ and $\widetilde\eta_\alpha$ are another such
choice, there is a unique 2-morphism
$\zeta:\Phi\Rightarrow\widetilde\Phi$ satisfying
\[
 \zeta|_{S_\alpha}\odot\eta_\alpha=\widetilde\eta_\alpha
 \qquad\text{for every }\alpha\in A.
\]
If all the $\Phi_\alpha$ are coordinate changes, then $\Phi$ is
also a coordinate change.
\end{enumerate}

For the definitions of 1-morphisms and 2-morphisms between Kuranishi
spaces, we also refer to \cite[Definitions~4.20--4.21]{JoyceKur}.

\subsection{Kuranishi atlases in the sense of Borisov and Joyce}\label{ss:3.3}
 We recall the following description of Kuranishi atlases in
\cite[Definitions 2.13--2.17]{BJ}.

\begin{dfn}\label{r:3.7}
A \emph{Kuranishi atlas} of virtual dimension $n$ on $X$ consists of a
partially ordered set $(A,\prec)$, Kuranishi neighbourhoods
$\mathcal U_J=(V_J,E_J,s_J,\psi_J)$ of virtual dimension $n$, and coordinate
changes $\Phi_{JK}=(V_{JK},\phi_{JK},\chi_{JK})$ for $J\prec K$, subject
to the following conditions.
\begin{enumerate}
\item[(i)] The open subsets $S_J=\im\psi_J$ cover $X$. If $J\ne K$ and
$S_J\cap S_K\ne\varnothing$, then $J\prec K$ or $K\prec J$.
\item[(ii)] The coordinate change $\Phi_{JK}$ is defined over the full
intersection $S_J\cap S_K$, satisfies
\[
 \chi_{JK}s_J=\phi_{JK}^*s_K,\qquad
 \psi_K\phi_{JK}=\psi_J\quad\text{on }s_J^{-1}(0)\cap V_{JK},
\]
and satisfies the exactness condition \eqref{e:3.4}. Note that
$s_J^{-1}(0)\cap V_{JK}=\psi_J^{-1}(S_J\cap S_K)$.
\item[(iii)] If $J\prec K\prec L$, then
\begin{equation}\label{e:3.11}
 \phi_{JL}=\phi_{KL}\circ\phi_{JK},\qquad
 \chi_{JL}=(\phi_{JK}^*\chi_{KL})\circ\chi_{JK}
\end{equation}
on the common open domain
\begin{equation}\label{e:3.12}
 V_{JK}\cap V_{JL}\cap\phi_{JK}^{-1}(V_{KL}).
\end{equation}
This domain contains every zero over $S_J\cap S_K\cap S_L$.
\end{enumerate}
The atlas is \emph{locally finite} if the open subsets $S_J$ form a locally finite
cover. We use locally finite atlases throughout this paper.
\end{dfn}

Let $W_J\subseteq V_J$ be open subsets such that
$\{\psi_J(s_J^{-1}(0)\cap W_J):J\in A\}$ covers $X$. Restrict the
$J$th Kuranishi neighbourhood to $W_J$, and its coordinate changes to
\[
 W_{JK}=W_J\cap V_{JK}\cap\phi_{JK}^{-1}(W_K).
\]
These restrictions form a Kuranishi atlas on $X$.

We shall also use the following purely topological result.

\begin{prop}\label{r:3.8}\label{prop:closed-refinement}
Let $Z$ be a paracompact Hausdorff space and $\{R_i:i\in I\}$ an open cover.
For each finite nonempty subset $J\subseteq I$ there is a closed subset
$C_J\subseteq Z$ such that
\[
 C_J\subseteq\bigcap_{j\in J}R_j,\qquad
 Z=\bigcup_{\varnothing\ne J\subseteq I\text{ finite}}C_J.
\]
The family $\{C_J\}$ is locally finite, and
\[
 C_J\cap C_K\ne\varnothing\quad\Longrightarrow\quad
 J\subseteq K\text{ or }K\subseteq J.
\]
\end{prop}

This is \cite[Proposition 3.13]{BJ}, proved in \cite[\S6.1]{BJ}.

\subsection{Derived smooth manifolds and orientations}\label{ss:3.4}
We use Joyce's 2-category $\mathbf{dMan}$ of \emph{d-manifolds}
\cite[\S4]{JoyceSummary}, which we call derived smooth manifolds.
Every d-manifold is locally equivalent to a \emph{standard model}
$\boldsymbol S_{V,E,s}$ associated to a smooth manifold $V$, a smooth
vector bundle $E\to V$, and a smooth section $s\in C^\infty(V,E)$.
We recall its definition \cite[Example~4.4]{JoyceSummary}.

A $C^\infty$-ring is a set $A$ with an operation $A^k\to A$ for every
smooth function $\RR^k\to\RR$, respecting projections and composition.

\Needspace{10\baselineskip}
Write $\mathcal C_V$ for the sheaf of smooth functions on $V$ and
$\mathcal E_V$ for the sheaf of smooth sections of $E^*$. Let
$I_s\subseteq\mathcal C_V$ be the image of the evaluation map
$\mathcal E_V\to\mathcal C_V$, $\alpha\mapsto\alpha(s)$.
Put $Z=s^{-1}(0)$, write $i:Z\hookrightarrow V$ for the inclusion, and define
\[
 \mathcal O_Z=i^{-1}(\mathcal C_V/I_s),\qquad
 \mathcal O'_Z=i^{-1}(\mathcal C_V/I_s^2),\qquad
 \mathcal E_Z=i^{-1}(\mathcal E_V/I_s\mathcal E_V).
\]
Then $\underline Z=(Z,\mathcal O_Z)$ is \textit{a $C^\infty$-scheme} in the sense of \cite[Definition~2.4]{JoyceSummary}.
There is a natural surjection
$\iota_Z:\mathcal O'_Z\to\mathcal O_Z$ with square-zero kernel
$\mathcal I_Z=i^{-1}(I_s/I_s^2)$. There is also a natural map
\[
 \jmath_Z:\mathcal E_Z\longrightarrow\mathcal I_Z,
 \qquad [\alpha]\longmapsto[\alpha(s)].
\]
The tuple
\[
 \boldsymbol S_{V,E,s}
   =(\underline Z,\mathcal O'_Z,\mathcal E_Z,\iota_Z,\jmath_Z)
\]
is a \textit{d-manifold} of virtual dimension $\dim V-\rk E$, called the
\emph{standard model} of $(V,E,s)$.

Let $\Omega^1_V$ denote the sheaf of smooth one forms on $V$, and choose
any smooth connection $\nabla$ on $E$. The virtual cotangent bundle of
$\boldsymbol S_{V,E,s}$ is represented by the complex of
$\mathcal O_Z$-modules
\begin{equation}\label{e:3.13}
 \bigl[\mathcal E_Z\xrightarrow{(\nabla s)^*}
       i^{-1}(\Omega^1_V/I_s\Omega^1_V)\bigr]
\end{equation}
in degrees $-1,0$. The induced map in \eqref{e:3.13} is independent of
the choice of $\nabla$: 
see \cite[Example~4.4]{JoyceSummary}.

\begin{thm}\label{r:3.9}\label{thm:atlas-reconstruction}
Let $X$ be Hausdorff and second countable, with a Kuranishi atlas $\mathcal K$
of virtual dimension $n$ in Definition~\ref{r:3.7}. Then $\mathcal K$
determines a Kuranishi space in $\mathbf{mKur}$ and a derived smooth manifold
$M$ in $\mathbf{dMan}$, each with underlying space $X$ and virtual dimension
$n$, canonical up to equivalence. 

Orientations of $\mathcal K$ correspond naturally to
orientations of each of these objects.

If $Z$ is a manifold, $\pi:X\to Z$ is continuous, and there are smooth maps
$\pi_J:V_J\to Z$ satisfying
\begin{equation}\label{e:3.14}
 \pi_J|_{s_J^{-1}(0)}=\pi\circ\psi_J,\qquad
 \pi_K\circ\phi_{JK}=\pi_J\quad\text{on }V_{JK},
\end{equation}
then these data determine a morphism $M\to Z$, canonical up to 2-isomorphism,
inducing $\pi$ on underlying spaces.
\end{thm}

This follows from \cite[Theorem~2.18]{BJ} and
\cite[Theorem~4.67]{JoyceKur}. 

We now explain \textit{orientations}. For a chart $\mathcal U_J$ put
\[
 K_J=(\det T^*V_J\otimes\det E_J)|_{s_J^{-1}(0)}.
\]
Taking determinants in \eqref{e:3.4}, with the sign convention of
\cite[Definition 2.16]{BJ}, we obtain an isomorphism from $K_J$ to the pullback
of $K_K$ on the zero locus of an overlap. These isomorphisms are compatible
with composition, and they glue the $K_J$ to a real
line bundle $K_M$ on $X$, the \emph{canonical bundle} of $M$.

\begin{dfn}\label{r:3.10}
An \emph{orientation} of a Kuranishi atlas is an orientation of each $K_J$,
continuous on $s_J^{-1}(0)$, preserved by the isomorphisms on overlaps. An
\emph{orientation} of $M$ is an orientation of $K_M$. Reversing it, one gets \textit{the
oppositely oriented derived manifold}, written $-M$.
\end{dfn}

The tangent and obstruction spaces in \eqref{e:3.1} are identified by
coordinate changes and give spaces $T_xM$ and $O_xM$ independent of the
chart. The exact sequence
\[
 0\longrightarrow T_xM\longrightarrow T_vV_J
 \xrightarrow{\mathrm ds_J}E_{J,v}\longrightarrow O_xM\longrightarrow0
\]
implies
\begin{equation}\label{e:3.15}
 K_M|_x\cong\det T_x^*M\otimes\det O_xM.
\end{equation}
If $E_J=0$, this recovers the ordinary
orientation of a manifold. 

\subsection{Families, boundaries, and collars}\label{ss:3.5}

\subsubsection{Relative atlases and pullbacks}\label{sss:3.5.1}
Let $Z$ be a manifold without boundary and $\pi:X\to Z$ a continuous
map. A \emph{relative Kuranishi atlas} for $\pi$ is a Kuranishi atlas
on $X$ for which each $V_J$ is a manifold without boundary, together with
smooth maps $\pi_J:V_J\to Z$ satisfying
\eqref{e:3.14}. Write $M\to Z$ for the morphism given by
Theorem~\ref{r:3.9}. Suppose that every $\pi_J$ is a submersion, and
put $d=\vdim M-\dim Z$. The relative virtual cotangent bundle is
represented on the $J$th chart by
\[
 [E_J^*|_{s_J^{-1}(0)}\longrightarrow
       (\ker\mathrm d\pi_J)^*|_{s_J^{-1}(0)}].
\]
Its determinant lines glue to a real line bundle $K_{M/Z}$ on $X$.
A \emph{relative orientation} is an orientation of $K_{M/Z}$.
We use the determinant isomorphism
\[
 K_M\cong\pi^*\det T^*Z\otimes K_{M/Z},
\]
and note that an orientation of $Z$
and a relative orientation determine an orientation of $M$.

\begin{prop}\label{r:3.11}\label{prop:family-pullback}
Let $M\to Z$ be given by a relative Kuranishi atlas with submersions
$\pi_J:V_J\to Z$, and let $g:T\to Z$ be a smooth map
of manifolds without boundary. Then
\begin{equation}\label{e:3.16}
 V_J^T=V_J\times_Z T,\qquad E_J^T=\operatorname{pr}_1^*E_J,
 \qquad s_J^T=\operatorname{pr}_1^*s_J,
\end{equation}
with the pulled back coordinate changes, form a relative Kuranishi
atlas on $X\times_Z T$. It determines the fibre product
$M_T=M\times_Z T$ in $\mathbf{dMan}$, with
\[
 |M_T|=X\times_Z T,\qquad \vdim M_T=d+\dim T.
\]
If $p:M_T\to M$ is the projection, then
$K_{M_T/T}\cong p^*K_{M/Z}$. 
\end{prop}
\begin{proof}
The submersion hypothesis makes each $V_J\times_Z T$ an ordinary
manifold. Note that this base change preserves the coordinate change identities and the exactness condition
\eqref{e:3.4}. 

The associated standard
models are the local fibre products in $\mathbf{dMan}$ by
\cite[Theorems~4.22--4.23]{JoyceSummary}, hence they present
$M\times_Z T$. 

The rest is clear.

See also \cite[Theorem~3.24(iii)]{BJ}.
\end{proof}

\subsubsection{Derived smooth manifolds with boundary}\label{sss:3.5.2}
A \emph{derived smooth manifold with boundary} means a d-manifold
with boundary in Joyce's 2-category $\mathbf{dMan}^{\mathrm b}$
\cite[Definition~7.1]{JoyceSummary}. Its local models are
$\boldsymbol S_{V,E,s}$ with $V$ a manifold with boundary, $E\to V$
a smooth vector bundle, and $s$ a smooth section. The structure
includes a boundary $\partial M$ and a boundary morphism
$i_M:\partial M\to M$. In a standard model these are
\[
 \partial\boldsymbol S_{V,E,s}
 \cong\boldsymbol S_{\partial V,E|_{\partial V},s|_{\partial V}}
\]
and the morphism induced by $\partial V\hookrightarrow V$
\cite[Example~7.3]{JoyceSummary}. Notice that $\partial M$ is a derived
smooth manifold without boundary, of virtual dimension $\vdim M-1$.

An orientation of $M$ is an orientation of its canonical line bundle,
as in Definition~\ref{r:3.10}. We orient $\partial M$ by placing the
outward normal first. More precisely, on each chart the isomorphism
\[
 (\det T^*V\otimes\det E)|_{\partial V}
 \cong N_V^*\otimes
       (\det T^*(\partial V)\otimes\det E|_{\partial V})
\]
determines the boundary orientation, where $N_V^*$ is the conormal
line, oriented by covectors positive on outward pointing normal
vectors. These orientations agree on overlaps
\cite[Theorem~7.25]{JoyceSummary}.

\subsubsection{Kuranishi atlases with boundary}\label{sss:3.5.3}
Let $X$ be Hausdorff and second countable. A \emph{Kuranishi atlas
with boundary} on $X$ consists of the data in
Definition~\ref{r:3.7}, with the $V_J$ allowed to be manifolds with
boundary, satisfying the conditions of that definition and the
following additional conditions on each coordinate map:
\begin{enumerate}
\item[(i)] $\phi_{JK}^{-1}(\partial V_K)=V_{JK}\cap\partial V_J$;
\item[(ii)] near each point of $V_{JK}\cap\partial V_J$, local
boundary defining functions $\rho_J,\rho_K$ satisfy
\[
 \rho_K\circ\phi_{JK}=a_{JK}\rho_J
 \qquad\text{for a smooth function }a_{JK}>0.
\]
\end{enumerate}
(Here a boundary defining function is a nonnegative smooth function
whose zero set is the boundary and whose differential is nonzero
there. Conditions (i)--(ii) say that $\phi_{JK}$ is \textit{simple} and \textit{flat}
in the terminology of \cite[Definition~5.4]{JoyceSummary}.)

Put
\[
 \partial X=\bigcup_J\psi_J(s_J^{-1}(0)\cap\partial V_J),\qquad
 \partial V_{JK}=V_{JK}\cap\partial V_J.
\]
The restrictions
\[
 \mathcal U_J^\partial=
 (\partial V_J,E_J|_{\partial V_J},s_J|_{\partial V_J},
   \psi_J|_{s_J^{-1}(0)\cap\partial V_J})
\]
and the restricted coordinate changes form a Kuranishi atlas on
$\partial X$, of virtual dimension one less than the original atlas.

By \cite[Theorem~7.10]{JoyceSummary}, the atlas determines a derived
smooth manifold $M$ with boundary, with $|M|=X$, canonical up to
equivalence. Its boundary is represented by the atlas just described.

\subsubsection{Collars}\label{sss:3.5.4}
A \emph{collar} of a manifold $V$ with boundary is a diffeomorphism
\[
 c:\partial V\times[0,\epsilon)\longrightarrow U\subseteq V,
 \qquad c(v,0)=v,
\]
where $U$ is an open neighbourhood of $\partial V$.

A collar of a \textit{derived} smooth manifold $M$ with boundary is an
equivalence from $\partial M\times[0,\epsilon)$ to an open
neighbourhood $N$ of $\partial M$ in $M$, such that its composite
with $N\hookrightarrow M$, restricted at $0$, is 2-isomorphic
to $i_M$.

A Kuranishi atlas on a derived manifold $M$ with boundary has \emph{compatible collars}
if, for some $\epsilon>0$, there are collars $c_J$ of all $V_J$ and
bundle identifications over them, restricting to the identity at
$t=0$, with the following properties.
Write $E_J^\partial=E_J|_{\partial V_J}$,
$s_J^\partial=s_J|_{\partial V_J}$, and let $q_J$ be projection
from $\partial V_J\times[0,\epsilon)$ to $\partial V_J$. We require
\begin{enumerate}
\item[(i)] $c_J^*E_J\cong q_J^*E_J^\partial$ and
$c_J^*s_J=q_J^*s_J^\partial$;
\item[(ii)] $c_J^{-1}(V_{JK})=\partial V_{JK}\times[0,\epsilon)$;
\item[(iii)] on this domain, the coordinate changes take the form
\[
 \phi_{JK}(x,t)=(\phi_{JK}^\partial(x),t),\qquad
 \chi_{JK}(x,t)=\chi_{JK}^\partial(x)
\]
in the collar coordinates and bundle identifications.
\end{enumerate}
Here $\phi_{JK}^\partial,\chi_{JK}^\partial$ are the restricted
coordinate changes.

It is easy to see that in such a situation gluing gives a
collar of $M$.

\subsubsection{Pullbacks with boundary}\label{sss:3.5.5}
We can now extend Proposition~\ref{r:3.11} to  manifolds
with boundary.

\begin{cor}\label{r:3.12}\label{cor:family-bordism}\label{prop:collared-pullback}
Let $M\to Z$ satisfy the hypotheses of Proposition~\ref{r:3.11},
and let $g:T\to Z$ be a smooth map from a manifold with boundary.
Then the following hold.
\begin{enumerate}
\item[(i)] The data \eqref{e:3.16} form a Kuranishi atlas with
boundary. They determine $M_T=M\times_Z T$ in
$\mathbf{dMan}^{\mathrm b}$, with
\[
 |M_T|=X\times_Z T,\qquad \vdim M_T=d+\dim T,
 \qquad \partial M_T\simeq M\times_Z\partial T.
\]

Given a relative orientation of
$M\to Z$ and an orientation of $T$, we orient $M_T$ by
\[
 K_{M_T}\cong\pi_T^*\det T^*T\otimes p^*K_{M/Z},
\]
where $\pi_T:M_T\to T$ and $p:M_T\to M$ are the projections.

\item[(ii)] Suppose $T$ has a collar $c$ such that
$g(c(y,t))=g(y)$ for $y\in\partial T$ and $0\leqslant t<\epsilon$.
Then the pulled back atlas has compatible collars. They
give a collar
\[
 (M\times_Z\partial T)\times[0,\epsilon)
 \simeq M_T|_{|\pi_T|^{-1}(\im c)}.
\]

\item[(iii)] Suppose $M\to Z$ is relatively oriented and proper
on underlying spaces. Let $\gamma:[0,1]\to Z$ be a smooth path
constant near its endpoints, and orient $[0,1]$ by $\mathrm dt$.
Then $W=M\times_Z[0,1]$ is compact, has virtual dimension $d+1$,
and has collars with boundary
\begin{equation}\label{e:3.17}
 \partial W\simeq -M_{\gamma(0)}\amalg M_{\gamma(1)}.
\end{equation}
Consequently, $W$ can be considered as an oriented bordism between these fibres.
\end{enumerate}
\end{cor}
\begin{proof}
This is a very straightforward check.
\end{proof}

\subsection{Bordism and virtual fundamental classes}\label{ss:3.6}
We finish by recalling the results of \cite[\S\S2.6.4--2.6.5]{BJ}. 

\begin{dfn}\label{r:3.13}
For a manifold $Y$ and $n\in\ZZ$, consider pairs $(M,f)$, where $M$ is a
compact oriented derived smooth manifold without boundary of virtual
dimension $n$, and $f:M\to Y$ is a morphism. Two pairs $(M_0,f_0)$ and
$(M_1,f_1)$ are \emph{bordant} if there are a compact oriented derived
manifold $W$ with boundary, of virtual dimension $n+1$, a morphism
$F:W\to Y$, and an oriented equivalence
$j:-M_0\amalg M_1\to\partial W$, together with a 2-isomorphism from
$f_0\amalg f_1$ to $F\circ i_W\circ j$. Here $i_W:\partial W\to W$ is the
boundary morphism.

The bordism classes form an abelian group $\mathrm dB_n(Y)$: addition is
disjoint union, the inverse reverses orientation, and the empty derived
manifold represents zero. 

The ordinary oriented bordism group $B_n(Y)$ is
defined in the same way using manifolds of dimension $n$, smooth maps, and
oriented diffeomorphisms of boundaries.
\end{dfn}

\begin{thm}\label{r:3.14}
For every manifold $Y$, inclusion of ordinary manifolds induces an isomorphism
$B_n(Y)\xrightarrow{\sim}\mathrm dB_n(Y)$ for $n\geqslant0$. For $n<0$,
$\mathrm dB_n(Y)=0$.
\end{thm}

This is \cite[Theorem 2.21]{BJ}. (Its geometric content is that a compact
oriented derived manifold can be perturbed, within its bordism class, to an
ordinary oriented manifold.)

Applying the fundamental class homomorphism on
ordinary bordism gives a map
\begin{equation}\label{e:3.18}
 \Pi_Y:\mathrm dB_n(Y)\longrightarrow H_n(Y;\ZZ),\qquad
 [N,h]\longmapsto h_*[N]
\end{equation}
when $N$ is an ordinary manifold. For negative $n$ this map is zero.
Consequently, bordant derived manifolds over $Y$ have the same virtual class
in $H_n(Y;\ZZ)$.

To place a virtual class in the homology of the underlying space itself,
we use the following topological hypothesis. A space is a \emph{Euclidean
neighbourhood retract}, or ENR, if it is homeomorphic to a retract of an
open subset of some $\RR^N$. Complex analytic spaces are locally ENRs, by
local analytic triangulation \cite{Lojasiewicz}. A Hausdorff space covered
by finitely many open ENRs is an ENR \cite[Theorem 18.4.6]{tomDieck}.
Consequently every compact subset of a Hausdorff complex analytic space has
an open ENR neighbourhood: take the union of finitely many such
neighbourhoods. In particular, a compact Hausdorff complex analytic space
is an ENR.

\begin{thm}\label{r:3.15}\label{thm:virtual-class}
Let $M$ be a compact oriented derived smooth manifold without boundary of
virtual dimension $n$, whose underlying space $X$ is an ENR. Then there is
a \emph{virtual fundamental class}
\[
 [M]^{\mathrm{vir}}\in H_n(X;\ZZ).
\]
It is preserved by oriented equivalences, and for a morphism $f:M\to Y$ to
a manifold one has
\begin{equation}\label{e:3.19}
 |f|_*[M]^{\mathrm{vir}}=\Pi_Y([M,f]).
\end{equation}
For an ordinary compact oriented manifold it is the fundamental class.
If $n<0$, it is zero.
\end{thm}

Here is the construction in \cite[\S2.6.5]{BJ}. Choose an embedding
$e:M\to\RR^N$ for $N$ sufficiently large
\cite[Corollary 4.31]{JoyceSummary}. By the ENR hypothesis, we may choose an open
neighbourhood $U$ of $e(X)$ and a retraction $r:U\to X$, where
$r\circ|e|=\id_X$. Define
\begin{equation}\label{e:3.20}
 [M]^{\mathrm{vir}}=r_*\Pi_U([M,e]).
\end{equation}
The class is independent of the embedding, neighbourhood, and retraction.
Formula~\eqref{e:3.19} follows by making the construction compatible with
the graph of $f$.

\begin{prop}\label{r:3.16}\label{prop:virtual-bordism-enr}
Let $W$ be a compact oriented derived smooth manifold of virtual dimension
$n+1$, with boundary $\partial W=-M_0\amalg M_1$. Suppose the underlying
spaces $X_a$ of $M_a$ are ENRs, and let $q:|W|\to T$ be a continuous map
to an ENR. Then
\[
 (q|_{X_0})_*[M_0]^{\mathrm{vir}}=(q|_{X_1})_*[M_1]^{\mathrm{vir}}
 \quad\text{in }H_n(T;\ZZ).
\]
\end{prop}
\begin{proof}
This follows from \cite[\S\S2.6.4--2.6.5]{BJ}, after embedding $T$
in Euclidean space as a neighbourhood retract and applying smooth
approximation.
\end{proof}

\section{Relative analytic Darboux coordinates}\label{sec:local}\label{sec:4}

We adapt to the analytic context here the  argument of Brav, Bussi and Joyce
\cite[\S5]{BBJ}. 

Notice that all coordinate changes and differentials in this
section are over $B$.
In particular $d_B$ annihilates functions pulled back from the base.

\begin{prop}\label{prop:presentation}\label{r:4.1}
Suppose that $\bL^{\mathrm{an}}_{\bX/B}$ has Tor-amplitude contained in
$[-2,0]$. Near a point $p$ there is a presentation over $B$ of the form
\begin{equation}\label{e:4.1}
 A^\sharp=\cO_U[z_1,\ldots,z_r,y_1,\ldots,y_n],
 \qquad |z_a|=-1,\quad |y_j|=-2,
\end{equation}
where $U\to B$ is a holomorphic submersion, with relative coordinates
$x_1,\ldots,x_m$, and
\[
 \delta z_a=f_a,\qquad \delta y_j=\sum_a b_{ja}z_a,
 \qquad \sum_a b_{ja}f_a=0.
\]
Its analytic realisation is equivalent to a neighbourhood of $p$ in
$\bX$, over $B$. It can be chosen minimal relative to $B$ at $p$:
$f_a(p)=0$, $d_Bf_a(p)=0$, and $b_{ja}(p)=0$.

The numbers $m,r,n$ are the dimensions of the degree $0,-1,-2$
cohomology of $\bL^{\mathrm{an}}_{\bX/B}\Lotimes\C_p$, respectively.

On a
sufficiently small relative polydisc, the de Rham complex is
\eqref{e:2.2}.
\end{prop}

\begin{proof}
 Write $\C_p$ for the residue field at $p$, and put
\[
 \begin{aligned}
 m&=\dim_{\C}H^0(\bL^{\an}_{\bX/B}|_p),&
 r&=\dim_{\C}H^{-1}(\bL^{\an}_{\bX/B}|_p),\\
 n&=\dim_{\C}H^{-2}(\bL^{\an}_{\bX/B}|_p),
 \end{aligned}
\]
where $|_p$ denotes the derived fibre. We construct the presentation
by adjoining $r$ generators in degree $-1$ and then $n$ generators
in degree $-2$, following the argument of \cite[\S4.1]{BBJ} with
analytic cotangent complexes.

\medskip
\noindent\emph{Step 1.}
By \cite[Lemma 12.13]{LurieIX}, a neighbourhood $\bX_0$ of $p$
admits a derived closed embedding into a complex manifold. Taking
its graph with the map to $B$, and choosing coordinates, one gets
a closed embedding $j:\bX_0\to V$ over $B$, where
$V\subset B\times\C^N$ is open. After a linear change of
coordinates in $\C^N$, we may suppose that the differentials of
the first $m$ coordinates $x_1,\ldots,x_m$ form a basis of
$H^0(\bL^{\an}_{\bX/B}|_p)$. Let
$\pi:B\times\C^N\to B\times\C^m$ be the projection which
saves only the base coordinates and $x_1,\ldots,x_m$, and put
$g=\pi\circ j$.

The conormal sequence and the holomorphic implicit function theorem
show that, after shrinking $\bX_0$ near $p$, the classical
truncation $t_0(g)$ is an embedding. Choose a sufficiently small
relative polydisc $U\subset B\times\C^m$ containing $g(p)$ such
that the induced map $t_0(g)^{-1}(U)\to U$ is a closed embedding.
Define
\[
 \bX_1=\bX_0\times^{\R}_{B\times\C^m}U.
\]
Since $U$ is open, this is the open derived subspace of $\bX_0$
on $t_0(g)^{-1}(U)$, with the analytic structure restricted from
$\bX_0$. Thus $g$ restricts to a morphism $\bX_1\to U$ over
$B$ whose classical truncation is a closed embedding. We write
$\bX$ for $\bX_1$ throughout the rest of the proof.

\medskip
\noindent\emph{Step 2.}
We first work with germs at $p$. 

Set $R=\cO_{U,p}$ and
$T=\cO_{\bX,p}^{\mathrm{alg}}$, with the map $R\to T$ induced
by Step 1. Write $K_0=(\bL^{\an}_{\bX/U})_p$.
The transitivity gives a distinguished triangle
\[
 T\Lotimes_R\Omega^1_{U/B,p}
 \longrightarrow (\bL^{\an}_{\bX/B})_p
 \longrightarrow K_0.
\]

Put $Q_0=\cofib(R\to T)$. The map $R\to H^0(T)$ is
surjective, so $Q_0$ has cohomology only in degrees at most $-1$.
The analytic Hurewicz comparison
\cite[Proposition 5.34]{PortaYu}, followed by tensoring with
$\C_p$, gives
\[
 H^{-1}(Q_0)\otimes_R\C_p
 \cong H^{-1}(K_0\Lotimes_T\C_p).
\]
Choose classes $q_1,\ldots,q_r\in H^{-1}(Q_0)$ lifting a basis
of the right hand side. Each class is represented by a function
$f_a\in R$ and a specified homotopy from its image in $T$ to
zero. Adjoin $z_a$ in degree $-1$ with $\delta z_a=f_a$.
The homotopies give a morphism $A_1=R[z_1,\ldots,z_r]\to T$
extending $R\to T$.

The finite $R$-module $H^{-1}(Q_0)$ is generated by the $q_a$,
by Nakayama's lemma. Its natural map onto
$I=\ker(R\to H^0(T))$ sends $q_a$ to $f_a$, so the $f_a$
generate $I$. Consequently $H^0(A_1)\to H^0(T)$ is an
isomorphism. 

If $K_1$ denotes the new relative analytic cotangent
complex at $p$, one has
\[
 T^{\oplus r}[1]\longrightarrow K_0\longrightarrow K_1.
\]

Note that $K_1$ has cohomology only in degrees
at most $-2$, and $K_1\Lotimes_T\C_p$ has cohomology only
in degree $-2$, of dimension $n$.

\medskip
\noindent\emph{Step 3.}
Since $H^0(A_1)\cong H^0(T)$, \cite[Corollary 5.35]{PortaYu}
implies that $Q_1=\cofib(A_1\to T)$ has cohomology only in
degrees at most $-2$. Applying the Hurewicz comparison in this
degree, choose $n$ classes in $H^{-2}(Q_1)$ whose images form
a basis of $H^{-2}(K_1\Lotimes_T\C_p)$.
Their representatives consist of closed elements
$g_j\in A_1^{-1}$ and specified homotopies from their images
in $T$ to zero.

Write $g_j=\sum_a b_{ja}z_a$, with $b_{ja}\in R$.
Adjoin $y_j$ in degree $-2$ with $\delta y_j=g_j$ to obtain
$A_2=A_1[y_1,\ldots,y_n]\to T$. Closedness of $g_j$ says
precisely that $\sum_a b_{ja}f_a=0$.

For the new relative cotangent complex $K_2$, one has, once again,
$$T^{\oplus n}[2]\to K_1\to K_2.$$ Note that $K_2\Lotimes_T\C_p=0$.
As $K_2$ is perfect, by the derived Nakayama,  $K_2=0$.

After shrinking, perhaps, $U$, one gets in this manner a
sheaf of cdgas $A$ on $U$ and a morphism
$\bX\to\Specan A$ over $B$.
After a further restriction near $p$, this morphism induces an
isomorphism on classical truncations and has zero relative analytic
cotangent complex. It is therefore an
equivalence by \cite[Corollary 5.37]{PortaYu}.

\medskip
\noindent\emph{Step 4.}
Let $Z\subset U$ be the classical locus, identified with
$t_0(\bX)$, and write $i:Z\hookrightarrow U$ for its inclusion.

Note that $i^{-1}\Omega^1_{\hol}(A/B)$ is isomorphic $\bL^{\an}_{\bX/B}$,
as modules over $i^{-1}A\simeq\cO_{\bX}^{\mathrm{alg}}$.
Here $\Omega^1_{\hol}(A/B)$ is the free graded $A$-module on
$d_Bx_i,d_Bz_a,d_By_j$, with differential induced by $\delta$.
Its derived fibre at $p$ is
a complex
\[
 \C_p^{\oplus n}\longrightarrow\C_p^{\oplus r}
 \longrightarrow\C_p^{\oplus m}
\]
in degrees $-2,-1,0$. By the equivalence just proved, its
cohomology has dimensions $n,r,m$ in those degrees. Both
differentials must therefore vanish. Their matrices are, up to
the differential convention, $(b_{ja}(p))$ and
$(\partial f_a/\partial x_i(p))$. Also $f_a(p)=0$, since
$f_a\in I$. This proves the minimality.

\medskip
\noindent\emph{Step 5.}
The assertion about de Rham complex is now straightforward.
\end{proof}

\begin{thm}
\label{thm:darboux}\label{r:4.2}
Let $\omega$ be a relative $(-2)$-shifted symplectic form on $\bX/B$.
Suppose that it satisfies Definition~\ref{r:2.7}.
Then near every point there is a presentation over $B$ with
\begin{equation}\label{e:4.2}
 \begin{gathered}
 A^\sharp=\cO_U[z_1,\ldots,z_r,y_1,\ldots,y_m],\\
 \omega_{\mathrm{st}}=\sum_i d_Bx_i\,d_By_i+
                  \sum_a d_Bz_a\,d_Bz_a,\\
 \delta z_a=\tfrac12 s_a,\qquad
 \delta y_i=\sum_a\frac{\partial s_a}{\partial x_i}z_a,
 \qquad \sum_a s_a^2=0,
 \end{gathered}
\end{equation}
where the $s_a$ are holomorphic in both base and relative coordinates.
The given closed relative form is homotopic to
$(\omega_{\mathrm{st}},0,\ldots)$ as a closed relative form.
Conversely, a chart with these formulas is locally exact in the sense
of Definition~\ref{r:2.7}.
\end{thm}

\begin{proof}
Note that the perfect
relative cotangent complex has amplitude $[-2,0]$. Apply
Proposition~\ref{prop:presentation}, choosing a presentation
minimal at the given point. 

Passing to a relative neighborhood $U$, choose, by exactness,
$\eta\in C_B(A;U)^{-1}$ with
$D\eta=-2\omega$. Decompose it by form degree as
$\eta=\Phi+\phi+\eta_{\geq2}$. The form degree $0$ and $1$
equations are
\begin{equation}\label{e:4.3}
 \delta\Phi=0,\qquad d_B\Phi+\delta\phi=0.
\end{equation}
Set $\theta=-\frac12d_B\phi$. Then
\[
 (\theta,0,\ldots)=\omega+\tfrac12D\eta_{\geq2},
 \qquad \delta\theta=d_B\theta=0.
\]

Let $E$ be the Euler derivation, so that $E(z_a)=-z_a$,
$E(y_j)=-2y_j$, and $E$ vanishes on all holomorphic degree zero
functions, including the base coordinates. 

On a form of internal
degree $k$, the Cartan identity reads
$(d_B\iota_E+\iota_Ed_B)\alpha=k\alpha$.
Consequently $\phi_E=\iota_E\theta$ satisfies
$d_B\phi_E=-2\theta$. The difference $\beta=\phi_E-\phi$ is a
$d_B$-closed $1$-form of internal degree $-2$, so
$\psi=-\frac12\iota_E\beta$ satisfies $d_B\psi=\beta$.
Replace $\phi$ by $\phi+d_B\psi$ and $\Phi$ by $\Phi+\delta\psi$.
Since $\delta d_B=-d_B\delta$, these replacements preserve both
equations in
\eqref{e:4.3}. We have therefore arranged
$\phi=\phi_E$.

Minimality and nondegeneracy imply $n=m$. They also imply that the
$d_Bx_i\,d_By_j$ coefficient matrix $P$ and the symmetric
$d_Bz_a\,d_Bz_b$ coefficient matrix $Q$ of $\theta$ are invertible
at the chosen point, and hence on a neighbourhood. Its Euler primitive has the form
\[
 \phi_E=2\sum_{i,j}P_{ij}y_jd_Bx_i
       -2\sum_{a,b}Q_{ab}z_ad_Bz_b
       +\sum_i C_i(z)d_Bx_i,
\]
where each $C_i$ is quadratic in the odd variables and all coefficients
are holomorphic in the base and relative coordinates. Locally there
is a holomorphic invertible matrix $M$ with $M^tM=Q$: conjugate $Q$ at
the point to the identity by a constant matrix and take the normally
convergent square-root series on a sufficiently small neighbourhood.
This operation is holomorphic in all parameters.

Put $Z=Mz$. Then
$\sum_a Z_ad_BZ_a-\sum_{a,b}Q_{ab}z_ad_Bz_b$
has no $d_Bz$ terms, and so equals $\sum_i B_i(z)d_Bx_i$ with
quadratic $B_i$. Define
$Y_i=\sum_jP_{ij}y_j+\frac12C_i+B_i$. These are invertible finite
graded changes of coordinates over $B$, and they give
\[
 \phi_E=2\sum_iY_id_Bx_i-2\sum_aZ_ad_BZ_a,
 \qquad
 \theta=\sum_i d_Bx_i\,d_BY_i+\sum_a d_BZ_a\,d_BZ_a.
\]

Write $\Phi=\sum_as_aZ_a$, $\delta Z_a=f_a$, and
$\delta Y_i=\sum_ag_{ia}Z_a$. In the second equation of \eqref{e:4.3},
comparison of the coefficients of $d_BZ_a$ gives $s_a=2f_a$.
Comparing  the coefficients of $Z_ad_Bx_i$, we get
$g_{ia}=\partial_i s_a$. Here we use $\delta d_B=-d_B\delta$
and the total degrees $|d_Bx_i|=1$, $|d_BZ_a|=0$.
The first equation then implies
$0=\delta\Phi=\frac12\sum_as_a^2$ in $\cO_U$.
Renaming $Z,Y$ proves \eqref{e:4.2}.

Conversely, for \eqref{e:4.2} put
$\Phi=\sum_as_az_a$ and
$\phi=2\sum_iy_id_Bx_i-2\sum_az_ad_Bz_a$.
The same calculation gives
$\delta\Phi=0$, $d_B\Phi+\delta\phi=0$, and
$d_B\phi=-2\omega_{\mathrm{st}}$.
Thus $D(\Phi+\phi)=-2\omega_{\mathrm{st}}$, proving local exactness.
\end{proof}

\section[Compatible cell presentations]{Compatible \textcolor{black}{cell} presentations}\label{sec:comparisons}\label{sec:5}

\begin{lem}\label{lem:cells}\label{r:5.1}
Let $W''\Subset W'\Subset W$ be Stein open subsets of a complex
manifold smooth over $B$.

Let $\bY$ be a complex-analytic space. Suppose that $Y=t_0(\bY)$ is a closed
analytic subspace of $W$ and that $\bY$ has coherent homotopy sheaves.

Let $P$ be a cell presentation over $W$, with cells of degrees
$-1,-2$, and a map $\bY\to\Specan P$ over $W$.

Then, after restriction to $W''$, this map factors through a cell
extension $P\to P'$ with new cells of degrees $-1,-2$, for which
$\cofib(P'\to i_*\cO_{\bY})$ is $3$-connective.

\end{lem}
\begin{proof}
\textit{Step 1.} Put $Q=\cofib(P\to i_*\cO_{\bY})$. Since $Y\hookrightarrow W$ is a closed embedding,
$Q$ is $1$-connective, and its cohomology sheaves are coherent.
Now, by Cartan's theorem B,
$H^j\R\Gamma(W,Q)=\Gamma(W,\cH^jQ)$, and, by Cartan's theorem A and compactness,
finitely many sections of $\cH^{-1}Q$ generate it near
$\overline{W'}$ \cite{GR}.

Write $T=i_*\cO_{\bY}^{\mathrm{alg}}$ and let $\phi:P\to T$ be
the morphism above. Let $q_1,\ldots,q_r$ be the sections of
$\cH^{-1}Q$ chosen above. Under the preceding identification, each
$q_a$ determines a class in $H^{-1}\R\Gamma(W,Q)$. Since
$Q=\cofib(\phi)$, this class can be represented by a function
$f_a\in\Gamma(W,P^0)=\Gamma(W,\cO_W)$ together with a homotopy
from $\phi(f_a)$ to zero. More explicitly, in a complex computing
$\R\Gamma(W,T)$, the homotopy is represented by an element $h_a$
of degree $-1$ satisfying $\delta h_a=\phi(f_a)$.

On $W'$ adjoin generators $u_1,\ldots,u_r$ and define
\[
 P_1=(P|_{W'})[u_1,\ldots,u_r],\qquad
 |u_a|=-1,\qquad \delta u_a=f_a|_{W'}.
\]
The differential on $P_1$ extends that on $P|_{W'}$. The
homotopies we just chose, give an extension $\phi_1:P_1\to T|_{W'}$ of $\phi$,
with a homotopy between its restriction and the original
morphism, by the defining property of the analytic cell extension
\cite[Lemma 12.19]{LurieIX}.

\textit{Step 2.} Put $Q_1=\cofib(\phi_1)$ and
$R=\cofib(P|_{W'}\to P_1)$. As a complex, $R$ is represented by
the quotient $P_1/(P|_{W'})$, which is zero in degrees $0$ and
above. Each $u_a$ therefore defines a class in $\cH^{-1}R$.

Consider the distinguished triangle
\[
 R\longrightarrow Q|_{W'}\longrightarrow Q_1\longrightarrow R[1],
\]
in which the class of $u_a$ maps to $q_a|_{W'}$. These sections
generate $\cH^{-1}Q|_{W'}$, so the exact sequence
\[
 \cH^{-1}R\longrightarrow\cH^{-1}Q|_{W'}
 \longrightarrow\cH^{-1}Q_1\longrightarrow0
\]
shows that $\cH^{-1}Q_1=0$. From the same triangle,
$\cH^jQ_1=0$ for $j\geq0$. Thus $Q_1$ is $2$-connective.
In particular, $\cH^0(P_1)\to i_*\cO_Y|_{W'}$ is an
isomorphism, so the induced map on classical truncations is an
isomorphism.

\textit{Step 3.} Repeat with generators of the degree $-2$ cohomology of this
cofibre, chosen near $\overline{W''}$. Their lifts are closed
degree $-1$ elements with specified nullhomotopies. Adjoining degree
$-2$-cells kills the next relative cotangent classes, and produces a $3$-connective cofibre. 
\end{proof}

\begin{prop}\label{prop:diagram}\label{r:5.2}
There is a countable locally finite cover of $X$ by relative Darboux
charts and, for every nonempty finite intersection $X_J$, a finite
cell presentation $P_J$ over a Stein manifold $V_J$ smooth over $B$.
For $K\subset J$ there are maps
\begin{equation}\label{e:5.1}
 \Phi_{JK}:P_K\longrightarrow P_J,\qquad
 \Phi_{JL}=\Phi_{JK}\Phi_{KL}\quad(L\subset K\subset J),
\end{equation}
Each $\Phi_{JK}$ pulls back holomorphic coefficients along the product
projection $V_J\to V_K$ and sends every cell generator of $P_K$
to the corresponding generator of $P_J$.
The maps $\bX_J\to\Specan P_J$  identify classical truncations, and the maps
$P_J\to(i_J)_*\mathcal O_{\boldsymbol X_J}^{\mathrm{alg}}$
have $3$-connective cofibres.
\end{prop}
\begin{proof}

For a similar argument, see
\cite[Theorem~3.1 and \S4]{BJ}.

\emph{Step 1.}
Choose relative Darboux charts $A_i$, indexed by a countable set $I$,
with open chart domains $X_i^0\subseteq X$ and closed analytic
embeddings
\[
 \iota_i:X_i^0\hookrightarrow U_i^0,
\]
where $U_i^0$ is a relative polydisc and its map to $B$
is projection to the base coordinates.
Choose these charts so that the $X_i^0$ form a locally finite cover, and choose smaller relative polydiscs
$U_i'\Subset U_i^0$ such that the sets $\iota_i^{-1}(U_i')$
still cover $X$. Such choices can be made using a compact exhaustion
of $X$ and finitely many coordinate neighbourhoods at each stage.

Put
\[
 N_i=\#\{j\in I:\overline{X_i^0}\cap X_j^0\ne\varnothing\}.
\]
Each $N_i$ is finite. Insert intermediate relative polydiscs to obtain
$U_i^{2N_i}=U_i'\Subset\cdots\Subset U_i^0$, and define
\[
 X_i^a=\iota_i^{-1}(U_i^a),\qquad 0\leq a\leq2N_i.
\]
In particular, $X_i^a$ is open in $X$, and its image is a closed
analytic subspace of $U_i^a$.

For a finite nonempty $J\subseteq I$ with
$\bigcap_{i\in J}X_i^0\ne\varnothing$, put
\[
 W_J^a=\prod_{B,\,i\in J}U_i^a,
 \qquad X_J^a=\bigcap_{i\in J}X_i^a,
\]
whenever all the factors are defined. Write $\bX_J^a$ for the
corresponding open derived subspace of $\bX$. The spaces $W_J^a$
are smooth over $B$ and Stein, being closed analytic subspaces of
products of polydiscs. The map
$x\mapsto(\iota_i(x))_{i\in J}$ embeds $X_J^a$ as a closed
analytic subspace of $W_J^a$, by the closed diagonal of $X$.
Moreover, $|J|\leq N_i$ for every $i\in J$, so all levels
$0\leq a\leq2|J|$ are available.

\medskip
\noindent\emph{Step 2.}
We construct sheaves of cdgas $P_J$ by induction on $|J|$.
For $J=\{i\}$, define $P_{\{i\}}=A_i|_{U_i^2}$, with its
given differential and map from $\bX_{\{i\}}^2$, and let
$\mathcal{G}_{\{i\}}$ be its finite set of generators. 

Fix $k\geq2$. Suppose that for every nonempty $K$ with $|K|<k$
we have constructed $P_K$ on $W_K^{2|K|}$ and a finite set $\mathcal{G}_K$
of generators introduced at stage $K$, so that
\[
 P_K^\sharp=\cO_{W_K^{2|K|}}
 [\,e\in \mathcal{G}_L:\varnothing\ne L\subseteq K\,].
\]
Here $\sharp$ denotes the underlying graded algebra. 

We also suppose that the maps
$\bX_K^{2|K|}\to\Specan P_K$ have been chosen with coherent
homotopies for all chains of inclusions, on the relevant restrictions.

Fix $J$ with $|J|=k$, and put $W=W_J^{2k-2}$.
We first define $\widetilde P_J$ on $W$. Its underlying graded
algebra is
\[
 \widetilde P_J^\sharp
 =\cO_W[\,e\in \mathcal{G}_K:\varnothing\ne K\subsetneq J\,],
\]

For $e\in \mathcal{G}_K$, define $\widetilde\delta e$ by taking its
already defined differential $\delta_K e$ in $P_K$, and pulling
each holomorphic coefficient back along
$p_{JK}:W\to W_K^{2|K|}$.

Set $\widetilde\delta=0$ on $\cO_W$ and extend by the graded
Leibniz rule. This defines the sheaf of cdgas
$\widetilde P_J=(\widetilde P_J^\sharp,\widetilde\delta)$.

The defining property of the algebra $\widetilde P_J$  determines the map
$\bX_J^{2k-2}\to\Specan\widetilde P_J$.

We next construct $P_J$. Put $W'=W_J^{2k}$ and
$R=\widetilde P_J|_{W'}$. Apply Lemma~\ref{lem:cells} on
\[
 W_J^{2k}\Subset W_J^{2k-1}\Subset W_J^{2k-2}.
\]
Choose the cell extension supplied by this lemma and
write it explicitly as
\[
 P_J=\bigl(R^\sharp[u_1,\ldots,u_r,v_1,\ldots,v_s],\delta_J\bigr),
 \qquad |u_\alpha|=-1,\quad |v_\beta|=-2.
\]
The differential extends that of $R$ and is defined on the new
generators by
\[
 \delta_Ju_\alpha=f_\alpha,\qquad
 \delta_Jv_\beta=g_\beta.
\]
Here $f_\alpha\in\Gamma(W',\cO_{W'})$, and $g_\beta$ is a
degree $-1$ section of $R[u_1,\ldots,u_r]$, closed for the
differential already specified by $\delta_Ju_\alpha=f_\alpha$.
In particular, $\delta_J^2=0$. The lemma gives these data
together with a map
\[
 P_J\longrightarrow(i_J)_*\cO_{\bX_J^{2k}}^{\mathrm{alg}}
\]
whose cofibre is $3$-connective, where
$i_J:X_J^{2k}\hookrightarrow W'$ is the closed embedding.
Its restriction to $R$ agrees, up to the specified homotopy,
with the map induced above. We now set
$\mathcal{G}_J=\{u_1,\ldots,u_r,v_1,\ldots,v_s\}$.
\end{proof}

\begin{lem}\label{lem:combined-presentations}\label{r:5.3}
For $a=0,1$, let
\[
 \{P_J^a,\Phi_{JK}^a:\varnothing\ne K\subseteq J\subseteq I_a,
                         \ J\text{ finite}\}
\]
be a system of presentations satisfying Proposition~\ref{r:5.2}.
Use the notation $U_i^r,U_i',W_J^r$ from its proof, so that $P_J^a$
is defined on $W_J^{2|J|}$. Suppose that finite and pairwise disjoint
sets $\mathcal G_L^a$ of generators of degrees $-1,-2$ are given,
with
\[
 (P_J^a)^\sharp=
 \cO_{W_J^{2|J|}}
 [\,e\in\mathcal G_L^a:\varnothing\ne L\subseteq J\,].
\]
Require each $\Phi_{JK}^a$ to pull back holomorphic functions by
the product projection and to send every generator $e$ to the
same generator $e$. Consider
\[
 \epsilon_J^a:P_J^a\longrightarrow
 (i_J)_*\cO_{\bX_J^{2|J|}}^{\mathrm{alg}},
\]
 and fix homotopies
$\epsilon_J^a\circ\Phi_{JK}^a\simeq\epsilon_K^a$, together with
their compatibilities for every chain of inclusions.

Put $R_i=\iota_i^{-1}(U_i')$.
For each $a$, suppose that compact subsets $C_J^a\subseteq X$
are given such that
\[
 C_J^a\subseteq\bigcap_{i\in J}R_i,\qquad
 X=\bigcup_J C_J^a,
\]
and the family $\{C_J^a\}_J$ is locally finite.

Then there are relative polydiscs $\widehat U_i\Subset U_i'$
such that the sets $R_i''=\iota_i^{-1}(\widehat U_i)$,
$i\in I_a$, cover $X$ for each $a$. There are also presentations
$\widehat P_J$ on smaller product domains $\widehat V_J$, indexed
by finite nonempty subsets of $I_0\amalg I_1$, and maps
$\widehat\Phi_{JK},\widehat\epsilon_J$, satisfying the conclusions
of Proposition~\ref{r:5.2}, with the following additional properties.
\begin{enumerate}
\item[(i)] If $J\subseteq I_a$, then
\[
 \widehat P_J=P_J^a|_{\widehat V_J},\qquad
 \widehat\epsilon_J=\epsilon_J^a|_{\widehat V_J}.
\]
For $K\subseteq J\subseteq I_a$, the map $\widehat\Phi_{JK}$
is the restriction of $\Phi_{JK}^a$.  We also restrict to $\widehat V_J$ the homotopies $\epsilon_J^a\circ\Phi_{JK}^a\simeq\epsilon_K^a$
and their higher compatibilities.

\item[(ii)] If $J\subseteq I_a$, the open subset
\[
 \widehat X_J=
 \{x\in\bigcap_{i\in J}R_i:
                  (\iota_i(x))_{i\in J}\in\widehat V_J\}
\]
satisfies
\[
 C_J^a\subseteq\bigcap_{i\in J}R_i''\subseteq\widehat X_J.
\]
\end{enumerate}
We will omit below indices with $\bigcap_{i\in J}R_i=\varnothing$.
\end{lem}

\begin{proof}
\textit{Step 1.}
For $i\in I_a$, put
\[
 D_i=\bigcup_{J\ni i}C_J^a.
\]
Every nonempty set in this union meets the compact set
$\overline R_i$. Local finiteness implies that only finitely many
occur, so $D_i$ is a compact subset of $R_i$. Choose a relative
polydisc $\widehat U_i\Subset U_i'$ containing $\iota_i(D_i)$,
and put $R_i''=\iota_i^{-1}(\widehat U_i)$. Then
\[
 C_J^a\subseteq\bigcap_{i\in J}R_i''.
\]
In particular, the $R_i''$ with $i\in I_a$ cover $X$ for each $a$.

The family $\{R_i:i\in I_0\amalg I_1\}$ is locally finite, so
\[
 M_i=\#\{j\in I_0\amalg I_1:
                    \overline R_i\cap R_j\ne\varnothing\}
\]
is finite. Choose intermediate relative polydiscs
\[
 \widehat U_i=\widehat U_i^{2M_i}\Subset\cdots
       \Subset\widehat U_i^0=U_i'.
\]
If $\bigcap_{i\in J}R_i\ne\varnothing$, then $|J|\leqslant M_i$
for every $i\in J$. Consider the products
\[
 \widehat W_J^r=\prod_{B,\,i\in J}\widehat U_i^r
       \quad(0\leqslant r\leqslant2|J|),\qquad
 \widehat V_J=\widehat W_J^{2|J|}.
\]
hey are Stein manifolds smooth over $B$, and the
diagonal maps from the corresponding intersections in $X$ are
closed analytic embeddings, as in Proposition~\ref{r:5.2}.

\smallskip
\noindent\textit{Step 2.}
We construct the algebras $\widehat P_J$ and the maps
$\widehat\Phi_{JK}$ simultaneously by induction on $|J|$. If $J\subseteq I_a$,
define $\widehat P_J=P_J^a|_{\widehat V_J}$ and restrict $\epsilon_J^a$, the homotopies
$\epsilon_J^a\circ\Phi_{JK}^a\simeq\epsilon_K^a$, and their higher
compatibilities to $\widehat V_J$. This is possible since every factor of $\widehat V_J$
lies in $U_i'$, hence in the original domain of $P_J^a$.
Set $\mathcal G_J=\mathcal G_J^a$ for these subsets.

Now suppose that $J$ meets both $I_0$ and $I_1$, put $k=|J|$, and
assume the construction has been made for every proper subset of
$J$. On $W=\widehat W_J^{2k-2}$ define
\[
 \widetilde P_J^\sharp=
 \cO_W[\,e\in\mathcal G_K:\varnothing\ne K\subsetneq J\,].
\]
For $e\in\mathcal G_K$, define $\widetilde\delta e$ by pulling
back its differential in $\widehat P_K$ along the product
projection $W\to\widehat V_K$. Extend this differential by zero
on $\cO_W$ and by the graded Leibniz rule. The maps $\widehat\Phi_{KL}$, for
$\varnothing\ne L\subseteq K\subsetneq J$, commute with the
differentials by the induction hypothesis.

The pullbacks to $W$ of the algebras $\widehat P_K$, for
$\varnothing\ne K\subsetneq J$, form a diagram of $\cO_W$-algebras
with colimit $\widetilde P_J$. Since each $\widehat P_K$ is obtained
from $\operatorname{colim}_{\varnothing\ne L\subsetneq K}\widehat P_L$
by adjoining the cells $\mathcal G_K$, this colimit is also a
homotopy colimit (see also \cite[\S4]{BJ}).

By its universal property, the maps $\widehat\epsilon_K$,
together with the homotopies
$\widehat\epsilon_K\circ\widehat\Phi_{KL}\simeq\widehat\epsilon_L$
and their higher compatibilities, determine a morphism
\[
 \widetilde P_J\longrightarrow
 (\widehat i_J)_*\cO_{\widehat{\bX}_J^{2k-2}}^{\mathrm{alg}}.
\]
Here $\widehat{\bX}_J^r$ is the restriction of $\bX$ to
$\bigcap_{i\in J}\iota_i^{-1}(\widehat U_i^r)$, and
$\widehat i_J$ is the induced closed embedding of
$t_0(\widehat{\bX}_J^{2k-2})$ into $W$.

\smallskip
\noindent\textit{Step 3.}
For $K\subseteq J$, define $\widehat\Phi_{JK}$ by pullback of
holomorphic coefficients along the product projection and by
$e\mapsto e$ for every generator of $\widehat P_K$. The definition
of the differential makes this a morphism of cdgas. For
$L\subseteq K\subseteq J$, we have the equality
\[
 \widehat\Phi_{JL}
   =\widehat\Phi_{JK}\circ\widehat\Phi_{KL}.
\]
Indeed, both sides pull back coefficients by the same product
projection and send each generator of $\widehat P_L$ to itself.
If $J\subseteq I_a$, these maps agree with the restrictions of
the original maps by their specified action on coefficients and
generators. This proves (i).

If $J\subseteq I_a$, the cofibre of $\widehat\epsilon_J$ is
$\cofib(\epsilon_J^a)|_{\widehat V_J}$, which is $3$-connective.
Otherwise, $3$-connectivity follows from Lemma~\ref{r:5.1}. Finally, if $x\in C_J^a$ and $i\in J$, then
\[
 \iota_i(x)\in\widehat U_i
 \subseteq\widehat U_i^{2|J|}.
\]
Hence $x\in\widehat X_J$, proving (ii).
\end{proof}

The next statement is now straightforward from $3$-connectivity.

\begin{lem}\label{lem:form-lift}\label{r:5.4}
Let $g_J:\bX_J\to\Specan P_J$ be a comparison map from
Proposition~\ref{prop:diagram}. For $x\in X_J$, let $\C_x$ be
the residue field, regarded as a $P_J$-module through the map
$P_J\to\cO_{X_J}$. Define
\[
 \begin{gathered}
 D_{J,x}=\operatorname{Hom}_{\C}\bigl(
   \Omega^1_{\hol}(P_J/B)\Lotimes_{P_J}\C_x,\C\bigr),\\
 T_{J,x}=H^0(D_{J,x}),\quad O_{J,x}=H^1(D_{J,x}).
 \end{gathered}
\]
The complex $D_{J,x}$ represents the relative tangent complex of
$\Specan P_J\to B$ at the point corresponding to $x$.

 The map $g_J$ induces isomorphisms
\[
 T_x=H^0(\bT_{\bX/B}|_x)\xrightarrow{\sim}T_{J,x},
 \qquad
 O_x=H^1(\bT_{\bX/B}|_x)\xrightarrow{\sim}O_{J,x},
\]
where $|_x$ denotes the derived fibre. These identifications are
compatible with the transition maps $\Phi_{JK}$.

There exists a relative $2$-form
\[
 \eta_J\in\Gamma\bigl(V_J,\Omega^2_{\hol}(P_J/B)^{-2}\bigr),
 \qquad \delta\eta_J=0,
\]
where $\delta$ is the differential induced by that of $P_J$,
whose cohomology class maps to the underlying $2$-form class of
$\omega|_{\bX_J}$. The symmetric bilinear form induced by
$\eta_J$ on $O_{J,x}$ agrees with $Q_x$ under the 
isomorphism $O_x\cong O_{J,x}$.
\end{lem}

\section{Construction of the Kuranishi atlas}\label{sec:reduction}\label{sec:6}
Using the presentations of
Proposition~\ref{r:5.2}, we now construct the relative Kuranishi atlas. We follow here
\cite[\S\S3.3--3.4 and \S5]{BJ} and
\cite[\S6.2]{BJ}, and  use the notation of the proof of
Proposition~\ref{r:5.2}, so that
\[
 V_J=W_J^{2|J|},\qquad X_J=X_J^{2|J|}.
\]
When choosing $E_J^-$ and forming $E_J^+$, we regard
$E_J$ as a real vector bundle.

\subsection{Admissible subbundles}\label{ss:6.1}
Write $u_\alpha,v_\beta$ for the generators of $P_J$ of degrees $-1,-2$,
and $r_J,q_J$ for the numbers of these respective generators. Their differentials have the form
\[
 \delta u_\alpha=f_\alpha,\qquad
 \delta v_\beta=\sum_{\alpha=1}^{r_J}g_{\beta\alpha}u_\alpha,
\]
where $f_\alpha,g_{\beta\alpha}$ are holomorphic functions on $V_J$.
Let $E_J=V_J\times\C^{r_J}$ and $F_J=V_J\times\C^{q_J}$
be the trivial holomorphic vector bundles with standard bases
$e_\alpha$ and $h_\beta$.

Define a holomorphic
section $s_J$ of $E_J$ and a holomorphic bundle morphism
$t_J:E_J\to F_J$ by
\begin{equation}\label{e:6.1}
 s_J=\sum_\alpha f_\alpha e_\alpha,\qquad
 t_J(e_\alpha)=\sum_\beta g_{\beta\alpha}h_\beta,\qquad t_J(s_J)=0.
\end{equation}
The last identity follows from $\delta^2v_\beta=0$.

Proposition~\ref{r:5.2} identifies $s_J^{-1}(0)$ with
$X_J$, and we will write $\psi_J:s_J^{-1}(0)\to X_J$ for this identification.
At a zero $v$ corresponding to $x=\psi_J(v)$, differentiation along
the fibres of $V_J\to B$ results in a complex
\[
 T_v(V_J/B)\xrightarrow{D_{J,v}}E_{J,v}
 \xrightarrow{t_{J,v}}F_{J,v},\qquad
 D_{J,v}=(\mathrm d_{V_J/B}s_J)_v,
\]
in degrees $0,1,2$. By Lemma~\ref{r:5.4}, its cohomology in degrees
$0$ and $1$ is identified with the relative tangent and obstruction spaces:
\[
 \ker D_{J,v}\cong T_x=H^0(\bT_{\bX/B}|_x),\qquad
 \ker t_{J,v}/\im D_{J,v}\cong O_x=H^1(\bT_{\bX/B}|_x).
\]
For $K\subsetneq J$, by \eqref{e:5.1}, we may consider holomorphic
submersions $\phi_{JK}:V_J\to V_K$ over $B$ and surjective bundle morphisms
$\chi_{JK}:E_J\to\phi_{JK}^*E_K$,
$\xi_{JK}:F_J\to\phi_{JK}^*F_K$. They satisfy
\[
 \chi_{JK}s_J=\phi_{JK}^*s_K,\qquad
 \xi_{JK}t_J=(\phi_{JK}^*t_K)\chi_{JK},
\]
and induce morphisms between the complexes above at corresponding zeros.

 Moreover, for $L\subsetneq K\subsetneq J$,
\begin{equation}\label{e:6.2}
 \begin{aligned}
 \phi_{JL}&=\phi_{KL}\circ\phi_{JK},\\
 \chi_{JL}&=(\phi_{JK}^*\chi_{KL})\circ\chi_{JK},\\
 \xi_{JL}&=(\phi_{JK}^*\xi_{KL})\circ\xi_{JK}.
 \end{aligned}
\end{equation}

The induced maps on $H^0,H^1$ are the
identifications with $T_x,O_x$ from Lemma~\ref{r:5.4}, and, in particular,
they preserve $Q_x$.

\begin{dfn}\label{def:negative}\label{r:6.1}
Fix one of these presentations, and omit the subscript $J$.
Let $U\subset V$ be open. A smooth real vector subbundle
$E^-\subset E|_U$ is called \emph{admissible} if, for every
$v\in s^{-1}(0)\cap U$, with $x=\psi(v)$ and
$D_v=(\mathrm d_{V/B}s)_v$, the following conditions hold:
\begin{enumerate}
\renewcommand{\labelenumi}{(\roman{enumi})}
\item $E^-_v\cap\im D_v=0$;
\item $t_v(E^-_v)=\im t_v$;
\item the image $N_x$ of $E^-_v\cap\ker t_v$ in
$O_x=\ker t_v/\im D_v$ is a maximal negative real subspace
for $\operatorname{Re}Q_x$.
\end{enumerate}
Here $E$ and $F$ are regarded as real vector bundles.
Let $E^+=E|_U/E^-$ be the quotient real vector bundle, and let
$s^+$ be the image of $s|_U$ in $E^+$. (See \cite[Definition 3.21]{BJ}.)
\end{dfn}
At a zero $v$, put $n=\dim_{\C}(V/B)$, $r=\rk_{\C}E$,
$a=\rk_{\C}D_v$ and $c=\rk_{\C}t_v$.
Then $\dim_{\C}T_x=n-a$ and $\dim_{\C}O_x=r-a-c$. 

The definition and \eqref{e:2.5} give
\begin{equation}\label{e:6.3}
 \begin{gathered}
 d=2n-r-a+c,\qquad \rk_{\RR}E^-=r+c-a,\\
 \dim_{\RR}(U/B)-\rk_{\RR}E^+=d.
 \end{gathered}
\end{equation}
In particular, the required rank of $E^-$ is $2r-2n+d$, which
is locally constant.

\begin{lem}\label{lem:admissible-extension}\label{r:6.2}
The conditions of Definition~\ref{r:6.1} are open along $s^{-1}(0)$
for a fixed smooth subbundle of the required rank. If $C\subset s^{-1}(0)$
is compact, there is an admissible pair $(U,E^-)$ with $C\subset U$.
If an admissible subbundle is already given on an open neighbourhood
of a closed subset $D\subset C$, the new subbundle may be chosen to agree
with it on an open neighbourhood of $D$.
\end{lem}
\begin{proof}
Put $Z=s^{-1}(0)$. We first prove openness, following \cite[\S5.1]{BJ}.
Let $E^-\subset E|_U$ be a smooth subbundle of rank $2r-2n+d$.
Suppose $v_j\to v$ in $Z\cap U$, with $E^-_v$ admissible and each
$E^-_{v_j}$ not admissible. After passing to a subsequence, we may
suppose that the ranks of $D_{v_j}$ and $t_{v_j}$ are constant, say
$a',c'$, and that
\[
 I=\lim_j\im D_{v_j},\qquad K=\lim_j\ker t_{v_j}
\]
exist in the appropriate Grassmannians. Put $A=\im D_v$, and write
$a,c$ for the ranks at $v$. By continuity,
$A\subset I\subset K\subset\ker t_v$.

The representative in Lemma~\ref{r:5.4} defines a symmetric form
$\beta$ on $E$ whose restriction to $\ker t_v$ induces $Q_x$.
Its closedness for $\delta$ gives
$\beta_w(\im D_w,\ker t_w)=0$ at every zero $w$.
Taking limits yields $\beta_v(I,K)=0$, and hence
\[
 I/A\subseteq(K/A)^{\perp_{Q_x}}\subseteq O_x.
\]
The spaces $I/A$ and $(K/A)^{\perp_{Q_x}}$ have dimensions
$a'-a$ and $c'-c$. These agree by \eqref{e:6.3}, so the inclusion
is an equality. In particular, $I/A$ is isotropic.

Let $L=E^-_v$ and let $N\subset O_x$ be the image of
$L\cap\ker t_v$. Negativity implies $N\cap(I/A)=0$.
Also $N+(K/A)=O_x$ as real spaces: its orthogonal complement is
$N^\perp\cap(I/A)$, which is zero since $N^\perp$ is positive
definite. Thus
\[
 L\cap I=0,\qquad L+K=E_v.
\]
For large $j$ this implies $E^-_{v_j}\cap\im D_{v_j}=0$ and
$E^-_{v_j}+\ker t_{v_j}=E_{v_j}$. The spaces
$E^-_{v_j}\cap\ker t_{v_j}$ tend to $L\cap K$, on which
$\operatorname{Re}\beta_v$ is negative definite. Thus
$\operatorname{Re}\beta_{v_j}$ is negative definite on
$E^-_{v_j}\cap\ker t_{v_j}$ for all sufficiently large $j$,
and, by \eqref{e:6.3},
\[
 \dim_{\RR}(E^-_{v_j}\cap\ker t_{v_j})
 =r-a'-c'=\dim_{\C}O_{\psi(v_j)}.
\]
All three conditions hold for large $j$, a contradiction.

We now prove the extension assertion, as in \cite[\S5.2]{BJ}.
Choose a relatively compact open neighbourhood $V_0$ of $C$ in $V$.
After making $V_0$ smaller, choose a finite stratification of $Z\cap V_0$
by smooth complex submanifolds on which the ranks of $D$ and $t$ are
constant, such that the frontier of each stratum is a union of strata
of smaller dimension. Such a stratification exists for complex analytic
sets and their analytic rank loci.

See also \cite[\S18]{Whitney}.

For each $m$, let $S_m$ be the union
of the strata of complex dimension $m$, and put
$Z_m=\bigcup_{j\le m}S_j$. Thus $Z_m$ is closed in $V_0$, and $S_m$
is a closed smooth submanifold of $V_0\setminus Z_{m-1}$, possibly
disconnected. We construct an admissible subbundle near $D\cup Z_m$
by induction on $m$, leaving it unchanged on a smaller neighbourhood of
$D\cup Z_{m-1}$ at each step.

Suppose $P\subset E|_O$ has already been chosen on an open
neighbourhood $O$ of $A=D\cup Z_{m-1}$, and write $S=S_m$.
Over $S$, put $I=\im D$ and $K=\ker t$. Choose a maximal negative
real subbundle $N\subset K/I$ for $\operatorname{Re}Q$, a subbundle
$\widetilde N\subset K$ mapping isomorphically to $N$, and a
complement $L$ to $K$ in $E|_S$. Then
$P_S=\widetilde N\oplus L$ is admissible on $S$.

On $S\cap O$, define $N$ to be the image of $P\cap K$ in $K/I$,
take $\widetilde N=P\cap K$, and choose $L$ inside $P$.
The spaces of choices at each of these three steps are contractible.
After replacing $O$ by a smaller neighbourhood of $A$, we can therefore
extend these choices over $S$, with $P_S=P$ on $S\cap O$.

Now extend $P_S$ to a smooth subbundle near $S\cup A$, agreeing with $P$
on a smaller neighbourhood of $A$. By the openness proved above,
this extension is admissible (after maybe restricting its domain).

Starting with the given subbundle near $D$, or the empty open set if
$D=\varnothing$, this proves the induction, and the lemma follows.
\end{proof}

\subsection{Kuranishi neighbourhoods and coordinate changes}\label{ss:6.2}
\begin{lem}\label{lem:compatible-complements}\label{r:6.3}
Let $K\subsetneq J$, and let $(U_K,E_K^-)$ be admissible.
At corresponding zeros $v,w$, put $D_J=(\mathrm d_{V_J/B}s_J)_v$
and define $D_K$ similarly. Then $D_J$ is injective on
$\ker\mathrm d\phi_{JK}$, and
\begin{equation}\label{e:6.4}
 \ker\chi_{JK}\cap\ker t_J=D_J(\ker\mathrm d\phi_{JK}),\qquad
 t_J(\ker\chi_{JK})=\ker\xi_{JK}\cap\im t_J.
\end{equation}
On $\phi_{JK}^{-1}(U_K)$ put
\[
 P_{JK}=\chi_{JK}^{-1}(\phi_{JK}^*E_K^-),\qquad
 R_{JK}=(\mathrm ds_J)(\ker\mathrm d\phi_{JK}).
\]

Near the zero locus, $R_{JK}$ is a subbundle of $P_{JK}$.
A subbundle $E_J^-\subset P_{JK}$ is admissible if and only if,
after restriction to a neighbourhood of the zero locus,
\begin{equation}\label{e:6.5}
 P_{JK}=E_J^-\oplus R_{JK}.
\end{equation}
In particular, admissible subbundles $E_J^-\subset P_{JK}$
exist near the zero locus. For compact sets
$D\subseteq C\subset s_J^{-1}(0)\cap\phi_{JK}^{-1}(U_K)$,
any such subbundle defined near $D$ extends to an admissible
subbundle near $C$, agreeing with the given subbundle on
an open neighbourhood of $D$.
\end{lem}
\begin{proof}
At a zero, the morphism of the three-term tangent complexes is
surjective in each degree and induces isomorphisms on $H^0,H^1$, hebce 
$D_J$  is injective on $\ker\mathrm d\phi_{JK}$, and the first identity in
\eqref{e:6.4} follows.

For the second, suppose $\xi_{JK}t_J(e)=0$.
Then $\chi_{JK}(e)\in\ker t_K$. Surjectivity on $H^1$, followed
by surjectivity of $\mathrm d\phi_{JK}$, gives $e'\in\ker t_J$
with $\chi_{JK}(e')=\chi_{JK}(e)$. Therefore
$e-e'\in\ker\chi_{JK}$ and $t_J(e-e')=t_J(e)$, as required.

Now, note that $R_{JK}\subset\ker\chi_{JK}\subset P_{JK}$
on $\phi_{JK}^{-1}(U_K)$. At a zero, $\mathrm ds_J$ is injective
on $\ker\mathrm d\phi_{JK}$ by \eqref{e:6.4}, and hence $R_{JK}$
is a subbundle near the zero locus.

We prove the remaining assertions at a zero
$v\in s_J^{-1}(0)\cap\phi_{JK}^{-1}(U_K)$, with
$w=\phi_{JK}(v)$ and $x=\psi_J(v)=\psi_K(w)$.

Put $A_J=\im D_J$ and $A_K=\im D_K$. Then
\begin{equation}\label{e:6.6}
 P_{JK}\cap A_J=R_{JK},\qquad
 \ker\chi_{JK}\cap\ker t_J=R_{JK},\qquad
 t_J(P_{JK})=\im t_J.
\end{equation}
For the first identity, if $D_Ju\in P_{JK}$, then
$D_K\mathrm d\phi_{JK}(u)\in E_K^-\cap A_K=0$.
The isomorphism on $H^0$ gives $u_0\in\ker D_J$ with
$\mathrm d\phi_{JK}(u_0)=\mathrm d\phi_{JK}(u)$, so
$D_Ju=D_J(u-u_0)\in R_{JK}$. The second identity is
\eqref{e:6.4}. For the third, take $e\in E_J$.
Since $t_K(E_K^-)=\im t_K$, there is $e_K^-\in E_K^-$ with
$t_K(e_K^-)=t_K\chi_{JK}(e)$. Lift $e_K^-$ to $p\in P_{JK}$.
Then $t_J(e-p)\in\ker\xi_{JK}\cap\im t_J=t_J(\ker\chi_{JK})$
by \eqref{e:6.4}. Adding a suitable element of $\ker\chi_{JK}$
to $p$ proves the third identity.

It follows that any complement in \eqref{e:6.5} is disjoint from
$A_J$ and maps onto $\im t_J$ under $t_J$. Furthermore,
$\chi_{JK}$ induces an isomorphism
\[
 E_J^-\cap\ker t_J\longrightarrow E_K^-\cap\ker t_K.
\]
Indeed, its kernel is zero by \eqref{e:6.6}. For surjectivity,
lift an element of $E_K^-\cap\ker t_K$ to $p\in P_{JK}$.
Since $t_J(p)\in\ker\xi_{JK}\cap\im t_J=t_J(\ker\chi_{JK})$,
subtract an element of $\ker\chi_{JK}$ with the same $t_J$ image,
and then project along $R_{JK}$. This isomorphism identifies their
images in the common obstruction space $O_x$. Hence the image is
negative definite and has real dimension $\dim_{\C}O_x$.

Conversely, an admissible $E_J^-\subset P_{JK}$ is disjoint from
$R_{JK}$, and, by \eqref{e:6.3},
\[
 \rk_{\RR}P_{JK}-\rk_{\RR}R_{JK}
 =2r_J-2\dim_{\C}(V_J/B)+d=\rk_{\RR}E_J^-.
\]
Thus \eqref{e:6.5} holds at each zero and, by openness of
transversality, on a neighbourhood of the zero locus.

Choose a smooth complement $H_0$ to $R_{JK}$ in $P_{JK}$ near $C$. A complement given near $D$ is then the graph of a bundle map
$a:H_0\to R_{JK}$.

To construct such object, choose a smooth function $\rho$, equal to one near $D$ and
supported where $a$ is defined, and extend $\rho a$ by zero near
$C$: its graph is the required extension.
\end{proof}

\begin{lem}\label{lem:reduction}\label{r:6.4}
The following statements hold.
\begin{enumerate}
\renewcommand{\labelenumi}{(\roman{enumi})}
\item For an admissible pair $(U,E^-)$ there is an open
$U'\subseteq U$ containing $s^{-1}(0)\cap U$ such that,
after restriction to $U'$,
\begin{equation}\label{e:6.7}
 (s^+)^{-1}(0)=s^{-1}(0)\cap U'.
\end{equation}
Thus $(U',E^+,s^+,\psi|)$ is a Kuranishi neighbourhood of
relative virtual dimension $d$ over $B_{\RR}$.
\item Suppose $K\subsetneq J$ and admissible pairs
$(U_J,E_J^-),(U_K,E_K^-)$ satisfy \eqref{e:6.7}.
On an open $U_{JK}\subseteq U_J$, suppose
\[
 \phi_{JK}(U_{JK})\subseteq U_K,\qquad
 \chi_{JK}(E_J^-|_{U_{JK}})\subseteq\phi_{JK}^*E_K^-.
\]
Then $\chi_{JK}$ induces a quotient map $\chi_{JK}^+$, and
$(U_{JK},\phi_{JK}|_{U_{JK}},\chi_{JK}^+)$ is a coordinate
change over $\psi_J(s_J^{-1}(0)\cap U_{JK})$.
\end{enumerate}
\end{lem}
\begin{proof}
We first prove (i). For a relative Darboux presentation minimal
at a zero $v$, let $q$ be the nondegenerate symmetric
form on $E$ given by the Darboux expression. We have $D_v=t_v=0$.
Admissibility says that $\operatorname{Re}q_v$ is negative definite
on $E_v^-$, and this remains true on $E^-$ near $v$. The identity
$q(s,s)=0$ of Theorem~\ref{r:4.2} holds on $V$. If $s^+=0$, then $s\in E^-$, and hence $s=0$.

We next show how to transfer the equality \eqref{e:6.7} between
two presentations. Write
$\phi:V_J\to V_K$, $\chi:E_J\to\phi^*E_K$ and
$\xi:F_J\to\phi^*F_K$ for a comparison as in
Lemma~\ref{r:6.3}. Fix a zero $v$ and put $w=\phi(v)$.
First let $E_J^-$ be an admissible subbundle near $v$, we shall
construct an admissible subbundle $E_K^-$ near $w$. At $v$ put
$G=t_J(\ker\chi)$, extend it to a local subbundle of $F_J$,
and choose a projection $\pi:F_J\to G$. The restrictions
of $\pi t_J$ to $\ker\chi$ and to $E_J^-$ are surjective
at $v$, hence nearby. Consequently
\[
 H=\ker(\pi t_J),\qquad R=H\cap\ker\chi,\qquad L=H\cap E_J^-
\]
are smooth real subbundles near $v$. The map
$\chi:H\to\phi^*E_K$ is surjective with kernel $R$.
At $v$, \eqref{e:6.4} gives $R_v=D_J(\ker\mathrm d\phi)$
and $L_v\cap R_v=0$. Choose a decomposition $H=E_0\oplus R$
near $v$ with $L\subset E_0$. Then $\chi|_{E_0}$ is an isomorphism.

Since $t_Js_J=0$, write $s_J=s_0+s_R$ in this decomposition.
The derivative of $s_R$ on $\ker\mathrm d\phi$ is an
isomorphism at $v$. By the implicit function theorem, a
submanifold $S=s_R^{-1}(0)$ near $v$ on which $\phi$ is a
diffeomorphism onto a neighbourhood of $w$: we write $\theta$
for its inverse, and on that
neighbourhood define $E_K^-=\theta^*\chi(L)$.

We check admissibility at $w$. If $l\in L_v$ and
$\chi(l)\in\im D_K$, lift its preimage along $\mathrm d\phi$.
Subtracting the corresponding vector in $\im D_J$ from $l$
leaves a vector in $\ker\chi\cap H_v=R_v$, since
$\pi t_J(l)=0$ and $t_JD_J=0$. Thus $l\in\im D_J$,
and $l=0$ by admissibility of $E_J^-$. 

Now, for any $e_K\in E_K$,
choose a lift $h\in H_v$ and then $l\in E_J^-$ with
$t_J(l)=t_J(h)$. Since $\pi t_J(h)=0$, we have $l\in L_v$.
It follows that $t_K(E_K^-)=\im t_K$.
Finally, if $l\in L_v$ and $t_K\chi(l)=0$, then
$t_J(l)\in G\cap\ker\pi=0$.

Hence $\chi$ identifies
$L_v\cap\ker t_J=E_{J,v}^-\cap\ker t_J$ with
$E_{K,w}^-\cap\ker t_K$, and their images in the common
obstruction space agree. This proves negativity, the required
dimension, and (together with Lemma~\ref{r:6.2}) proves admissibility near $w$.

As $s_J\in H$ and $H\cap E_J^-=L\subseteq E_0$,
\[
 s_J\in E_J^-\quad\Longleftrightarrow\quad s_R=0\text{ and }s_0\in L.
\]
On $S$, the isomorphism $\chi|_{E_0}$ takes $s_0$ to $s_K$.
Thus $\phi|_S$ identifies both the original zero loci and the
zero loci of the quotient sections near $v,w$.

If instead an admissible subbundle $E_K^-$ near $w$ is given,
choose a compatible $E_J^-$ near $v$ by Lemma~\ref{r:6.3}.
The preceding construction recovers the original $E_K^-$.
Indeed, compatibility gives $\chi(L)\subseteq\phi^*E_K^-$,
and these subbundles have the same rank: if
$h=\dim_{\C}t_J(\ker\chi)$ at $v$, then
\[
 \rk_{\RR}L=\rk_{\RR}E_J^--2h=\rk_{\RR}E_K^-,
\]
by \eqref{e:6.4} and \eqref{e:6.3}.
Thus $\chi(L)=\phi^*E_K^-$ near $v$, and
$\theta^*\chi(L)=E_K^-$ near $w$.

Indeed, $\chi(L)$ is contained
in it and the ranks agree: if $h=\dim_{\C}t_J(\ker\chi)$ at $v$, then
\[
 \rk_{\RR}L=\rk_{\RR}E_J^--2h=\rk_{\RR}E_K^-.
\]
This follows from \eqref{e:6.4} and \eqref{e:6.3}.
Let $P$ be the cell presentation on $V$ defining $E$ and $s$,
and choose a relative Darboux presentation $P'$ on $V'$,
minimal at the same point $x\in X$. For the corresponding zeros $v,v'$, choose Stein open
neighbourhoods $W\Subset W_1\Subset W_0\subset V\times_B V'$
of $(v,v')$. The two inclusions of zero loci define a closed
embedding $i:Y\hookrightarrow W_0$ of a common open neighbourhood
$Y\subset X$ of $x$. Let $Q$ be the tensor product over
$\cO_{W_0}$ of the pullbacks of $P,P'$. Their maps to the
structure sheaf determine $Q\to i_*\cO_{\bX|_Y}^{\mathrm{alg}}$.

By Lemma~\ref{r:5.1}, we can factor the morphism on $W$ as
\(
 Q|_W\longrightarrow\widehat P
 \longrightarrow
 \left(i_*\cO_{\bX|_Y}^{\mathrm{alg}}\right)|_W,
\)
where the first arrow is a cell extension and the cofibre of the
second is concentrated in degrees at most $-3$. As in
Lemma~\ref{r:5.4}, the comparisons with $P$ and $P'$ identify
the relative tangent and obstruction spaces at corresponding
zeros, preserving the quadratic forms on the obstruction spaces.

By Lemma~\ref{r:6.3}, choose an admissible subbundle
$\widehat E^-$ for $\widehat P$ whose image under the comparison
with $P$ lies in the pullback of $E^-$. Apply the preceding
construction to the comparison with $P'$ to obtain an admissible
subbundle $E'^-$. The first paragraph proves \eqref{e:6.7} for
$P'$ and $E'^-$, and then it follows for
$\widehat P$ and $\widehat E^-$, and hence for $P$ and $E^-$
near $v$. The union of these neighbourhoods in $U$, for all
$v\in s^{-1}(0)\cap U$, is the required $U'$.
The dimension follows from \eqref{e:6.3}.

For (ii), let $x=\psi_J(v)=\psi_K(w)$ be a point of
$\psi_J(s_J^{-1}(0)\cap U_{JK})$. The first two admissibility
conditions give
\begin{equation}\label{e:6.8}
 \ker D_J^+\cong T_{x,\RR},\qquad
 \operatorname{coker}D_J^+\cong O_{x,\RR}/N_{J,x},\qquad
 D_J^+=(\mathrm d_{V_J/B}s_J^+)_v.
\end{equation}
The inclusion of negative subbundles identifies $N_{J,x}$
with a subspace of $N_{K,x}$, and their dimensions are equal,
so they are the same subspace of $O_{x,\RR}$. Thus the
comparison of the two complexes $[T(V/B)\to E^+]$ induces
isomorphisms on both cohomology groups. Thus, one has
\begin{equation}\label{e:6.9}
 0\longrightarrow T_v(V_J/B)
 \xrightarrow{(D_J^+,\mathrm d\phi_{JK})}
 E_{J,v}^+\oplus T_w(V_K/B)
 \xrightarrow{-\chi_{JK}^++D_K^+}E_{K,w}^+
 \longrightarrow0.
\end{equation}
Passing to total tangent spaces preserves exactness, since the
quotient map on $T_{f(x)}B_{\RR}$ is the identity. Together
with $\chi_{JK}^+s_J^+=\phi_{JK}^*s_K^+$ and the identity
$\psi_K\circ\phi_{JK}=\psi_J$ on the zero locus, this produces a
1-morphism satisfying the exactness condition \eqref{e:3.4}.
It is therefore a coordinate change by Theorem~\ref{r:3.3}.
\end{proof}

\begingroup
\color{black}
\subsection{Compatible choices and cocycles}\label{ss:6.3}

\begin{prop}\label{r:6.5}\label{prop:compatible-reductions}
The presentations of Proposition~\ref{r:5.2} admit admissible pairs
$(U_J,E_J^-)$ whose Kuranishi neighbourhoods form a relative Kuranishi atlas on $X$.
\end{prop}

\begin{proof}
We follow \cite[\S6.2]{BJ}. Write $R_i=X_i^{2N_i}=\iota_i^{-1}(U_i')$
for the innermost open subsets in the proof of Proposition~\ref{r:5.2}.
They cover $X$, and each has compact closure in $X_i^0$.
Apply Proposition~\ref{r:3.8} to this cover. 

We have a locally finite
closed cover $\{C_J\}$ with
\[
 C_J\subset R_J:=\bigcap_{i\in J}R_i\subset X_J,\qquad
 C_J\cap C_K\ne\varnothing\quad\Longrightarrow\quad
 J\subseteq K\text{ or }K\subseteq J.
\]

\noindent\emph{Step 1.} We first construct, by induction on $|J|$,
an admissible pair $(\widetilde U_J,\widetilde E_J^-)$
near $\psi_J^{-1}(C_J)$. For every $K\subsetneq J$ we also require
an open neighbourhood $\widetilde U_{JK}$ of $\psi_J^{-1}(C_J\cap C_K)$ such that
\begin{equation}\label{e:6.10}
 \widetilde U_{JK}\subset\widetilde U_J,\qquad
 \phi_{JK}(\widetilde U_{JK})\subset\widetilde U_K,\qquad
 \chi_{JK}(\widetilde E_J^-|_{\widetilde U_{JK}})
     \subset\phi_{JK}^*\widetilde E_K^-.
\end{equation}
For $|J|=1$, Lemma~\ref{r:6.2} gives the required pair immediately.

Suppose now the choices have been made for $|K|<|J|$, and fix $J$.
For each proper subset $K$ occurring in the cover, put
\[
 A_K=\psi_J^{-1}(C_J\cap C_K),\qquad
 P_K=\chi_{JK}^{-1}(\phi_{JK}^*\widetilde E_K^-),\qquad
 R_K=(\mathrm ds_J)(\ker\mathrm d\phi_{JK}).
\]
The last two bundles are defined near $A_K$, after restriction as in
Lemma~\ref{r:6.3}. If $L\subsetneq K\subsetneq J$, by the induction hypothesis,

\begin{equation}\label{e:6.11}
 P_K\subset P_L\quad\text{near }A_K\cap A_L,\qquad R_K\subset R_L
\end{equation}
on the domains in question. (The second inclusion follows from
$\ker\mathrm d\phi_{JK}\subset\ker\mathrm d\phi_{JL}$.)

By descending induction on $|K|$, choose complements $P_K=H_K\oplus R_K$
near $A_K$, agreeing near $A_K\cap A_L$ whenever $K,L$ are comparable.
For maximal $K$, choose any complement. At the step for $L$, the previously
chosen $H_K$, for $L\subsetneq K\subsetneq J$, lie in $P_L$ near $A_K\cap A_L$
by \eqref{e:6.11}. They are admissible, so Lemma~\ref{r:6.3} shows that they
complement $R_L$ there (maybe after restriction.)

These complements agree near intersections of their compact sets.
Indeed, incomparable $K_1,K_2$ have $A_{K_1}\cap A_{K_2}=\varnothing$,
and comparable indices have already been treated by the descending induction.
As only finitely many subsets of $J$ occur, their neighbourhoods may be chosen
so that our complements agree on every overlap. They then define
one complement near the compact set
$\bigcup_{L\subsetneq K\subsetneq J}(A_K\cap A_L)$.
Lemma~\ref{r:6.3} extends it to a complement near
$A_L$.
The descending induction is proven.

Choose smaller neighbourhoods of the finitely many $A_K$ so that the $H_K$
agree on all overlaps. They give an admissible subbundle near the compact set
$\bigcup_{K\subsetneq J}A_K$. Lemma~\ref{r:6.2}, with
$C=\psi_J^{-1}(C_J)$ and $D=\bigcup_{K\subsetneq J}A_K$, extends it to
an admissible $\widetilde E_J^-$ near $C$, and
\eqref{e:6.10} holds. This completes the induction on $|J|$.

\noindent\emph{Step 2.} We now choose the open sets $U_J$ and $U_{JK}$
on which the Kuranishi neighbourhoods and their coordinate changes
are defined. First use Lemma~\ref{r:6.4}(i) near each $\psi_J^{-1}(C_J)$,
so that $(s_J^+)^{-1}(0)=s_J^{-1}(0)$ on the resulting domains.
The images of these zero loci under $\psi_J$ are open neighbourhoods of $C_J$ in $X$.

Choose a locally finite family of relatively compact open neighbourhoods $O_J$. We can make these neighbourhoods smaller so that $O_J\cap O_K=\varnothing$
for incomparable $J,K$, and
\begin{equation}\label{incl}
 O_J\cap O_K\subset
 \psi_J\bigl(s_J^{-1}(0)\cap\widetilde U_{JK}\bigr)
 \quad(K\subsetneq J).
\end{equation}
These restrictions are possible because the respective right hand sides contain
$C_J\cap C_K$. If two compact sets have intersection in an open set $W$, they have
open neighbourhoods whose intersection lies in $W$. Apply this to each pair for which
the original $O_J,O_K$ meet. Each index occurs in only finitely many such pairs.
Intersecting the corresponding neighbourhoods, we get open sets
$C_J\subset S_J\subset O_J$ such that $S_J\cap S_K=\varnothing$
when neither index contains the other, and the inclusion \ref{incl} holds with $S_J,S_K$ in place of $O_J,O_K$. 

Choose $U_J$ inside the domains where the equality of zero loci was established,
with $\psi_J(s_J^{-1}(0)\cap U_J)=S_J$, and put
$E_J^-=\widetilde E_J^-|_{U_J}$. For $K\subsetneq J$ define
\[
 U_{JK}=U_J\cap\widetilde U_{JK}\cap\phi_{JK}^{-1}(U_K).
\]
Then
\[
 s_J^{-1}(0)\cap U_{JK}=\psi_J^{-1}(S_J\cap S_K).
\]
The pairs are compatible on $U_{JK}$ by \eqref{e:6.10}, so
Lemma~\ref{r:6.4}(ii) gives the coordinate changes
$\Phi_{JK}^+=(U_{JK},\phi_{JK}|_{U_{JK}},\chi_{JK}^+)$.

\noindent\emph{Step 3.} It remains to check the cocycle condition.
Let $L\subsetneq K\subsetneq J$ and put
\[
 U_{JKL}=U_{JK}\cap U_{JL}\cap\phi_{JK}^{-1}(U_{KL}).
\]
This open set contains $\psi_J^{-1}(S_J\cap S_K\cap S_L)$.
On $U_{JKL}$ the maps $\phi_{JK}$ satisfy $\phi_{JL}=\phi_{KL}\circ\phi_{JK}$.
If $q_J:E_J\to E_J^+$ is the quotient map, then
\[
\begin{aligned}
 \chi_{JL}^+q_J&=\phi_{JL}^*q_L\,\chi_{JL}\\
 &=\phi_{JK}^*(\chi_{KL}^+)\,\chi_{JK}^+q_J.
\end{aligned}
\]
As $q_J$ is surjective,
\begin{equation}\label{e:6.12}
 \chi_{JL}^+=\phi_{JK}^*(\chi_{KL}^+)\circ\chi_{JK}^+
 \quad\text{on }U_{JKL}.
\end{equation}
The zero bundle map on $U_{JKL}$ therefore represents a $2$-isomorphism
$\Lambda_{JKL}:\Phi_{KL}^+\circ\Phi_{JK}^+\Rightarrow\Phi_{JL}^+$.
For $H\subsetneq L\subsetneq K\subsetneq J$, all comparisons are defined on
\[
 U_{JK}\cap U_{JL}\cap U_{JH}\cap
 \phi_{JK}^{-1}(U_{KL}\cap U_{KH})\cap\phi_{JL}^{-1}(U_{LH}).
\]
This contains $\psi_J^{-1}(S_J\cap S_K\cap S_L\cap S_H)$.
Both sides of
\begin{equation}\label{e:6.13}
 \Lambda_{JLH}\odot(\id_{\Phi_{LH}^+}*\Lambda_{JKL})
 =\Lambda_{JKH}\odot(\Lambda_{KLH}*\id_{\Phi_{JK}^+})
\end{equation}
have the zero bundle map as representative there, by the formulas
of Section~\ref{ss:3.2}. Thus \eqref{e:3.7} is satisfied.

The open subsets $S_J$ cover $X$, since they contain the $C_J$, and intersect only
for comparable indices, moreover, with order $J\prec K$ when $K\subsetneq J$,
these charts and coordinate changes form the desired Kuranishi atlas.

The relative assertion follows from the fact that restrictions $\pi_J:U_J\to B_{\RR}$ are submersions, and
$\pi_K\phi_{JK}=\pi_J$ holds on every coordinate domain.

\end{proof}

 We record the following technical statement needed in
Section~\ref{sec:7}.

\begin{lem}\label{r:6.6}\label{lem:prescribed-reductions}
Let $P_A$ be presentations as in Proposition~\ref{r:5.2}, indexed by
the nonempty finite subsets $A$ of a set $I$.
Write $V_A$ for their domains, $\mathcal G_A$ for
their sets of generators, and
$\psi_A:s_A^{-1}(0)\to X_A\subseteq X$ for their zero locus
identifications. For $\varnothing\ne A'\subsetneq A$, let
$\phi_{AA'}:V_A\to V_{A'}$ be a holomorphic submersion over $B$
and let
\[
 \Phi_{AA'}:\phi_{AA'}^{-1}P_{A'}\longrightarrow P_A
\]
be a morphism of sheaves of cdgas on $V_A$. Suppose that
$\mathcal G_{A'}\subseteq\mathcal G_A$ and
\[
 \Phi_{AA'}(h)=h\circ\phi_{AA'},\qquad
 \Phi_{AA'}(e)=e
\]
for local holomorphic functions $h$ on $V_{A'}$ and generators
$e\in\mathcal G_{A'}$. Require
$\psi_{A'}\circ\phi_{AA'}=\psi_A$ on $s_A^{-1}(0)$.
Let $\chi_{AA'}:E_A\to\phi_{AA'}^*E_{A'}$ and
$\xi_{AA'}:F_A\to\phi_{AA'}^*F_{A'}$ be the projections
onto the coordinates corresponding to the degree $-1$ and
degree $-2$ generators of $P_{A'}$, respectively.
Require $\phi_{AA'},\chi_{AA'},\xi_{AA'}$ to satisfy
\eqref{e:6.2}. At corresponding zeros, the morphism of
relative tangent complexes
$(\mathrm d\phi_{AA'},\chi_{AA'},\xi_{AA'})$
must also induce the identifications of tangent and obstruction
spaces in Lemma~\ref{r:5.4}, preserving the quadratic form
$Q_x$ on the obstruction space.

Let $\{C_A\}$ be a locally finite family of compact subsets of $X$,
with $C_A\subseteq X_A$ and
\[
 C_A\cap C_{A'}\ne\varnothing
 \quad\Longrightarrow\quad A\subseteq A'\text{ or }A'\subseteq A.
\]
Let $\mathcal F$ be a collection of indices closed under taking
nonempty subsets. For each $A\in\mathcal F$, suppose an admissible
subbundle $E_A^-$ is given near $\psi_A^{-1}(C_A)$, and that
\[
 \chi_{AA'}(E_A^-)\subseteq\phi_{AA'}^*E_{A'}^-
 \quad\text{near }\psi_A^{-1}(C_A\cap C_{A'})
 \quad(A'\subsetneq A,\ A,A'\in\mathcal F).
\]
Then there are admissible pairs $(U_A,E_A^-)$ for all $A$, with
$\psi_A^{-1}(C_A)\subseteq U_A$. For $A\in\mathcal F$, $E_A^-$ is
the restriction of the given subbundle to $U_A$.
The Kuranishi neighbourhoods
$(U_A,E_A^+,s_A^+,\psi_A|_{s_A^{-1}(0)\cap U_A})$,
with the induced coordinate changes, form a relative Kuranishi atlas
on an open neighbourhood of $\bigcup_A C_A$ in $X$.
\end{lem}

\begin{proof}
Repeat the construction in the proof of Proposition~\ref{r:6.5}.
\end{proof}

\par
\endgroup

\subsection{The derived manifold and its orientation}\label{ss:6.4}

\begin{prop}\label{r:6.7}\label{prop:manifold}
The atlas of Proposition~\ref{r:6.5} determines a derived smooth
manifold $M$ and a morphism $f^{\rm sm}:M\to B_{\RR}$ with underlying
map $f:X\to B$ and virtual dimension $d+2b$. The relative orientation
$o$ induces an orientation on $M$.
\end{prop}

\begin{proof}
The existence of $M$ and $f^{\rm sm}$ follows from Theorem~\ref{r:3.9}, and the dimension formula is given by
\eqref{e:6.3}.

We describe now the orientation, following \cite[\S6.3]{BJ}. At $x\in S_J$,
the relative tangent and obstruction spaces of the Kuranishi neighbourhood are
\[
 T_{x,\RR},\qquad O_{x,\RR}/N_{J,x}.
\]
Orient $T_{x,\RR}$ using its complex structure. Under the determinant convention
of Section~\ref{ss:2.5}, $o_x$ determines an orientation of the complex
quadratic space $(O_x,Q_x)$. Choose a $Q_x$ orthonormal complex basis
$e_1,\ldots,e_m$. If its determinant represents the orientation specified
by $o_x$, declare
\begin{equation}\label{e:6.14}
 e_1+N_{J,x},\ldots,e_m+N_{J,x}
\end{equation}
to be a positive real basis of $O_{x,\RR}/N_{J,x}$, and otherwise take the
opposite orientation.

The orientation is independent of the choice of basis.

We now check continuity of the pointwise orientation defined
above. Suppose first that the presentation is as in
Theorem~\ref{r:4.2} and minimal at a zero $v_0$.
On a sufficiently small neighbourhood $V_0$ of $v_0$,
$\operatorname{Re}q$ is negative definite on $E^-$ and
$(E,q)$ has a holomorphic orthonormal frame.
Its images in $E^+$, together with the complex orientation
of $T(V_0/B)$, define a continuous orientation of
$\det T^*(V_0/B)\otimes\det E^+$.
Under $K_{\bX/B}\cong\det E$, the determinant of the frame
represents an orientation which differs from $o$ by a
locally constant sign on $s^{-1}(0)\cap V_0$.
It suffices to show that the pointwise orientation defined
above differs from this continuous orientation by the same sign. This is a straightforward explicit check.
\end{proof}

\section{Independence of choices and virtual classes}\label{sec:7}\label{sec:classes}

We now prove independence of the choices in Proposition~\ref{r:6.7},
and deduce the assertions about virtual classes in Theorem~\ref{r:1.1}.
Throughout this section we fix $f:\bX\to B$, the relative form $\omega$,
and its orientation $o$, satisfying the hypotheses of that theorem.
We put $b=\dim_{\C}B$ and $d=\vdim_{\C}(\bX/B)$.

\begingroup
\color{black}
\subsection{Independence of choices}\label{ss:7.1}

\begingroup
\color{black}
The argument follows \cite[\S6.4]{BJ}.
\par
\endgroup

\begin{prop}\label{r:7.1}\label{prop:independence}
Let $M_0,M_1\to B_{\RR}$ be obtained by two choices in
Proposition~\ref{r:6.7}, with the orientations induced by $o$.
There is an oriented derived smooth manifold $W$ with boundary,
of virtual dimension $d+2b+1$, with the following properties.
\begin{enumerate}
\renewcommand{\labelenumi}{(\roman{enumi})}
\item Its underlying space is $X\times[0,1]$, and it has a morphism
$W\to B_{\RR}\times[0,1]$ inducing $(x,t)\mapsto(f(x),t)$.
\item There is an oriented equivalence over $B_{\RR}$
\begin{equation}\label{e:7.1}
 \partial W\simeq -M_0\amalg M_1,
\end{equation}
identifying the two boundary components with $X\times\{0\}$ and
$X\times\{1\}$.
\item Near each boundary component the construction has product
collars. The projection $|W|\to X$ is proper.
\end{enumerate}
If $X$ is compact, $W$ is an oriented bordism over $B_{\RR}$.
For every $u\in B$, the same construction produces a derived manifold
$W_u$ with underlying space $X_u\times[0,1]$ and boundary
$-(M_0)_u\amalg(M_1)_u$. If $f$ is proper, $W_u$ is a compact
bordism of virtual dimension $d+1$.
\end{prop}

\begin{proof}
\textit{Step 1.}
For $a=0,1$, let $I_a$ be the indexing set of the initial 
charts. Write $P_J^a$, $C_J^a$, $(U_J^a,E_J^{a,-})$, and $S_J^a$
for the presentations, compact closed sets, admissible pairs, and
open subsets used to construct $M_a$. We regard $I_0,I_1$ as disjoint.
The open subsets $S_J^a$ form a locally finite cover of $X$, have compact
closures, and intersect only for comparable indices $J$.

Apply Lemma~\ref{r:5.3}. It results in a system $P_L$ indexed by finite
nonempty $L\subset I_0\amalg I_1$, whose restriction to each $I_a$
is the old system on smaller domains containing the respective
$C_J^a$. Restrict $U_J^a$ and its coordinate domains accordingly,
and shrink $S_J^a$ further, still containing $C_J^a$, so that
\begin{equation}\label{e:7.2}
 S_J^a\subseteq\bigcap_{i\in J}R_i'',
\end{equation}
where $R_i''$ are the smaller open sets in that lemma. These
restrictions still cover $X$, since the $C_J^a$ do. They give the
same oriented derived manifold by Theorem~\ref{r:3.9}, and we use the
same notation for the restricted atlases.

Put
\[
 \widetilde{\bX}=\bX\times\C,\qquad
 \widetilde f=f\times\id_{\C}:\widetilde{\bX}\longrightarrow B\times\C,
 \qquad \widetilde\omega=\operatorname{pr}_{\bX}^*\omega.
\]
The relative cotangent complex is the pullback of
$\bL^{\an}_{\bX/B}$. Thus $\widetilde\omega$ is a relative
$(-2)$ shifted symplectic form, locally exact by
Definition~\ref{r:2.7}, with relative virtual dimension $d$.
Its orientation is $\widetilde o=\operatorname{pr}_{\bX}^*o$.
The presentations $P_L$, their comparisons, and their quadratic
forms extend to this family by pullback along the projection.

\smallskip
\noindent\textit{Step 2.}
Choose open subsets $D_0,D_1\subset\C$ which cover $\C$, such that
$D_0$ is disjoint from a neighbourhood of $1$ and $D_1$ is disjoint
from a neighbourhood of $0$. For example, take
\[
 D_0=\{z:\operatorname{Re}z<2/3\},\qquad
 D_1=\{z:\operatorname{Re}z>1/3\}.
\]
Put
\[
 I=\{(a,J):a\in\{0,1\},\ \varnothing\ne J\subseteq I_a,
                    \ J\text{ finite},\ S_J^a\ne\varnothing\}.
\]
For the remainder of the proof, $A$ denotes a nonempty finite
subset of $I$.

Apply Proposition~\ref{r:3.8} to the cover
\[
 \bigl\{(S_J^a\times D_a)\cap(X\times[0,1]):(a,J)\in I\bigr\}
\]
of $X\times[0,1]$. We obtain compact sets $\widetilde C_A$ with
\begin{equation}\label{e:7.3}
 \widetilde C_A\subseteq\bigcap_{(a,J)\in A}(S_J^a\times D_a),
 \qquad \bigcup_A\widetilde C_A=X\times[0,1].
\end{equation}
Indeed, each set is closed in $X\times[0,1]$ and contained in
$S_J^a\times[0,1]$ for any $(a,J)\in A$, whose closure is compact.
They are locally finite, and two of them intersect only for comparable
$A$. Since $X\times[0,1]$ is closed in $X\times\C$, the family
is locally finite in the latter space as well.

\begingroup
\color{black}

For each $A$, put
\[
 L(A)=\bigcup_{(a,J)\in A}J,\qquad
 D_A=\bigcap_{(a,J)\in A}D_a,\qquad
 Q_A=\bigcap_{(a,J)\in A}(S_J^a\times D_a).
\]
We consider those $A$ for which $Q_A\ne\varnothing$.
Let $\widetilde P_A$ be the pullback of $P_{L(A)}$ to
$V_{L(A)}\times D_A$ along the first projection. By
\eqref{e:7.2}, the image of its zero locus contains $\widetilde C_A$.
For $A'\subset A$, define the comparison by pulling back
$\Phi_{L(A)L(A')}$, with underlying map
\[
 (v,z)\longmapsto\bigl(\phi_{L(A)L(A')}(v),z\bigr).
\]
Here $D_A\subseteq D_{A'}$, and when $L(A')=L(A)$ we set
$\Phi_{L(A)L(A')}=\id$ and $\phi_{L(A)L(A')}=\id$.
The resulting comparisons satisfy \eqref{e:6.2}.
\par
\endgroup

\smallskip
\noindent\textit{Step 3.}
\begingroup
\color{black}

Fix $a\in\{0,1\}$, and let
$A=\{(a,J_1),\ldots,(a,J_r)\}$ satisfy $Q_A\ne\varnothing$.
Since the open sets $S_{J_1}^a,\ldots,S_{J_r}^a$ have nonempty
intersection, the sets $J_1,\ldots,J_r$ form a chain under inclusion.
Write $J(A)$ for the largest member. Then $L(A)=J(A)$, and
$\widetilde P_A$ is the pullback of the restriction of
$P_{J(A)}^a$ to $V_{L(A)}$. On
$(U_{J(A)}^a\cap V_{L(A)})\times D_a$, take
\[
 \widetilde E_A^-=\operatorname{pr}_1^*E_{J(A)}^{a,-}.
\]
For $\varnothing\ne A'\subsetneq A$, we have
$J(A')\subseteq J(A)$. If this inclusion is strict, compatibility
follows by pulling back the compatible subbundles
$E_{J(A)}^{a,-}$ and $E_{J(A')}^{a,-}$ and their coordinate
change in the atlas of $M_a$. If $J(A')=J(A)$, the comparison
is the identity and both subbundles equal
$\operatorname{pr}_1^*E_{J(A)}^{a,-}$.

If $\widetilde E_A^-$ is defined above,
so is $\widetilde E_{A'}^-$ for every
$\varnothing\ne A'\subseteq A$. By Lemma~\ref{r:6.6}, they extend
to compatible admissible pairs for all $A$ near $\widetilde C_A$.
Let $\widetilde S_A^0$ be the image in $X\times\C$ of the zero
locus of the resulting chart, and define
\[
 \widetilde S_A=\widetilde S_A^0\cap Q_A,\qquad
 \mathcal N=\bigcup_A\widetilde S_A.
\]
Restrict the chart domains and coordinate changes to these open
subsets. 
\endgroup

On the slice at $a$, every chart is an open restriction of a chart
in the given atlas for $M_a$. These restrictions cover the slice by \eqref{e:7.3}, and they consequently reconstruct $M_a$, with its map to $B_{\RR}$
and its orientation.

\smallskip
\noindent\textit{Step 4.}
Write $\widetilde U_{AA'}$ for the coordinate domain when $A'\subset A$.
For $A''\subset A'\subset A$, the two composites are compared on
\[
 \widetilde U_{AA'}\cap\widetilde U_{AA''}\cap
 \widetilde\phi_{AA'}^{-1}(\widetilde U_{A'A''}).
\]
This contains all zeros over
$\widetilde S_A\cap\widetilde S_{A'}\cap\widetilde S_{A''}$.
On this domain we have
\[
 \widetilde\phi_{AA''}
 =\widetilde\phi_{A'A''}\circ\widetilde\phi_{AA'}.
\]
\begingroup\color{black}
The induced bundle maps satisfy
\[
 \widetilde\chi_{AA''}^+
 =\widetilde\phi_{AA'}^*(\widetilde\chi_{A'A''}^+)
   \circ\widetilde\chi_{AA'}^+
\]
by the calculation in \eqref{e:6.12}. We take the $2$-isomorphism
comparing these coordinate changes to be represented by the zero
bundle map. For $A'''\subsetneq A''\subsetneq A'\subsetneq A$,
the cocycle identity \eqref{e:6.13} follows as in Step~3 of
Proposition~\ref{r:6.5}. The same argument applies when some sets
$L(A)$ coincide.
\par\endgroup

The maps $(V_{L(A)}\times D_A)_{\RR}\to(B\times\C)_{\RR}$ are
submersions. Using Proposition~\ref{r:6.7}, we thus get a derived smooth manifold
$\widetilde M$ with underlying space $\mathcal N$, relative virtual dimension
$d$, and a relative orientation induced by $\widetilde o$.
Its fibres at $z=0,1$, as derived manifolds over $B_{\RR}$,
are $M_0,M_1$.

\smallskip
\noindent\textit{Step 5.}
Choose a smooth map $\gamma:[0,1]\to[0,1]\subset\C$ equal to $0$
near $0$ and to $1$ near $1$. Apply Corollary~\ref{r:3.12} to the map
\[
 h:B_{\RR}\times[0,1]\longrightarrow(B\times\C)_{\RR},
 \qquad h(u,t)=(u,\gamma(t)),
\]
and put
\begin{equation}\label{e:7.4}
 W=\widetilde M\times_{(B\times\C)_{\RR}}(B_{\RR}\times[0,1]).
\end{equation}
Since $\mathcal N$ contains $X\times[0,1]$, its underlying space identifies
with $X\times[0,1]$ by
$(x,t)\mapsto((x,\gamma(t)),(f(x),t))$. The virtual dimension is
$d+2b+1$. Constancy of $\gamma$ near its endpoints makes the
pulled back data of the required product form there.

Orient $B_{\RR}$ by its complex structure and $[0,1]$ by
$\partial/\partial t$. Use the relative orientation of $\widetilde M$
and the convention of Section~\ref{ss:3.5}: we get therefore
\eqref{e:7.1}, with precisely the original orientations of $M_0,M_1$.

The projection $X\times[0,1]\to X$ is proper. Pulling the atlas of
$W$ back to $\{u\}\times[0,1]$, one gets $W_u$, with the asserted
underlying space, boundary and virtual dimension. If $f$ is proper,
$X_u$ is compact, and hence so is $W_u$.
\end{proof}
\par
\endgroup

\subsection{Compatibility with base change}\label{ss:7.2}

\begin{prop}\label{r:7.2}\label{prop:pullback}
Fix a relative Kuranishi atlas constructed in Proposition~\ref{r:6.7},
and write $M\to B_{\RR}$ for its associated derived smooth manifold.
For a smooth map $g:S\to B_{\RR}$, with $S$ a manifold possibly
with boundary, the pulled back atlas represents
\[
 M_S=M\times_{B_{\RR}}S,\qquad |M_S|=X\times_B S,
 \qquad \vdim M_S=d+\dim_{\RR}S.
\]
It has the pulled back relative orientation.

If $S$ is complex and $g$ is holomorphic,
this atlas is obtained by applying the construction of
Proposition~\ref{r:6.5} to the pulled back cell presentations,
with the pulled back admissible subbundles. Its relative form and orientation are
the pullbacks of $\omega,o$.
\end{prop}

The first assertion follows from Proposition~\ref{r:3.11} and
Corollary~\ref{r:3.12}(i). We prove the second one below.

\begingroup
\color{black}
\begin{proof}
For each $J$, let $p:V_J^S=V_J\times_B S\to V_J$ be the
projection, and put
\[
 P_J^S=\cO_{V_J^S}\otimes_{p^{-1}\cO_{V_J}}p^{-1}P_J,
\]
with the induced differential. Since $V_J\to B$ is a holomorphic
submersion, $V_J^S$ is a complex manifold smooth over $S$.
By the analytic cell construction in Section~\ref{sec:2},
$\Specan P_J^S\simeq\Specan P_J\times_B^{\R}S$.
Base change therefore identifies its relative tangent
and obstruction spaces with those of $\bX_J\times_B^{\R}S$,
as in Lemma~\ref{r:5.4}, with the pulled back quadratic forms.

In the notation of \eqref{e:6.1}, we have
\[
 \begin{gathered}
 T(V_J^S/S)=p^*T(V_J/B),\qquad E_J^S=p^*E_J,
 \qquad F_J^S=p^*F_J,\\
 s_J^S=p^*s_J,\qquad t_J^S=p^*t_J.
 \end{gathered}
\]
At each zero of $s_J^S$, its derivative along the fibres over
$S$ is $p^*D_J$. Hence $p^*E_J^-$ is admissible on
$U_J\times_B S$ by Definition~\ref{r:6.1}, and
\begin{equation}\label{e:7.5}
 E_J^S/p^*E_J^-=p^*(E_J/E_J^-),\qquad
 (s_J^S)^+=p^*s_J^+.
\end{equation}
Under these identifications, the coordinate changes on
$U_{JK}\times_B S$ are
$(v,s)\mapsto(\phi_{JK}(v),s)$, with bundle maps
$p^*\chi_{JK}^+$. Thus the Kuranishi neighbourhoods and coordinate
changes constructed from $P_J^S$ and $p^*E_J^-$ are precisely
those of the pulled back atlas.
The relative form is the pullback of $\omega$, and
Proposition~\ref{r:6.7} identifies the relative orientation
induced by the pullback of $o$ with that on $M_S$.
\end{proof}
\par
\endgroup

\subsection{Virtual fundamental classes}\label{ss:7.3}

Suppose now that $f$ is proper. For $u\in B$, write
\[
 \bX_u=\bX\times_B^{\R}\{u\},\qquad
 \omega_u=\omega|_{\bX_u},\qquad o_u=o|_{\bX_u}.
\]
The derived smooth fibre $M_u=M\times_{B_{\RR}}\{u\}$ has
virtual dimension $d$ and underlying space $X_u$. It is compact.
The space $X_u$ is an ENR, by the facts recalled in
Section~\ref{ss:3.6}. Theorem~\ref{r:3.15} therefore defines
\begin{equation}\label{e:7.6}
 [\bX_u,\omega_u,o_u]:=[M_u]^{\mathrm{vir}}\in H_d(X_u;\ZZ).
\end{equation}
This class is zero when $d<0$.

\begin{prop}\label{r:7.3}\label{prop:virtual-class-independence}
The class \eqref{e:7.6} is independent of all choices.
\end{prop}

\begingroup
\color{black}
\begin{proof}
Let $M_0,M_1\to B_{\RR}$ arise from two choices in
Proposition~\ref{r:6.7}. By Proposition~\ref{r:7.1},
there is a compact oriented bordism $W_u$ from
$(M_0)_u$ to $(M_1)_u$, with underlying space
$X_u\times[0,1]$. Applying Proposition~\ref{r:3.16}
to the projection $X_u\times[0,1]\to X_u$, we obtain
\[
 [(M_0)_u]^{\mathrm{vir}}
 =
 [(M_1)_u]^{\mathrm{vir}}
 \quad\text{in }H_d(X_u;\ZZ).
\]
Thus the class \eqref{e:7.6} is independent of all choices.
\end{proof}
\par
\endgroup

\begin{prop}\label{r:7.4}\label{prop:deformation-invariance}
Suppose $f$ is proper and $u_0,u_1$ belong to the same connected
component of $B$. Then their derived smooth fibres are oriented
bordant, and
\[
 (i_{u_0})_*[\bX_{u_0},\omega_{u_0},o_{u_0}]
 =(i_{u_1})_*[\bX_{u_1},\omega_{u_1},o_{u_1}]
 \quad\text{in }H_d(X;\ZZ).
\]
\end{prop}

\begin{proof}
Choose a smooth path $\gamma:[0,1]\to B$ from $u_0$ to $u_1$,
constant near its endpoints. The smooth pullback
\[
 W_{\gamma}=M\times_{B_{\RR}}[0,1]
\]
is a derived manifold with boundary, of virtual dimension $d+1$.
By the boundary and orientation conventions in Section~\ref{ss:3.5},
\[
 |W_{\gamma}|=X\times_B[0,1],\qquad
 \partial W_{\gamma}\simeq -M_{u_0}\amalg M_{u_1}.
\]
It is compact: $f^{-1}(\gamma([0,1]))$ is compact by properness,
and $X\times_B[0,1]$ is a closed subset of its product with
$[0,1]$. This proves the bordism assertion.

Let $q:|W_{\gamma}|\to X$ be the projection. Its compact image
has an open ENR neighbourhood $V\subseteq X$, as recalled in
Section~\ref{ss:3.6}. By Proposition~\ref{r:3.16}, applied to
$q:|W_{\gamma}|\to V$, the two pushed forward
classes in $H_d(V;\ZZ)$ are equal. Pushing forward by $V\hookrightarrow X$
then proves the displayed equality. Finally, Proposition~\ref{r:7.3}
identifies these classes with \eqref{e:7.6}.
\end{proof}

\begin{proof}[Proof of Theorem~\ref{r:1.1}]
Part (i) is Proposition~\ref{r:6.7}, and part (ii) is
Proposition~\ref{r:7.1}. The definition \eqref{e:7.6} and
Proposition~\ref{r:7.3} prove part (iii).
Proposition~\ref{r:7.4} proves part (iv).
\end{proof}

\section{Stable pairs with proper support}\label{sec:pt-support}

\subsection{The moduli problem}\label{ss:pt-moduli}

Let $\pi:Y\to B$ be a holomorphic submersion of relative complex
dimension four between Hausdorff second countable complex manifolds. The fibres of $\pi$ may be noncompact.
Fix a nowhere vanishing relative holomorphic volume form
$\Omega\in H^0(Y,K_{Y/B})$ (and, in particular, suppose that such form exists).

A \emph{stable pair} on $Y_b$ is a morphism
$s:\cO_{Y_b}\to F$, where $F$ is a nonzero coherent pure
one-dimensional sheaf with compact support and
$\operatorname{coker}s$ has support of dimension at most zero.
Its \textit{associated complex} is $I=[\cO_{Y_b}\xrightarrow{s}F]$,
with terms in degrees zero and one.

See also \cite{PT}.

By Theorem~\ref{app:pt-representability}, stable pairs
represented by a separated Hausdorff second countable complex
analytic space $P(Y/B)$. 

Let $P\subseteq P(Y/B)$ be an open
subspace (we will be interested in $\chi=n$-components later for various $n$). For a derived analytic space $T\to B$, let $\mathbf P(T)$ be
the $\infty$-groupoid of perfect complexes on $Y_T$ with trivialized
determinant whose restriction to $Y_{t_0(T)}$ corresponds to
some family classified by a map $t_0(T)\to P$.

\subsection{The deformation complex}\label{ss:pt-deformations}

On $Y\times_B\mathbf P$, write $\mathcal I$ for the universal
perfect complex and $p:Y\times_B\mathbf P\to\mathbf P$ for the
projection. Put
\begin{equation}\label{e:pt-deformation}
 \mathcal E=\operatorname{fib}\bigl(
 \R\mathcal Hom(\mathcal I,\mathcal I)
 \xrightarrow{\operatorname{tr}}\cO\bigr),
 \qquad \mathcal A=\R p_*\mathcal E.
\end{equation}

 Let $Z$ denote the universal classical support. Off $Z$,
the complex $\mathcal I$ is a line bundle, and its determinant
trivialization identifies it with $\cO$. 

\begin{prop}\label{prop:pt-derived-moduli}
The functor $\mathbf P$ defined in Section~\ref{ss:pt-moduli}
is represented by a derived complex analytic space. Its relative tangent complex is perfect of amplitude
$[0,2]$, with
\[
 \bT_{\mathbf P/B}\simeq\mathcal A[1],
 \qquad \mathcal A^\vee\simeq\mathcal A[4].
\]
\end{prop}

\begin{proof}
We apply \cite[Theorem~7.1]{PortaYu}. The classical truncation of
$\mathbf P$ is $P$ by Theorem~\ref{app:pt-representability} and
Proposition~\ref{app:pt-comparison}.
As $\mathcal E$ is perfect with
proper support, $\mathcal A$ is perfect and commutes with derived
base change, so the deformation theory of perfect complexes with
trivialized determinant identifies the relative analytic cotangent
complex with $(\mathcal A[1])^\vee$.

See also \cite[\S2.3]{PTVV} and
\cite[Proposition~3.2]{STV}.

Thus one may apply Porta--Yu's
representability theorem \cite[Theorem~7.1]{PortaYu}. 

Finally, by relative Serre duality
\cite{RR,RRV}, we have
$\mathcal A^\vee\simeq\mathcal A[4]$. On each fibre,
$\operatorname{Ext}^j(I,I)_0=0$ for $j\leq0$ by
Proposition~\ref{app:pt-comparison} and the absence of
automorphisms of stable pairs, and hence also for $j\geq4$ by duality, so
$\mathcal A[1]$ has amplitude $[0,2]$.
\end{proof}

\subsection{The relative form}\label{ss:pt-trace}

For a closed subset $Z\subseteq t_0(V)$, write
\[
 \R\Gamma_Z(V,K)=\operatorname{fib}\bigl(
 \R\Gamma(V,K)\longrightarrow\R\Gamma(V\setminus Z,K)\bigr).
\]
The canonical identification $\mathcal I|_{Y_{\mathbf P}\setminus Z}
\simeq\cO$ defines a (compactly supported) second Chern character
\begin{equation}\label{e:pt-supported-character}
 c_Z(\mathcal I)\in
 H^4\R\Gamma_Z(Y_{\mathbf P},F^2\Dhat_{Y_{\mathbf P}/B}).
\end{equation}

(The analytic stack of perfect complexes
is the analytification of the algebraic one by
\cite[Proposition~4.9]{HP}, and thus the construction in
\cite[\S2.3]{PTVV} proves \eqref{e:pt-supported-character} in the
analytic de Rham complex.)

\begin{prop}\label{prop:pt-trace}
There is a natural relative $(-2)$-shifted symplectic form $\omega$
on $\mathbf P/B$, whose underlying pairing is the Serre pairing of
Proposition~\ref{prop:pt-derived-moduli}. This form is compatible
with derived analytic base change.
\end{prop}

\begin{proof}
Let $T\to B$ be a derived analytic space, write $p_T:Y_T\to T$
and $q_T:Y_T\to Y$ for the projections, and let
$Z_T\subseteq t_0(Y_T)$ be closed and proper over $t_0(T)$.
Using the canonical splitting
\[
 \bL^{\an}_{Y_T/B}\simeq
 p_T^*\bL^{\an}_{T/B}\oplus q_T^*\Omega^1_{Y/B},
\]
we follow the construction in \cite[\S6.2]{Park}: project to
terms of degree zero in the $Y/B$ differentials, multiply by
$q_T^*\Omega$, and apply the relative analytic trace with supports
proper over $T$ \cite{RR,RRV}. We thus obtain a morphism of
filtered complexes
\begin{equation}\label{e:pt-integration}
 \mathcal T_\Omega:
 \R\Gamma_{Z_T}(Y_T,\Dhat_{Y_T/B})
 \longrightarrow \R\Gamma(T,\Dhat_{T/B})[-4].
\end{equation}

Taking $T=\mathbf P$ and $Z_T=Z$, set
$\omega=\mathcal T_\Omega(c_Z(\mathcal I))$.
As \eqref{e:pt-integration} preserves the Hodge filtration,
$\omega$ is a closed relative $2$-form of degree $-2$.
To check nondegeneracy, recall that the underlying $2$-form of the
Chern character is
$\tfrac12\operatorname{tr}(\operatorname{At}(\mathcal I)^2)$.
On $\bT_{\mathbf P/B}\simeq\mathcal A[1]$, its image under
$\mathcal T_\Omega$ is therefore the Serre pairing of
Proposition~\ref{prop:pt-derived-moduli}, which is nondegenerate.

\end{proof}

\section{Local exactness}\label{sec:pt-exactness}

We prove that the image of the relative trace form $\omega$ in the
full completed de Rham complex is locally zero. (We shall test this
image after Artinian base changes.)

\begin{lem}\label{pt:de-rham-detection}
Let $M\to B$ be a derived complex analytic space which locally has a
cell presentation
\[
 A=\cO_U[\xi_1,\ldots,\xi_r,\eta_1,\ldots,\eta_q],
 \qquad |\xi_a|=-1,\quad |\eta_i|=-2,
\]
where $U\to B$ is smooth, $B$ is classical, and
\[
 \delta\xi_a=s_a,\qquad
 \delta\eta_i=\sum_a g_{ia}\xi_a.
\]
Put $J=(s_1,\ldots,s_r)$, so that $M_0=t_0M$ is the closed analytic
subspace of $U$ defined by $J$.
\begin{enumerate}
\item There is a natural equivalence of complexes
\begin{equation}\label{pt:formal-dr}
 \Dhat_{M/B}\simeq
 \left[
 \widehat{\cO}_{U,M_0}\longrightarrow
 \widehat\Omega^1_{U/B,M_0}\longrightarrow\cdots
 \longrightarrow\widehat\Omega^n_{U/B,M_0}
 \right],
\end{equation}
where $n=\dim(U/B)$ and
$\widehat\Omega^p_{U/B,M_0}=\varprojlim_N
\Omega^p_{U/B}/J^N\Omega^p_{U/B}$, restricted to $M_0$.
In particular, $\Dhat_{M/B}$ has no negative cohomology sheaves.

\item The restriction map
\begin{equation}\label{pt:artin-tests}
 H^0\R\Gamma(M,\Dhat_{M/B})\longrightarrow
 \prod_{T\to B}H^0\R\Gamma(M\times_B T,\Dhat_{(M\times_B T)/T})
\end{equation}
is injective, where $T$ ranges over local Artinian complex analytic
spaces.

\item If $B=\Specan R$ for a local Artinian
$\C$-algebra $R$, then
\begin{equation}\label{pt:artin-constants}
 H^0\R\Gamma(M,\Dhat_{M/R})
 \cong H^0(M_0^{\mathrm{top}},\underline R).
\end{equation}
Restriction to a complex point $m$ is evaluation of this locally
constant function, with values in
$H^0\Dhat_{\Specan\C/R}\cong R$.
\end{enumerate}
All completions in \eqref{pt:formal-dr} are completions of sheaves.
\end{lem}

\begin{proof}
We give the local calculation for (i), which is the (analytic) relative of  \cite[Corollary~4.27]{Bhatt}.

Choose one Stein
coordinate neighbourhood and put $u_a=\mathrm d\xi_a$ and
$v_i=\mathrm d\eta_i$. Their total degrees are $0$ and $-1$.
With $D=\delta+\mathrm d$, we have
\[
 D\xi_a=s_a+u_a,\qquad Du_a=-\mathrm d s_a,
 \qquad D\eta_i=v_i+\sum_a g_{ia}\xi_a.
\]
Replace $v_i$ by $v'_i=D\eta_i$. Since $D\eta_i=v'_i$ and
$Dv'_i=0$, the quotient obtained by setting all $\eta_i=v'_i=0$
is quasi-isomorphic to the original complex. This remains true
after completion: apart from the $u_a$, the Hodge weight is bounded
by $n+q$, and in each total degree only finitely many monomials in
the $\eta_i$ occur. The quotient complex is
\[
 \bigl(\Omega^\bullet_{U/B}[[u_1,\ldots,u_r]]
       \otimes\Lambda(\xi_1,\ldots,\xi_r),D\bigr).
\]
Evaluation at $u_a=-s_a$, $\xi_a=0$ maps this complex to the
RHS of \eqref{pt:formal-dr}. To check that it is an
equivalence, filter by ordinary form degree on $U/B$. On each
associated graded term, it is the Koszul resolution for $u_a+s_a$,
completed along $u=0$. Before completion this is the resolution of
the graph $u=-s$ in a product with a polydisc. By Artin--Rees, completion along $(u_1,\ldots,u_r)$
preserves this Koszul resolution. On the graph $u=-s$,
this ideal is $J$, so in ordinary degree $p$ the
completed complex resolves $\widehat\Omega^p_{U/B,M_0}$.

Notice that evaluation commutes with the differentials: it sends
$D\xi_a=s_a+u_a$ to zero and $Du_a=-\mathrm d s_a$
to $\mathrm d(-s_a)$.
It is therefore a quasi-isomorphism on each associated
graded complex, and hence on the whole complex, since
ordinary form degree ranges only from $0$ to $n$.
This proves (i).

For a similar argument, see
\cite[Lemma~6.1.2(2)--(3)]{Park}.

For (ii), note that (i) identifies a degree $0$ class on a chart with a formal
function
\[
 h\in\Gamma(M_0,\widehat{\cO}_{U,M_0}),\qquad
 \mathrm d_{U/B}h=0.
\]
Suppose that all its restrictions in \eqref{pt:artin-tests} vanish.
For $x\in M_0$ and $k\geq1$, set
$T_{x,k}=\Specan(\cO_{U,x}/\mathfrak m_x^k)$.
Pull back further along
$M\times_U T_{x,k}\to M\times_B T_{x,k}$.
The pulled-back cell presentation has degree-zero part
$\cO_{T_{x,k}}$, so (i) applies with the identity
$T_{x,k}\to T_{x,k}$ in place of $U\to B$.

By (i), the full completed de Rham complex of
$M\times_U T_{x,k}$ over $T_{x,k}$ is quasi-isomorphic to
$\cO_{U,x}/\mathfrak m_x^k$, since the image of $J$
in this ring is nilpotent.
Thus the image of $h$ in
$\cO_{U,x}/\mathfrak m_x^k$ is zero for every $x,k$. For fixed $N$, the image
$h_N$ in $\cO_U/J^N$ has zero germ modulo $\mathfrak m_x^k$ for
every $k$. Krull intersection proves $(h_N)_x=0$ for every $x$,
and hence $h_N=0$. This holds for all $N$, proving $h=0$.
Since the complex in (i) is nonnegative, its degree $0$
hypercohomology is the space of sections of its degree $0$
cohomology sheaf, and the chartwise argument therefore proves (ii).

\begingroup\color{black}
For (iii), write $U=U'\times\Specan R$, where $U'$ is a polydisc,
and let $Z=(M_0)_{\mathrm{red}}$ be the reduced analytic space
underlying $M_0$, regarded as a closed subspace of $U'$.
Since $R$ is finite dimensional and completion depends only on
the radical of the ideal, the complex in \eqref{pt:formal-dr}
is $R\otimes_\C\widehat\Omega^\bullet_{U',Z}$.
For a formal function $h$ with $\mathrm dh=0$, its image in
$R\otimes_\C\cO_Z$ is locally constant: on each local irreducible component it is
constant on the smooth locus, and hence everywhere by continuity.

The Krull intersection argument in (ii) therefore shows that
$h$ equals this locally constant function.
Thus $\cH^0(\Dhat_{M/R})=\underline R$.
As the complex is nonnegative, taking global sections proves
\eqref{pt:artin-constants}.

For the closed immersion $\Specan\C\hookrightarrow\Specan R$,
the algebraic and analytic cotangent complexes agree
\cite[Corollary~5.33]{PortaYu}. The natural comparison of completed
de Rham complexes is therefore an equivalence: it is an equivalence
on each Hodge graded piece, and both complexes arecomplete.
Since the maximal ideal of $R$ is nilpotent,
\cite[Example~4.18]{Bhatt} now shows that the canonical unit map
\[
 R\xrightarrow{\sim}\Dhat_{\Specan\C/R}
\]
is an equivalence. Pullback along $m:\Specan\C\to M$ preserves
constants from $R$. Thus, identifying its source by
\eqref{pt:artin-constants} and its target by the displayed map,
it sends a locally constant function $h$ to $h(m)\in R$.
\par\endgroup
\end{proof}

\begin{prop}\label{pt:local-exactness}
The relative closed $(-2)$-shifted symplectic form $\omega$ on
$\mathbf P/B$ is locally exact in the full completed relative
analytic de Rham complex.
\end{prop}

\begin{proof}
The relative cotangent complex has amplitude $[-2,0]$ by the
deformation calculation above. Thus Proposition~\ref{prop:presentation}
provides the presentations in Lemma~\ref{pt:de-rham-detection},
before any use of the Darboux theorem. Let
\[
 a\in H^0\R\Gamma(\mathbf P,\Dhat_{\mathbf P/B})
\]
be the image of $\omega$. By Lemma~\ref{pt:de-rham-detection}(ii),
it is enough to prove that $a$ becomes zero after every local
Artinian base change $T=\Specan R\to B$.

Write $b$ for a closed point of $T$. 

Fix a complex
point $m$ corresponding to $I_m=[\cO_{Y_b}\to F_m]$, and put
$Z_m=\operatorname{Supp}F_m$. The  second Chern
character of $I_m$
is $-\operatorname{ch}_{2,Z_m}(F_m)$. Its de Rham class lies in
\[
 H^4\R\Gamma_{Z_m}(Y_b,\Omega^\bullet_{Y_b})
 \cong H^4(Y_b,Y_b\setminus Z_m;\C)
 \cong H^{\mathrm{BM}}_4(Z_m;\C)=0.
\]
Here the first isomorphism is the holomorphic Poincar\'e lemma,
the second is duality with supports on the oriented real
eight-manifold $Y_b$, and the last equality follows from
$\dim_\C Z_m\leq1$.

Change from absolute differentials to differentials relative to
$R$. The de Rham class remains zero in
$H^4\R\Gamma_{Z_m}(Y_b,\Dhat_{Y_b/R})$.
By naturality of the  character, this is precisely
the restriction at $m$ of the universal class defining $\omega$.
The integration map $\mathcal T_\Omega$ commutes with this
restriction, and therefore sends it to
\[
 m^*(a|_T)=0\quad\text{in}\quad
 H^0\Dhat_{\Specan\C/R}\cong R.
\]

Since $m$ was arbitrary, $a|_T=0$.
Lemma~\ref{pt:de-rham-detection}(ii) now shows that $a=0$ on each
analytic chart, and, thus, establishes  the required local exactness.

This argumeny reproduces, in fact, just the calculation of \cite[Appendix C]{CZZ},
combined with the  argument of
\cite[Proposition 6.2.1]{Park}.
\end{proof}

\section{Orientations and virtual classes for stable pairs}\label{sec:pt-orientations-classes}

\subsection{Orientations}\label{ss:PT-orientations}

We use the orientation convention in
\cite[Remark~3.18]{KarpovOrientations}.

\begin{prop}\label{prop:PT-family-orientation}
Let $\pi:Y\to B$ be a holomorphic submersion of relative complex
dimension four, with a nowhere vanishing relative holomorphic volume
form. The derived analytic space $\mathbf P/B$ of
Section~\ref{ss:pt-moduli} has a canonical relative orientation.
This orientation is compatible with
base change and with restriction to open subsets of $Y$ which
contain the support.
\end{prop}

\begin{proof}
Follows analogously to the proof of
\cite[Proposition~6.7]{KarpovOrientations}, applied in
neighbourhoods of the supports.
\end{proof}

\subsection{Virtual fundamental classes}\label{ss:pt-classes}

\begin{thm}\label{thm:pt-classes}
Let $\pi:Y\to B$ be a holomorphic submersion of relative
complex dimension four between Hausdorff second countable complex
manifolds, with a nowhere vanishing relative holomorphic volume form.
Let $P\subseteq P(Y/B)$ be open with $n=\chi(F)$ constant,
and let $\mathbf P$ be its derived analytic enhancement.
Then $\mathbf P/B$ has a canonical orientation and a locally exact relative
$(-2)$-shifted symplectic form. Theorem~\ref{thm:main} therefore
associates to it an oriented derived smooth manifold
\[
 P^{\mathrm{sm}}\longrightarrow B_{\RR},
 \qquad |P^{\mathrm{sm}}|=P,\qquad
 \vdim P^{\mathrm{sm}}=2n+2b.
\]
It is independent of choices up to the bordisms of
Theorem~\ref{thm:main}(ii).

For every compact fibre $P_u$, there is a canonical class
\begin{equation}\label{e:pt-vfc}
 [P_u]^{\mathrm{vir}}\in H_{2n}(P_u;\ZZ).
\end{equation}
These classes are compatible with holomorphic base change between
complex manifold bases. If $P\to B$ is proper and $u_0,u_1$ lie
in the same connected component of $B$, then
\[
 (i_{u_0})_*[P_{u_0}]^{\mathrm{vir}}
 =(i_{u_1})_*[P_{u_1}]^{\mathrm{vir}}
 \quad\text{in }H_{2n}(P;\ZZ).
\]
\end{thm}

\begin{proof}
 $P$ is Hausdorff and
second countable by Theorem~\ref{app:pt-representability}.

Propositions~\ref{prop:pt-derived-moduli}, \ref{prop:pt-trace},
\ref{pt:local-exactness} and \ref{prop:PT-family-orientation}
verify, respectively, the perfectness, symplectic, exactness and orientation
hypotheses of Theorem~\ref{thm:main}.

To compute the dimension, work on one fibre $Y_u$ and put
$I=[\cO_{Y_u}\to F]$. In $K$-theory with compact support,
\[
 [\R\mathcal Hom(I,I)_0]
 =[F^\vee\Lotimes F]-[F]-[F^\vee].
\]
The Chern character of $F$ starts in cohomological degree six,
so, by Riemann--Roch with compact supports, 
$\chi(F^\vee\Lotimes F)=0$. By Serre duality,
$\chi(F^\vee)=\chi(F)=n$. Consequently
\[
 \vdim_{\C}(\mathbf P/B)
 =-\chi\bigl(\R\mathcal Hom(I,I)_0\bigr)=2n.
\]

Independence of choices and deformation invariance follow from
Propositions~\ref{r:7.1}, \ref{r:7.3} and \ref{r:7.4}.
\end{proof}

\section{Non-compact situation}\label{sec:noncompact}

We now give a Borel-Moore variant of the constructions above.

\subsection{Borel--Moore virtual classes}\label{ss:proper-bordisms}

We use integral Borel--Moore homology, with the convention
\[
 H_k^{\mathrm{BM}}(T;\ZZ)=H^{-k}(T,\omega_T),
\]
where $\omega_T$ is the topological dualizing complex.

A metrizable space $T$ is an \emph{absolute neighbourhood retract}, or
\emph{ANR}, if every continuous map from a closed subset of a metrizable
space to $T$ extends to a neighbourhood of that subset.
Hausdorff second countable complex analytic spaces are locally compact
ANRs: they are locally ENRs by Section~\ref{ss:3.6}, and a metrizable
space covered by open ANRs is an ANR \cite[Theorem~3.3]{Hanner}.

Let $T$ be a locally compact second countable ANR. Define
$B_k^{\mathrm{prop}}(T)$ as the group of bordism classes of pairs
$(N,q)$, where $N$ is an oriented $k$-manifold without boundary and
$q:N\to T$ is continuous and proper. 

Define $\mathrm dB_k^{\mathrm{prop}}(T)$ in the same way, using
oriented derived smooth manifolds of virtual dimension $k$.
For each bordism $W$, we require a neighbourhood of its boundary
to be identified with $\partial W\times[0,\epsilon)$, so that
its map $q:|W|\to T$ satisfies $q(x,t)=q(x,0)$ in these coordinates.
Addition of such classes is disjoint union, and reversing orientation provides
the inverse.
\begin{lem}\label{lem:proper-bordism}
Inclusion of ordinary manifolds induces an isomorphism
\[
 B_k^{\mathrm{prop}}(T)\xrightarrow{\sim}
 \mathrm dB_k^{\mathrm{prop}}(T),\qquad k\geqslant0.
\]
For $k<0$, the latter group is zero.
\end{lem}

\begin{proof}
We adapt the proof of \cite[Theorem~13.11]{JoyceBook}, see also
\cite[Remark~13.6(c)]{JoyceBook}.
\smallskip\noindent\emph{Step 1.}
Let $q:|M|\to T$ be proper, and let $C\subseteq|M|$ be compact.
By \cite[Corollary~4.36(ii)]{JoyceBook}, there is an open
neighbourhood $M_0$ of $C$ with $M_0\simeq S_{V,E,s}$.
Identify $|M_0|$ with $s^{-1}(0)$. Since $T$ is an ANR,
$q|_{|M_0|}$ extends to a continuous map
$\widetilde q:V'\to T$ for some open neighbourhood $V'$ of
$s^{-1}(0)$ in $V$. Choose $V'$ small enough that the orientation
of $\det T^*V\otimes\det E$ also extends over $V'$, and replace
$V$ by $V'$.

Choose a compactly supported smooth section $\tau$ of $E$ such
that $s+\tau$ is transverse to zero near $C$. Let
$\rho:[0,1]\to[0,1]$ be smooth, equal to $0$ near $0$ and to $1$
near $1$. Writing $p:V\times[0,1]\to V$ for the projection,
define a section of $p^*E$ by
\[
 \widetilde s(v,t)=s(v)+\rho(t)\tau(v).
\]
Its standard model
$W_0=S_{V\times[0,1],p^*E,\widetilde s}$
is a derived manifold with boundary. Near $t=0$ it is the
product of $S_{V,E,s}$ with an interval, and near $t=1$ it is
the product of $S_{V,E,s+\tau}$ with an interval. Its boundary is
\[
 \partial W_0\simeq -S_{V,E,s}\amalg S_{V,E,s+\tau}.
\]

Choose a compact neighbourhood $K$ of $\operatorname{supp}\tau$
in $V$. Outside $K$, the section $\widetilde s$ is $p^*s$.
We may therefore glue $W_0$ to
$(M\setminus(K\cap s^{-1}(0)))\times[0,1]$ along their common
open subset. This produces a derived manifold $W$ with
$\partial W\simeq-M\amalg M'$, where $M'$ is obtained by
replacing $s$ with $s+\tau$ in $M_0$.
The maps $(v,t)\mapsto\widetilde q(v)$ and
$(x,t)\mapsto q(x)$ agree on the overlap and define a continuous
map $Q:|W|\to T$. It is proper because $q$ is proper and all
changes are contained in the compact set $K\times[0,1]$.

Where $s+\tau$ is transverse to zero, its standard model is
equivalent to the smooth manifold $(s+\tau)^{-1}(0)$.

Going further, we call a derived manifold \emph{ordinary} when it is equivalent
to a smooth manifold regarded as a derived manifold.

Suppose that $M$ is ordinary on an open subset $O$, and let
$D\subseteq O$ be closed in $|M|$. Then $s$ is transverse
at the zeros corresponding to $O$. Choose open neighbourhoods
$V_1,V_2$ of $D\cap|M_0|$ in $V$, with
$\overline V_1\subseteq V_2$, such that $s$ is transverse
near $\overline V_2$. Choose finitely many compactly supported
sections of $E$ which vanish on $V_1$ and whose values span
$E$ near $C\setminus V_2$. By Sard's theorem, an arbitrarily
small linear combination $\tau$ makes $s+\tau$ transverse
there. For sufficiently small coefficients, transversality
is also preserved near $C\cap\overline V_2$. Hence the
construction leaves $M$ unchanged near $D$.

\smallskip\noindent\emph{Step 2.}
Choose compact subsets $K_j\subseteq\operatorname{int}K_{j+1}$
whose interiors cover $T$, and put $K_{-1}=K_0=\varnothing$.
Starting with $(M_0,q_0)=(M,q)$, we construct $(M_j,q_j)$
so that $M_j$ is ordinary near $q_j^{-1}(K_j)$.
Apply Step~1 to the compact set $q_{j-1}^{-1}(K_j)$.
By induction, $M_{j-1}$ is already ordinary near
$q_{j-1}^{-1}(K_{j-1})$. The last paragraph of Step~1
therefore allows us to leave a neighbourhood of
$q_{j-1}^{-1}(K_{j-2})$ unchanged.

Using the extension $\widetilde q$ from Step~1, choose
the perturbing sections with compact support in
\[
 \widetilde q^{-1}
 \bigl(\operatorname{int}K_{j+1}\setminus K_{j-2}\bigr).
\]
As in Step~1, we can make $s+\tau$ transverse to zero on
an open neighbourhood $U$ of
$s^{-1}(0)\cap\widetilde q^{-1}(K_j)$ using arbitrarily
small coefficients. Let $L$ be the union of the supports
of the perturbing sections. Since $s$ is nowhere zero on
the compact set
\[
 (L\cap\widetilde q^{-1}(K_j))\setminus U,
\]
we may also choose the coefficients so that $|\tau|<|s|$
there. Thus $s+\tau$ has no zeros on this set. Outside $L$
it agrees with $s$, so every zero mapping into $K_j$ lies
in $U$, where transversality holds.

Concatenate the bordisms on intervals with endpoints tending to $1$.
Over every $\operatorname{int}K_j$, the construction is unchanged
after finitely many stages. We may therefore add at $1$ the resulting
ordinary manifold $N$. 

Its map to $T$ is  clearly proper. 

This proves surjectivity.

\smallskip\noindent\emph{Step 3.}
If two ordinary manifolds are properly derived bordant, apply the
same construction to the interior of the bordism, leaving smaller
boundary collars unchanged. The resulting ordinary proper bordism
has the original boundary. This proves injectivity.
\end{proof}

The fundamental class homomorphism consequently defines
\begin{equation}\label{e:proper-vfc-map}
 \Pi_T^{\mathrm{prop}}:\mathrm dB_k^{\mathrm{prop}}(T)
 \longrightarrow H_k^{\mathrm{BM}}(T;\ZZ),\qquad
 [N,q]\longmapsto q_*[N],
\end{equation}
where $N$ is ordinary. For an oriented derived smooth manifold $M$
whose underlying space $X$ is a locally compact second countable
ANR, define
\begin{equation}\label{e:bm-vfc}
 [M]^{\mathrm{vir}}_{\mathrm{BM}}
 :=\Pi_X^{\mathrm{prop}}([M,\id_X]).
\end{equation}
This class is preserved by oriented equivalences. If
$\partial W\simeq-M_0\amalg M_1$ and $q:|W|\to T$ is a proper
continuous map constant along the boundary collars, then
\begin{equation}\label{e:proper-bordism-invariance}
 (q|_{|M_0|})_*[M_0]^{\mathrm{vir}}_{\mathrm{BM}}
 =(q|_{|M_1|})_*[M_1]^{\mathrm{vir}}_{\mathrm{BM}},
\end{equation}
provided the underlying spaces of $M_0,M_1$ are ANRs.
Indeed, the two boundary pairs represent the same proper bordism
class over $T$.

When $X$ is compact, every proper bordism over $X$ is compact,
so \eqref{e:bm-vfc} agrees with Theorem~\ref{r:3.15}.

\subsection{Base change}\label{ss:noncompact-relative}

For a continuous map $f:X\to B$ between locally compact Hausdorff spaces,
\textit{Pardon's relative homology} group is
\begin{equation}\label{e:relative-bm}
 H_k^{\mathrm{rel}\infty}(X/B;\ZZ)
 :=H^{-k}(X,f^!\ZZ_B),
\end{equation}
where $f^!$ denotes the exceptional inverse image of sheaves
\cite[Definition~2.2]{PardonCurves}.
If $B$ is a complex manifold of dimension $b$, its complex
orientation identifies $f^!\ZZ_B$ with $\omega_X[-2b]$, and hence
\begin{equation}\label{e:relative-bm-identification}
 H_k^{\mathrm{rel}\infty}(X/B;\ZZ)
 \cong H_{k+2b}^{\mathrm{BM}}(X;\ZZ).
\end{equation}

Let $g:S\to B$ be holomorphic between complex manifolds, and put
$X_S=X\times_B S$. Under \eqref{e:relative-bm-identification},
the base-change map of \cite[Definition~2.3]{PardonCurves} is
\begin{equation}\label{e:relative-bm-pullback}
 g^*\alpha=\Gamma_g^!(\alpha\times[S]).
\end{equation}
Here $\Gamma_g=\{(g(s),s):s\in S\}\subseteq B\times S$
is the graph of $g$. Its normal bundle is $g^*TB$, with
its complex orientation. The refined Gysin map
\[
 \Gamma_g^!:H_\ell^{\mathrm{BM}}(X\times S;\ZZ)
 \longrightarrow H_{\ell-2b}^{\mathrm{BM}}(X_S;\ZZ)
\]
is defined by pulling back the Thom class of $\Gamma_g$
along $f\times\id_S:X\times S\to B\times S$ and taking
cap product with the resulting cohomology class with
support on $X_S$.
Its degree is $-2b$.

\begin{prop}\label{prop:bm-base-change}
Let $M\to B_{\RR}$ be an oriented derived smooth manifold of virtual
dimension $k$. Suppose $X=|M|$ and $X_S$ are locally compact second
countable ANRs. Then
\begin{equation}\label{e:bm-base-change}
 \Gamma_g^!\bigl([M]^{\mathrm{vir}}_{\mathrm{BM}}\times[S]\bigr)
 =[M\times_{B_{\RR}}S_{\RR}]^{\mathrm{vir}}_{\mathrm{BM}},
\end{equation}
where the fibre product has the orientation of Section~\ref{ss:3.5}.
\end{prop}

\begin{proof}
\emph{Step 1.}
By Lemma~\ref{lem:proper-bordism}, choose a proper derived bordism
$(W,Q)$ over $X$ from $(M,\id_X)$ to an ordinary manifold $(N,q)$.
Write $f:X\to B$ for the underlying map.
Since $X_S$ is an ANR and is closed in $X\times S$, there is a
retraction $r:V\to X_S$ from an open neighbourhood $V$ in
$X\times S$. Choose a proper metric on $X\times S$ and make $V$
smaller so that $d(z,r(z))<1$ for $z\in V$.
Put
\[
 P=Q\times\id_S:|W|\times S\to X\times S,\qquad
 U=P^{-1}(V),\qquad R=r\circ P:U\to X_S.
\]
Therefore
$R$ is proper on every closed subset of $|W|\times S$ contained
in $U$.

\smallskip\noindent\emph{Step 2.}
By smooth approximation, there is a homotopy
$F_t=(h_t,\id_S)$ from $(f\circ Q,\id_S)$ to a smooth
morphism $(h,\id_S)$, fixed near $M\times S$ and constant
in the collar coordinate at both boundaries.
The transversality theorem on the smooth manifold $N\times S$
allows us to choose $h$ so that $(h,\id_S)|_{N\times S}$
is transverse to $\Gamma_g$.
Since $\Gamma_g$ is closed and $F_0^{-1}(\Gamma_g)\subset U$,
the approximation and homotopy may be chosen sufficiently
close to $F_0$ at each point that
$F_t^{-1}(\Gamma_g)\subset U$ for all $t$.

The  inverse image $W'$ of $\Gamma_g$ under
$(h,\id_S)$ is a bordism from
$M\times_{B_{\RR}}S_{\RR}$ to the smooth manifold
\[
 N'=\{(n,s)\in N\times S:h(n,s)=g(s)\}.
\]
Its underlying space is closed in $|W|\times S$ and contained
in $U$. By Step~1, the map $R|_{|W'|}$ is therefore proper, and
its restriction to the underlying space $X_S$ of
$M\times_{B_{\RR}}S_{\RR}$ is the identity.
Thus \eqref{e:proper-bordism-invariance} implies
\[
 [M\times_{B_{\RR}}S_{\RR}]^{\mathrm{vir}}_{\mathrm{BM}}
 =(R|_{N'})_*[N'].
\]

\smallskip\noindent\emph{Step 3.}
Let $H:N\times S\times[0,1]\to B\times S$ be the homotopy
of Step~2, and put
\[
 Z=H^{-1}(\Gamma_g),\qquad
 Z_i=Z\cap(N\times S\times\{i\}),\quad i=0,1.
\]
Write $\tau$ for the Thom class of $\Gamma_g$.
The cap product of $H^*\tau$ with the relative fundamental
class of $N\times S\times[0,1]$ is a relative Borel--Moore
class for $(Z,Z_0\amalg Z_1)$. Its boundary is
\[
 H_1^*\tau\frown[N\times S]
 -
 H_0^*\tau\frown[N\times S].
\]
Since $Z$ is closed and contained in $U\times[0,1]$,
Step~1 and compactness of $[0,1]$ imply that
$(n,s,t)\mapsto R(n,s)$ defines a proper map $Z\to X_S$.
The two terms therefore have equal pushforwards to $X_S$.

At $t=0$, we have $R(n,s)=(q(n),s)$ on $Z_0$, so the
pushforward is $\Gamma_g^!(q_*[N]\times[S])$.
At $t=1$, transversality identifies the cap product with
$[N']$, with the orientation used in Step~2.
Consequently,
\[
 \Gamma_g^!(q_*[N]\times[S])=(R|_{N'})_*[N'].
\]
Together with Step~2 and
$q_*[N]=[M]^{\mathrm{vir}}_{\mathrm{BM}}$, this proves
the proposition.
\end{proof}

\subsection{Relative virtual classes}\label{ss:noncompact-classes}

\begin{thm}\label{thm:noncompact-relative-vfc}
Let $B$ be a complex manifold, and let $f:\bX\to B$ be a derived
complex analytic space, with
$X=t_0(\bX)$ Hausdorff and second countable. Suppose that $\bX/B$
has an oriented locally exact relative $(-2)$-shifted symplectic
form $(\omega,o)$ and constant relative virtual dimension
$d=\vdim_{\C}(\bX/B)$.
There is a canonical class
\begin{equation}\label{e:noncompact-relative-vfc}
 [\bX/B,\omega,o]^{\mathrm{vir}}\in
 H_d^{\mathrm{rel}\infty}(X/B;\ZZ).
\end{equation}
It is independent of the choices in Proposition~\ref{r:6.7},
restricts to the corresponding class on open subsets of $X$,
and is compatible with holomorphic base change between complex
manifold bases. Its pullback to $u\in B$ is
\[
 [\bX_u,\omega_u,o_u]^{\mathrm{vir}}_{\mathrm{BM}}
 \in H_d^{\mathrm{BM}}(X_u;\ZZ),
\]
which agrees with \eqref{e:7.6} when $X_u$ is compact.
If $u_0,u_1$ lie in the same connected component of $B$, then
\begin{equation}\label{e:noncompact-deformation}
 (i_{u_0})_*[\bX_{u_0},\omega_{u_0},o_{u_0}]^{\mathrm{vir}}_{\mathrm{BM}}
 =(i_{u_1})_*[\bX_{u_1},\omega_{u_1},o_{u_1}]^{\mathrm{vir}}_{\mathrm{BM}}
 \quad\text{in }H_d^{\mathrm{BM}}(X;\ZZ),
\end{equation}
where $i_u:X_u\hookrightarrow X$ is the inclusion.
\end{thm}

\begin{proof}
Let $M\to B_{\RR}$ be the oriented derived smooth manifold of
Proposition~\ref{r:6.7}. Define \eqref{e:noncompact-relative-vfc}
to be the class corresponding to $[M]^{\mathrm{vir}}_{\mathrm{BM}}$
under \eqref{e:relative-bm-identification}.
For two constructions of $M$, Proposition~\ref{r:7.1} provides an oriented
bordism whose underlying map to $X$ is the proper projection
$X\times[0,1]\to X$. Equation~\eqref{e:proper-bordism-invariance}
therefore proves independence.

For a holomorphic map $g:S\to B$, put
$M_S=M\times_{B_{\RR}}S_{\RR}$.
By Proposition~\ref{r:7.2}, $M_S$ is obtained from the pulled-back
cell presentations and admissible subbundles by the construction
of Proposition~\ref{r:6.5}.
Proposition~\ref{r:7.1}, applied over $S$, therefore compares it
with any other construction on $\bX\times_B^{\R}S$ by an oriented
bordism with underlying space $X_S\times[0,1]$.
Its projection to $X_S$ is proper.
Proposition~\ref{prop:bm-base-change} and
\eqref{e:proper-bordism-invariance} now prove compatibility with
base change. The same comparison for atlases on an open subset
of $X$, together with restriction in \eqref{e:bm-vfc}, proves
compatibility with restriction.
Taking $S=\{u\}$ proves the assertion about fibres.

Finally, choose a smooth path $\gamma:[0,1]\to B$ from $u_0$ to
$u_1$, constant near its endpoints. By Proposition~\ref{r:7.2},
the pullback $W_\gamma=M\times_{B_{\RR}}[0,1]$ is an oriented
derived bordism from $M_{u_0}$ to $M_{u_1}$.
Its underlying projection $X\times_B[0,1]\to X$ is proper:
$X\times_B[0,1]$ is closed in $X\times[0,1]$.
Equation~\eqref{e:proper-bordism-invariance} now proves
\eqref{e:noncompact-deformation}.
\end{proof}

\begin{cor}\label{cor:noncompact-pt}
In the situation of Theorem~\ref{thm:pt-classes}, the derived space
of stable pairs has a canonical relative virtual class
\begin{equation}\label{e:noncompact-pt}
 [\mathbf P/B]^{\mathrm{vir}}\in
 H_{2n}^{\mathrm{rel}\infty}(P/B;\ZZ).
\end{equation}

The class is compatible with holomorphic base change and with
restriction to open subsets of $P$. On a compact fibre it agrees
with \eqref{e:pt-vfc}.
\end{cor}

\begin{proof}
The hypotheses of Theorem~\ref{thm:noncompact-relative-vfc}, with
$d=2n$, were verified in the proof of Theorem~\ref{thm:pt-classes}.
\end{proof}

\appendix
\section{Analytic moduli spaces of stable pairs}\label{app:pt-moduli}

Let $\pi:Y\to B$ be a holomorphic submersion of relative dimension
four between Hausdorff second countable complex manifolds.
We use the definition of families of stable pairs in
Section~\ref{ss:pt-moduli}, but this time do not require Calabi--Yau structure.

\subsection{Analytic representability}\label{app:pt-representability-section}

\begin{thm}\label{app:pt-representability}
The functor of isomorphism classes of families of stable pairs on
$Y/B$ with proper support is represented by a complex analytic
space $P(Y/B)$, separated over $B$. 

This space is Hausdorff and
second countable and carries a universal stable pair.

For every holomorphic map $B'\to B$, the space
$P(Y/B)\times_B B'$ represents the corresponding functor for
$Y\times_B B'\to B'$.
\end{thm}

\begingroup\color{black}
\begin{proof}
\emph{Step 0.}
We apply the analytic representability criterion of
\cite[Theorem~4.3]{FlennerLubke}. For the functor in the theorem,
we must check four conditions: every stable pair admits a
semiuniversal deformation, formal versality is open in every
family, for two families over $T$, the functor of maps $T'\to T$
over which they become isomorphic is represented by a closed
analytic subspace of $T$, and isomorphism classes form a sheaf
for open covers.

\smallskip\noindent
\emph{Step 1.}

Coherent sheaves with proper support admit convergent
semiuniversal deformations by
\cite[Theorem~9.1]{BingenerRepresentability}, and versality is
open by \cite[Theorem~8.1]{BingenerVersality}.

We remind now that for a coherent $\cO_T$-module $G$, the \textit{linear fibre space} $V(G)$
is characterized by
\[
 \operatorname{Hom}_T(T',V(G))
 =\operatorname{Hom}_{\cO_{T'}}(g^*G,\cO_{T'})
\]
for every $g:T'\to T$.

Recall also that for a flat family $F$ over $T$, its sections are represented by
a linear fibre space: apply \cite[Theorem~2.3]{FlennerLubke}
to $\cO_Z,F$ on the proper space
$Z=V(\operatorname{Ann}F)\to T$.

To pass from sheaves to sheaves with a section, we apply
\cite[Corollary~10.3]{BingenerRepresentability} to the
functor which forgets the section. Besides the deformation
results above, this requires an infinitesimal gluing
condition, denoted by (S1). We verify it as follows.
Let $A_1\to A_0\leftarrow A_2$ be homomorphisms of local
Artinian $\mathbb C$-algebras, with $A_1\to A_0$ surjective.
Flat families of sheaves over $A_1$ and $A_2$, together
with an isomorphism of their pullbacks to $A_0$, glue
over $A_1\times_{A_0}A_2$, uniquely up to an isomorphism
compatible with these identifications. Sections agreeing
over $A_0$ also glue, so (S1) holds for sheaves with a
section. The corollary therefore proves that these pairs
admit semiuniversal deformations and that versality is
open. By \cite[Remark~4.3(3)]{BingenerVersality}, formal
versality is equivalent to versality.
The conditions that the fibres be nonzero and pure of dimension
one are open in flat families with proper support
\cite[Remark~6.2]{Toma}. Since the cokernel of the section has
proper support, the condition that its fibres have dimension
at most zero is also open, by upper semicontinuity of fibre
dimension. Hence stable pairs admit semiuniversal deformations,
and their formally versal loci are open.

\smallskip\noindent
\emph{Step 2.}
Let $(F_i,s_i)$, $i=1,2$, be families over $T$, and put
$Z=V(\operatorname{Ann}F_1\cap\operatorname{Ann}F_2)$.
Then $Z\to T$ is proper. By
\cite[Theorem~2.3]{FlennerLubke}, the relative Hom functors
from $F_1$ and $\cO_Z$ to $F_2$ are represented by linear
fibre spaces $V(A)$ and $V(C)$ for coherent
$\cO_T$-modules $A,C$.
Precomposition with $s_1$ corresponds to a map $C\to A$.
Its dual on each fibre is injective, since
\[
 \operatorname{Hom}_{Y_t}
 (\operatorname{coker}s_{1,t},F_{2,t})=0
\]
by purity. Thus $C\to A$ is surjective by Nakayama's lemma,
and $V(A)\to V(C)$ is a closed embedding. Pulling back along
the section defined by $s_2$, we obtain a closed subspace
$D_{12}\subseteq T$ representing morphisms of pairs.
In particular, after any base change there is at most one such
morphism, and every endomorphism of a pair is the identity.
Construct $D_{21}$ in the same way. Then
$D_{12}\times_TD_{21}$ represents isomorphisms, since the
compositions of the two universal morphisms are identities.

\smallskip\noindent
\emph{Step 3.}
Second countability is proved in Proposition~\ref{app:pt-countability}
below, and the rest is clear.
\end{proof}
\par\endgroup

\subsection{Second countability}\label{app:pt-countability-section}

\begin{prop}\label{app:pt-countability}
The representing space $P(Y/B)$ is second countable.
\end{prop}

\begin{proof}
Let $D$ be the open subspace of the relative Douady space
parametrizing nonempty compact Cohen--Macaulay curves of pure
dimension one in fibres of $Y/B$. Its universal curve
$\mathcal C\subseteq Y_D$ is proper and flat over $D$.
Fujiki's countability theorem and its relative version
\cite[pp.~431, 447]{Fujiki} show that $D$, and hence
$\mathcal C$, is second countable.

Now, we recall that for a smooth
relative fourfold $X\to T$ and a $T$-flat sheaf $M$ with
Cohen--Macaulay fibres of codimension $c$, where $c=3$ or $4$,
there are local free resolutions of length $c$. Their duals
have fibre cohomology only in degree $c$. The local flatness
criterion \cite[Tags 00MH, 00MI]{Stacks} shows that
\[
 D_c(M)=\mathcal Ext^c_X(M,\Omega^4_{X/T})
\]
is $T$-flat and commutes with arbitrary base change. 

Put $\omega_{\mathcal C/D}=D_3(\cO_{\mathcal C})$.

By Pourcin's relative Douady theorem \cite{Pourcin}, quotients
of $\omega_{\mathcal C/D}$ flat over the parameter space with
proper support are represented by a relative Quot space.
Let $Q$ be its open subspace of quotients with zero-dimensional
fibres, including the zero quotient. (Note that Fujiki's theorem makes
$Q$ second countable.) Write
\[
 0\longrightarrow K\longrightarrow\omega_{\mathcal C_Q/Q}
 \longrightarrow G\longrightarrow0
\]
for its universal quotient. The kernel $K$ is flat over $Q$ with
Cohen--Macaulay one-dimensional fibres. Duality gives the exact sequence
\[
 0\longrightarrow\cO_{\mathcal C_Q}\longrightarrow
 F:=D_3(K)\longrightarrow D_4(G)\longrightarrow0.
\]
Thus $\cO_{Y_Q}\to F$ is a stable-pair family, defining a map
$Q\to P(Y/B)$.

To prove surjectivity on points, let $(F_0,s_0)$ be a
stable pair, and let $C_0$ be defined by $\ker s_0$.
The subsheaf $(\ker s_0)F_0\subseteq F_0$ vanishes by
purity, since it is supported on the finite set
$\operatorname{Supp}(\operatorname{coker}s_0)$.
Thus $F_0$ is an $\cO_{C_0}$-module. Both $F_0$ and its
subsheaf $\cO_{C_0}$ are pure of dimension one, and hence
Cohen--Macaulay. Dualizing
$0\to\cO_{C_0}\to F_0\to\operatorname{coker}s_0\to0$
over $C_0$, we therefore obtain a quotient
\[
 \omega_{C_0}\twoheadrightarrow
 \mathcal Ext^1_{C_0}
 (\operatorname{coker}s_0,\omega_{C_0}).
\]
By biduality, the preceding construction applied to this
quotient recovers $(F_0,s_0)$, proving surjectivity.

Since $Q$ is locally compact and second countable, it is
a countable union of compact subsets, whose images cover
$P(Y/B)$. Covering each image by finitely many coordinate
neighbourhoods, we obtain a countable open cover of
$P(Y/B)$ by second countable spaces. The union of
countable bases for these neighbourhoods is therefore
a countable basis for $P(Y/B)$.-
\end{proof}

\subsection{Comparison with perfect complexes}\label{app:pt-complexes}

\begin{prop}\label{app:pt-comparison}
Let $q:X\to T$ be a smooth morphism of complex analytic spaces
of relative dimension four. Then
\[
 (\cO_X\to F)\longmapsto I=[\cO_X\to F]
\]
defines an equivalence between the groupoid of stable-pair
families with proper support and the groupoid of perfect
complexes $I$ with a trivialization $\eta:\det I\cong\cO_X$
satisfying the following conditions:
\begin{enumerate}
\renewcommand{\labelenumi}{(\roman{enumi})}
\item the derived restriction $I_t$ to each fibre is a
stable-pair complex, with terms in degrees zero and one;
\item the closed subset $Z(I)\subseteq X$ where $I$ is not
isomorphic to a line bundle in degree zero is proper over $T$.
\end{enumerate}

\end{prop}

\begingroup\color{black}
\begin{proof}
For a fixed smooth projective Calabi--Yau fourfold and algebraic
bases, the assertion follows from Cao and Toda
\cite[Theorem~1.4]{CaoToda} by restricting to stable pairs.
See also
\cite[Proposition~1.21 and Theorem~2.7]{PT}.
We give a local proof for analytic families.

\smallskip\noindent
\emph{Step 1.}
Let $Z\subseteq X$ be a closed analytic subset whose fibres
have dimension at most one. At $x\in Z$, put
$A=\cO_{X,x}$, $R=\cO_{T,q(x)}$, and let $J\subseteq A$
be the ideal of $Z$. It contains three elements forming a
regular sequence modulo $\mathfrak m_R A$. Since $A$ is
$R$-flat, they form an $A$-regular sequence
\cite[Tag 00MG]{Stacks}, so $H_J^i(A)=0$ for $i<3$.
For a perfect $A$-complex $L$ of Tor-amplitude $[-2,1]$,
\[
 \R\Gamma_J(L^\vee)
 \simeq L^\vee\Lotimes\R\Gamma_J(A)\in D^{\geq2}(A).
\]
Thus, writing $U_A=\operatorname{Spec}A\setminus V(J)$,
restriction induces
\begin{equation}\label{app:pt-extension}
 \operatorname{Hom}_{D(A)}(L,A)\xrightarrow{\sim}
 \operatorname{Hom}_{D(U_A)}(L|_{U_A},\cO_{U_A}).
\end{equation}
If $L|_{U_A}$ is a line bundle in degree zero, the corresponding
negative Hom groups vanish.

\smallskip\noindent
\emph{Step 2.}
Let $(\cO_X\to F)$ be a family of stable pairs. Locally on $X$,
choose an exact sequence
\[
 0\longrightarrow K\longrightarrow E_2\longrightarrow E_1
 \longrightarrow E_0\longrightarrow F\longrightarrow0,
\]
where the $E_i$ are finite free sheaves. Since $F$ and the $E_i$
are flat over $T$, all the kernels in this sequence are flat
over $T$, and the sequence remains exact on each fibre.
Each $F_t$ is Cohen--Macaulay of codimension three in the
smooth fourfold $X_t$, so it admits local free resolutions
of length three. Hence $K_t$ is locally free, and
\cite[Tag~00MH]{Stacks} implies that $K$ is locally free.
Thus $F$ is perfect and $I=[\cO_X\to F]$ has
Tor-amplitude $[-2,1]$.

On $X\setminus\operatorname{Supp}F$, the isomorphism
$I\cong\cO_X$ trivialises $\det I$. By Step~1, this
trivialization and its inverse extend uniquely to $X$.
With this trivialization, $I$ satisfies (i)--(ii).

\smallskip\noindent
\emph{Step 3.}
Conversely, let $(I,\eta)$ satisfy (i)--(ii), and put $Z=Z(I)$.
The criterion for Tor-amplitude in terms of residue fields
shows that $I$ has Tor-amplitude $[-2,1]$ and that
$Z\cap X_t$ is the support of the corresponding stable pair.

To see that $Z$ is analytic, represent $I$ locally by a finite
complex $(E^\bullet,d)$ of finite free sheaves. The alternating
sum of their ranks is one. Thus $x\in Z$ precisely when
the complex of vector spaces $E^\bullet\otimes\C_x$ has
nonzero cohomology in some degree $j\ne0$, or equivalently,
\[
 \operatorname{rank}d^{j-1}(x)+\operatorname{rank}d^j(x)
 <\operatorname{rank}E^j.
\]
Each such condition is a finite union of loci defined by
the vanishing of minors of the differential matrices.
Hence $Z$ is a closed analytic subset.

On $X\setminus Z$, the determinant trivialization defines
$I\cong\det I\xrightarrow{\eta}\cO_X$.
By \eqref{app:pt-extension}, this map extends uniquely at
each analytic local ring. A representative of an extension
and a homotopy identifying its restriction with the original
map involve only finitely many germs, so are defined on
a common analytic neighbourhood. Uniqueness and the negative
Hom vanishing in Step~1 show that these local extensions
glue to a morphism $a:I\to\cO_X$.

On each fibre, choose an identification of $I_t$ with a
stable-pair complex respecting the determinant trivializations.
By Step~1, $a_t$ is then the canonical map from that complex
to $\cO_{X_t}$. Thus the derived fibres of $F=\cofib(a)$
are the corresponding sheaves. Applying
\cite[Tag~00MI]{Stacks} to a local finite free complex
representing $F$, we find that $F$ is a coherent sheaf
flat over $T$. Its support lies in $Z$, so $\cO_X\to F$
is a family of stable pairs with proper support.

\smallskip\noindent
\emph{Step 4.}
The triangle $I\to\cO_X\to F$ shows that the two constructions
are inverse on objects. Their determinant trivializations agree
outside the support, and therefore agree everywhere by Step~1.

If an isomorphism $f:I\to I'$ preserves the determinant trivializations,
then $a'f=a$ by Step~1. Taking cofibres, we get an isomorphism
$F\to F'$ preserving the sections. It is independent of
the homotopy between $a'f$ and $a$, since
$\operatorname{Hom}(I,\cO_X[-1])=0$ by Step~1.
Conversely, an isomorphism of pairs induces an isomorphism
of their complexes preserving the determinant trivializations.
As $\operatorname{Hom}(\cO_X,F'[-1])=0$, these operations
are inverse on isomorphisms as well.

The rest is clear.
\end{proof}
\par\endgroup

\end{document}